\documentclass[preprint,12pt,authoryear]{elsarticle}

\usepackage{amssymb}
\usepackage{amsmath}
\usepackage{graphicx}
\usepackage{caption} 
\usepackage{subcaption} 
\usepackage{booktabs} 
\usepackage{longtable}
\usepackage{setspace}
\usepackage{color}
\usepackage{colortbl} 
\usepackage{tikz}
\usetikzlibrary{positioning,arrows.meta,shapes.geometric}
\usepackage{pdflscape}
\usepackage{tabularx}
\usepackage{adjustbox}

\usepackage[ruled,vlined]{algorithm2e}

\usepackage{siunitx}
\journal{Computers \& Operations Research}

\begin{document}

\begin{frontmatter}



\title{A priority-driven constructive heuristic for assigning and scheduling spontaneous volunteers in disaster response} 


\author{Martina Sperling} 

\affiliation{organization={Paderborn University},
            addressline={Warburger Strasse 100}, 
            city={Paderborn},
            postcode={33100}, 
            state={North Rhine-Westphalia},
            country={Germany}}

\begin{abstract}
    Large-scale disaster response operations frequently involve spontaneous volunteers who arrive independently at disaster sites and must be coordinated under severe time pressure. 
    Assigning such volunteers to relief activities constitutes a complex workforce assignment and scheduling problem with heterogeneous capabilities, dynamic arrivals, and operational constraints.    
    Recent work formulated the spontaneous volunteer coordination problem (SVCP) as a lexicographic multi-objective mixed-integer optimization model. However, solving this model to optimality becomes computationally challenging in large-scale and rolling-horizon disaster response settings.    
    This paper proposes a problem-specific constructive heuristic for the SVCP that explicitly leverages the lexicographic objective hierarchy, capability scarcity among volunteers, and workload balancing across activities.    
    A large-scale computational study based on empirically grounded disaster response scenarios derived from the 2013 flood response in Halle (Germany) evaluates the proposed approach. Across 3\,200 simulated instances with up to 10\,000 volunteers and more than 4\,000 activity–time combinations, the heuristic closely approximates optimal solutions for the primary objectives while achieving a median runtime speedup of approximately 28×. Whereas the exact solver exceeds operational decision time limits in more than 60\% of instances, the heuristic consistently produces solutions within minutes, enabling real-time decision support for spontaneous volunteer coordination.
\end{abstract}

\begin{graphicalabstract}
\end{graphicalabstract}

\begin{highlights}
    \item Priority-driven constructive heuristic for spontaneous volunteer coordination 
    \item Exploits lexicographic priorities and volunteer capability scarcity 
    \item Computational study with up to 10\,000 volunteers and 3\,200 instances 
    \item Median runtime speedup of about $28\times$ while preserving primary objectives 
\end{highlights}


\end{frontmatter}


\section{Introduction}
\label{sec:intro}

    Disaster response operations require the rapid coordination of diverse resources under severe time pressure and uncertainty. 
    From an operations research (OR) perspective, these coordination challenges give rise to complex resource allocation and workforce scheduling problems. 
    Effective allocation of personnel, equipment, and relief tasks is essential to ensure timely assistance to affected populations.  
    Within OR, these challenges are addressed in the field of disaster operations management, which develops analytical models and decision-support methods for planning and coordinating disaster response activities. 
    Early surveys such as \cite{Altay2006} and \cite{Galindo2013} highlight the potential of OR methods to support disaster response and emergency management, while more recent reviews emphasize the growing role of optimization, simulation, and heuristic approaches in humanitarian logistics and disaster response planning \citep{Farahani2020,Gutjahr2016,Alturki2024}. 
    Disaster response has been described as the “management of disorganization” \citep[p.~126]{Simpson2009}, reflecting the complex and dynamic decision environments faced by emergency managers.
    From an OR perspective, these coordination challenges naturally give rise to large-scale assignment and scheduling problems that require efficient optimization and heuristic solution approaches.
    
    A central operational challenge in disaster response concerns the allocation and scheduling of human resources. 
    Within OR, disaster response coordination problems are often modeled as personnel assignment and workforce scheduling problems that integrate task allocation, skill compatibility, and temporal scheduling decisions.  
    However, most existing models assume that responders are organized personnel affiliated with professional emergency services or humanitarian organizations.    

    In large-scale disaster events, response operations often attract large numbers of spontaneous volunteers (SVs) who arrive independently at disaster sites without prior coordination or formal affiliation.
    These volunteers can provide valuable additional manpower but also create substantial coordination challenges. 
    Unlike affiliated responders, SVs are typically not pre-registered in organizational systems, and their arrival times, availability, and skills are uncertain \citep{Fernandez2006,Whittaker2015}. 
    Coordinating such volunteers therefore represents a particularly challenging workforce assignment and scheduling problem within disaster operations management.
    
    The role of SVs has received increasing attention in the broader disaster management literature. 
    Empirical studies highlight both the benefits and the coordination challenges associated with SV participation in disaster response and recovery activities \citep{Whittaker2015,Daddoust2021}. 
    Despite this growing recognition, decision-support models for coordinating SVs remain comparatively scarce in the OR literature \citep{Hager2024}. 
    Recent research has also investigated digital platforms for managing SVs \citep{Betke2024HICSS}, but these studies mainly address information system design rather than optimization-based decision support for operational planning.
    
    Within OR, research on volunteer coordination has primarily focused on assignment and scheduling models for organized volunteers. 
    Early work on volunteer labor allocation already emphasized structural differences between volunteer workforce management and traditional labor scheduling problems \citep{Sampson2006}. 
    Subsequent studies developed optimization models that assign volunteers to tasks while balancing objectives such as task coverage, volunteer preferences, and operational efficiency \citep{Garcia2018}. 
    However, most existing approaches assume relatively stable workforce structures and do not explicitly address the dynamic and uncertain arrival patterns typical for SVs.

    To address this gap, \cite{Sperling2022} introduced a mixed-integer optimization model for coordinating SVs in disaster response operations. 
    The model integrates assignment and scheduling decisions and evaluates solutions using lexicographically ordered objective functions that reflect operational priorities in disaster management. 
    While this formulation captures key operational characteristics of SV coordination, solving the resulting optimization problem to optimality can become computationally demanding for realistically sized and dynamically evolving instances.

    In real disaster response operations, coordination decisions must be updated repeatedly as new information about volunteer arrivals, task demands, and operational conditions becomes available. 
    Decision cycles may occur at short intervals (e.g., every 30 minutes), requiring solution approaches capable of generating high-quality decisions within strict computational time limits. 
    In such rolling-horizon decision environments, the computational effort and solving includes the repeated construction of large mixed-integer programs for each planning instance, which can become a substantial computational overhead in large-scale coordination settings.

    Heuristic solution approaches therefore play an important role in enabling real-time decision support in complex disaster response environments. 
    However, while several optimization models for volunteer coordination have been proposed, existing work primarily focuses on exact optimization approaches.
    To the best of our knowledge, existing studies on the spontaneous volunteer coordination problem (SVCP) have primarily focused on exact optimization approaches, most notably the lexicographic multi-objective model proposed by~\cite{Sperling2022}, while heuristic solution methods remain largely unexplored.

    This paper addresses this gap by proposing a problem-specific constructive heuristic for the SVCP. 
    In particular, the heuristic exploits the lexicographic priority hierarchy, capability scarcity among volunteers, and workload balancing requirements across activities. 
    These characteristics motivate a priority-driven constructive heuristic that efficiently seeks to satisfy the lexicographic objective hierarchy while achieving the computational efficiency required for real-time decision support.

    A comprehensive computational study evaluates the performance of the proposed heuristic relative to an exact optimization solver across a wide range of simulated disaster scenarios. 
    The experiments analyze how operational factors such as SV inflow, task dynamics, and workforce heterogeneity influence the trade-off between solution quality and computational runtime.

    The contributions of this paper are threefold. 
    First, we exploit structural properties of the SVCP to design a priority-driven constructive heuristic for large-scale SV assignment and scheduling.
    Second, we conduct a systematic computational study comparing the heuristic with an exact optimization under lexicographic objective evaluation. 
    Third, we provide computational insights into the conditions under which heuristic approaches provide an effective alternative to exact optimization in time-critical disaster response environments.

    The remainder of the paper is organized as follows. 
    Section~\ref{sec:related_work} reviews the literature on volunteer coordination and decision-support models in disaster response operations. 
    Section~\ref{sec:problem_description_requirements} briefly summarizes the SVCP introduced in \cite{Sperling2022}. 
    Section~\ref{sec:model} presents the optimization model used in this study. 
    Section~\ref{sec:heuristic} introduces the proposed heuristic solution approach. 
    Section~\ref{sec:computational_study} presents the computational study and discusses the experimental results. 
    Finally, Section~\ref{sec:conclusion} concludes the paper and outlines directions for future research.

\section{Related literature}
\label{sec:related_work} 

    Operations research has produced a substantial body of work on volunteer coordination and resource allocation in disaster response operations.
    A central distinction in this literature concerns the type of volunteers considered.
    Many studies focus on volunteers affiliated with humanitarian organizations who can be coordinated as a structured workforce with known availability and formal coordination by relief organizations.
    In contrast, disaster response operations frequently involve large numbers of SVs who arrive independently at disaster sites without prior organizational affiliation.
    Coordinating such volunteers poses additional challenges because their arrival times, availability, and capabilities are uncertain and difficult to predict.

    Early operations research studies primarily investigated volunteer coordination in settings with affiliated volunteers.
    For example, \cite{Falasca2009} developed multi-objective optimization models that balance unmet task demand and volunteer preferences. 
    \cite{Falasca2012} extended this line of work by proposing an optimization model for assigning volunteers to humanitarian tasks.
    Subsequent research incorporated additional operational considerations such as training requirements and workforce flexibility.
    In particular, \cite{Lassiter2014,Lassiter2015} proposed optimization and robust optimization models that account for volunteer skills and uncertainty in workforce availability.
    Similarly, \cite{Garcia2018} developed a resource allocation framework for crisis volunteer management, while \cite{Li2019model} proposed a mathematical model for allocating rescue personnel across multiple disaster areas.
    These studies typically assume that volunteers are affiliated with organizations and can therefore be coordinated within relatively structured workforce management settings.

    A smaller stream of literature addresses the coordination of SVs arriving independently at disaster sites.
    Conceptual and operational models for integrating spontaneous volunteers into emergency response systems have also been proposed \citep{Nielsen2024}. 
    In such environments, volunteer availability and task demand often evolve dynamically during disaster response operations.
    Several studies therefore investigate assignment policies under uncertainty.
    For instance, \cite{Mayorga2017} formulate the SV assignment problem as a Markov decision process that captures stochastic SV arrivals and departures.
    Building on this approach, \cite{Paret2021} analyze assignment policies under uncertainty in both SV availability and task demand.
    Related work considers the management of SV convergence at disaster relief centers.
    For example, \cite{Zayas2020} model the allocation of arriving SVs to parallel service queues using a queueing framework and analyze optimal control policies using a Markov decision process.
    Other studies develop optimization models for assigning SVs to disaster response tasks. 
    For example, \cite{Pielorz2015} propose an integer programming formulation for the geospatial allocation of volunteers.
    More recent research has also investigated optimization-based resource allocation and volunteer dispatching models in disaster response contexts \citep{Kapukaya2025,Wang2026}.

    More recent work considers integrated assignment and scheduling decisions for SVs.
    Related decision-support approaches have also been proposed for coordinating volunteers in other disaster contexts.
    For example, \cite{Krstikj2021} develop a spatial decision-support approach for allocating volunteers to neighborhoods during large-scale lockdowns based on proximity and vulnerability indicators.
    In particular, \cite{Rauchecker2018} develop a mixed-integer programming model that simultaneously assigns SVs to disaster response tasks and schedules their activities over time. 
    The model is applied sequentially as new information becomes available during disaster response operations.

    Building on this work, \cite{Sperling2022} introduce a multi-objective mixed-integer optimization model for coordinating SVs. 
    Their formulation incorporates a lexicographic objective structure reflecting operational priorities in disaster management and supports rolling-horizon coordination in which assignments are repeatedly updated as new information about volunteer arrivals, task demands, and operational conditions becomes available.

    The present paper continues this line of research by developing a heuristic solution approach for the SVCP introduced by \cite{Sperling2022}. 
    While the existing models focus on exact optimization, the proposed heuristic explicitly addresses the computational challenges arising in large-scale disaster response scenarios and enables efficient solution of the model in repeated decision cycles.

    In parallel, several studies explore approximate or heuristic approaches for volunteer allocation problems in disaster settings. 
    For example, \cite{Abualkhair2020managing} analyze volunteer assignment policies using agent-based simulation of disaster relief centers. 
    Metaheuristic approaches based on evolutionary algorithms have also been proposed to address volunteer allocation problems \citep{Rabiei2023,Xue2024}. 
    However, these approaches primarily address assignment decisions and typically do not consider integrated assignment and scheduling problems for SVs.   

    Table~\ref{tab:model_comparison} summarizes key characteristics of the models proposed in the literature with respect to volunteer type, decision scope, coordination setting, and modeling approach.
    The coordination category distinguishes between static, stochastic, and sequential coordination settings depending on whether coordination decisions are made once for a fixed instance, under explicit uncertainty, or repeatedly updated as new information becomes available.
    The table shows that most existing studies either focus on affiliated volunteers or consider assignment decisions only.
    Only a small number of studies address coordination settings involving SVs together with integrated assignment and scheduling decisions.
    Moreover, these approaches rely primarily on exact optimization methods.
    
        {\footnotesize
            \begin{longtable}{p{3.75cm} c c c p{2.1cm}}
            \caption{Characteristics of volunteer coordination models in the literature.
            A\&S denotes assignment and scheduling decisions. Sequential coordination refers to rolling-horizon planning.}
            \label{tab:model_comparison}\\
            
            \toprule
            Reference & Vol. type & Decision & Coordination & Model type \\
            \midrule
            \endfirsthead
            
            \toprule
            Reference & Vol. type & Decision & Coordination & Model type \\
            \midrule
            \endhead
            
            \midrule
            \multicolumn{5}{r}{\textit{Continued on next page}} \\
            \midrule
            \endfoot
            
            \bottomrule
            \endlastfoot
            
            \cite{Falasca2009} & Affiliated & A\&S & Static & Multi-obj. MILP \\
            \cite{Falasca2012} & Affiliated & A\&S & Static & MILP \\
            \cite{Lassiter2014} & Affiliated & Assignment & Static & MILP \\
            \cite{Lassiter2015} & Affiliated & Assignment & Static & Bi-obj. robust optimization \\
            \cite{Garcia2018} & Affiliated & Assignment & Static & MILP \\
            \cite{Li2019model} & Affiliated & Assignment & Static & MILP \\
            \cite{Rabiei2023} & Affiliated & Assignment & Static & Evolutionary metaheuristic \\
            \cite{Xue2024} & Affiliated & Assignment & Static & Evolutionary metaheuristic \\
            
            \addlinespace
            
            \cite{Mayorga2017} & Spontaneous & Assignment & Stochastic & MDP \\
            \cite{Zayas2020} & Spontaneous & Assignment & Stochastic & Queueing/MDP \\
            \cite{Paret2021} & Spontaneous & Assignment & Stochastic & MDP \\
            \cite{Pielorz2015} & Spontaneous & Assignment & Static & MILP \\
            \cite{Abualkhair2020managing} & Spontaneous & Assignment & Static & Simulation \\
            
            \addlinespace
            
            \cite{Krstikj2021} & Spontaneous & Assignment & Sequential & Optimization-based DSS \\
            \cite{Rauchecker2018} & Spontaneous & A\&S & Sequential & MILP \\
            \textbf{\cite{Sperling2022}} & Spontaneous & A\&S & Sequential & Multi-obj. MIP \\
            
            \addlinespace
            
            \rowcolor{gray!10}
            \textbf{This paper} & Spontaneous & A\&S & Sequential & Heuristic \\
            
            \end{longtable}
        }

    Table~\ref{tab:model_comparison} further illustrates that the model proposed by \cite{Sperling2022} represents one of the few approaches that combine spontaneous volunteers with integrated assignment and scheduling decisions in a sequential coordination setting. 
    Nevertheless, applying exact optimization techniques to large-scale instances of this model can require substantial computational effort.

    In practical disaster response operations, coordination decisions must be updated frequently as new information becomes available. 
    Consequently, solution approaches are required that can generate high-quality solutions within the limited computation times available in operational disaster response environments.

    To address this challenge, this paper develops a problem-specific heuristic solution approach for the SV coordination model introduced by \cite{Sperling2022}. 
    The proposed method exploits structural properties of the model, including its lexicographic priority hierarchy, capability constraints, and workload balancing requirements, to construct high-quality assignment and scheduling decisions efficiently. 
    The heuristic is specifically designed for repeated solution of large-scale instances arising in rolling-horizon disaster response coordination, where exact optimization becomes computationally demanding.  

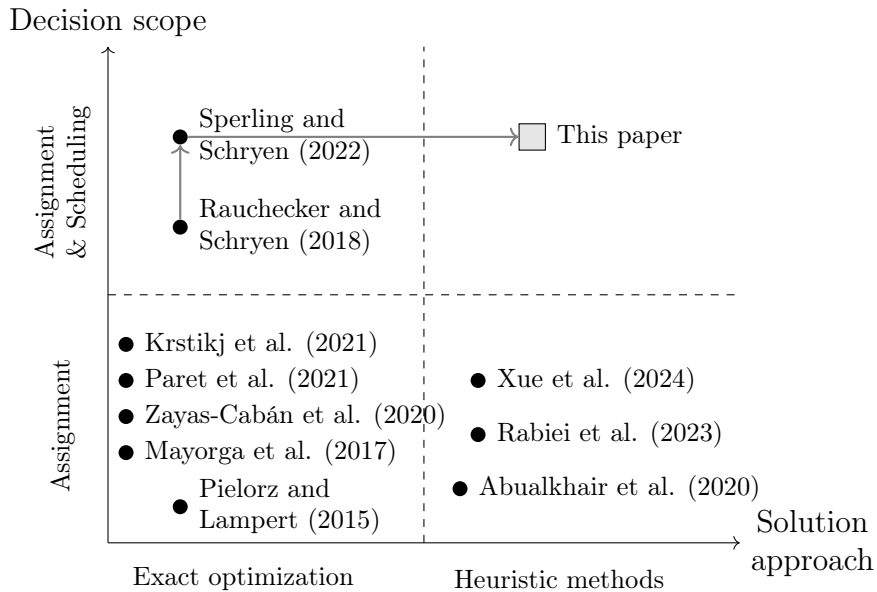
\begin{figure}[ht]
        \centering
            \begin{tikzpicture}[
            scale=1.2,
            point/.style={circle,fill=black,inner sep=2pt},
            paper/.style={font=\footnotesize,align=left}
            ]
            
            \draw[->] (0,0) -- (7,0) node[right,align=center]{Solution\\approach};
            \draw[->] (0,0) -- (0,5.5) node[above]{Decision scope};
            
            \node[paper] at (1.5,-0.4) {Exact optimization};
            \node[paper] at (5,-0.4) {Heuristic methods};
            
            \node[paper,rotate=90] at (-0.5,1.3) {Assignment};
            \node[paper,rotate=90] at (-0.5,4) {Assignment \\ \& Scheduling};
    
            \draw[dashed] (3.5,0) -- (3.5,5.5);
            \draw[dashed] (0,2.75) -- (7,2.75);
            
            \node[point] (krstikj) at (0.2,2.2) {};
            \node[paper,anchor=west] at (krstikj.east) {Krstikj et al. (2021)};
            \node[point] (paret) at (0.2,1.8) {};
            \node[paper,anchor=west] at (paret.east) {Paret et al. (2021)};
            \node[point] (zayas) at (0.2,1.4) {};
            \node[paper,anchor=west] at (zayas.east) {Zayas-Cabán et al. (2020)};            
            \node[point] (mayorga) at (0.2,1) {};
            \node[paper,anchor=west] at (mayorga.east) {Mayorga et al. (2017)};  
            \node[point] (pielorz) at (0.8,0.4) {};
            \node[paper,anchor=west] at (pielorz.east) {Pielorz and \\ Lampert (2015)}; 
            
            \node[point] (xue) at (4.1,1.8) {};
            \node[paper,anchor=west] at (xue.east) {Xue et al. (2024)};
            \node[point] (rabiei) at (4.1,1.2) {};
            \node[paper,anchor=west] at (rabiei.east) {Rabiei et al. (2023)};
            \node[point] (abualkhair) at (3.9,0.6) {};
            \node[paper,anchor=west] at (abualkhair.east) {Abualkhair et al. (2020)};

            \node[point] (rauchecker) at (0.8,3.5) {};
            \node[paper,anchor=west] at (rauchecker.east) {Rauchecker and \\ Schryen (2018)};            
            \node[point] (sperling) at (0.8,4.5) {};
            \node[paper,anchor=west] at (sperling.east) {Sperling and \\ Schryen (2022)};  
            
            \node[rectangle,draw,fill=gray!20,minimum width=0.35cm,minimum height=0.35cm] (thispaper) at (4.7,4.5) {};
            \node[paper,anchor=west] at (thispaper.east) {This paper};
            
            \draw[->,thick,gray] (rauchecker) -- (sperling);
            \draw[->,thick,gray] (sperling) -- (thispaper);
            
            \end{tikzpicture}        
        \caption{Positioning of operations research approaches for SV coordination with respect to decision scope and solution methodology.}
        \label{fig:literature_map}
        \end{figure}

    Figure~\ref{fig:literature_map} visually positions existing operations research approaches for coordinating SVs with respect to their decision scope and solution methodology.  
    Arrows indicate the methodological progression from the assignment and scheduling model of~\cite{Rauchecker2018} to the multi-objective formulation of \cite{Sperling2022}, which forms the basis for the heuristic solution approach proposed in this paper.

\section{Problem description and modeling framework}
\label{sec:problem_description_requirements}

    \subsection{Problem description}
    \label{subsec:problem_description}
        During disaster response operations, relief organizations must often coordinate large numbers of SVs who differ in their availability, capabilities, and willingness to participate.
        The resulting coordination problem consists of assigning SVs to appropriate relief activities under severe time pressure.
        The problem setting considered in this paper builds on the spontaneous volunteer coordination model introduced by \cite{Sperling2022} and focuses on the response phase of disasters.       
        At a given decision point, a disaster manager must assign available SVs to operational tasks while respecting operational constraints and organizational preferences.
        
        Tasks represent operational missions that must be performed within a specified time window at a given location. 
        Each task may consist of several \emph{task activities}, which represent specific operational actions such as filling sandbags, distributing supplies, or documenting relief efforts. 
        Each activity requires a certain number of SVs possessing a specific capability, while a single capability may enable SVs to perform multiple activities.
        
        SVs differ with respect to their capabilities and availability periods. 
        A SV may only be assigned to activities requiring a capability they possess and during time periods in which they are available. 
        Furthermore, each SV can perform at most one activity at a time.
        
        Tasks may differ in their urgency, which is expressed through priority levels defined by relief organizations. 
        When the number of available SVs is insufficient to fully staff all activities, decision makers must determine how SVs should be distributed across tasks with different priorities. 
        To support such decisions, we consider the concept of \emph{workload}, defined as the ratio between the number of SVs assigned to an activity and the number of SVs required for that activity.
        
\begin{figure}[t]
            	\centering
            	\includegraphics[scale=.3]{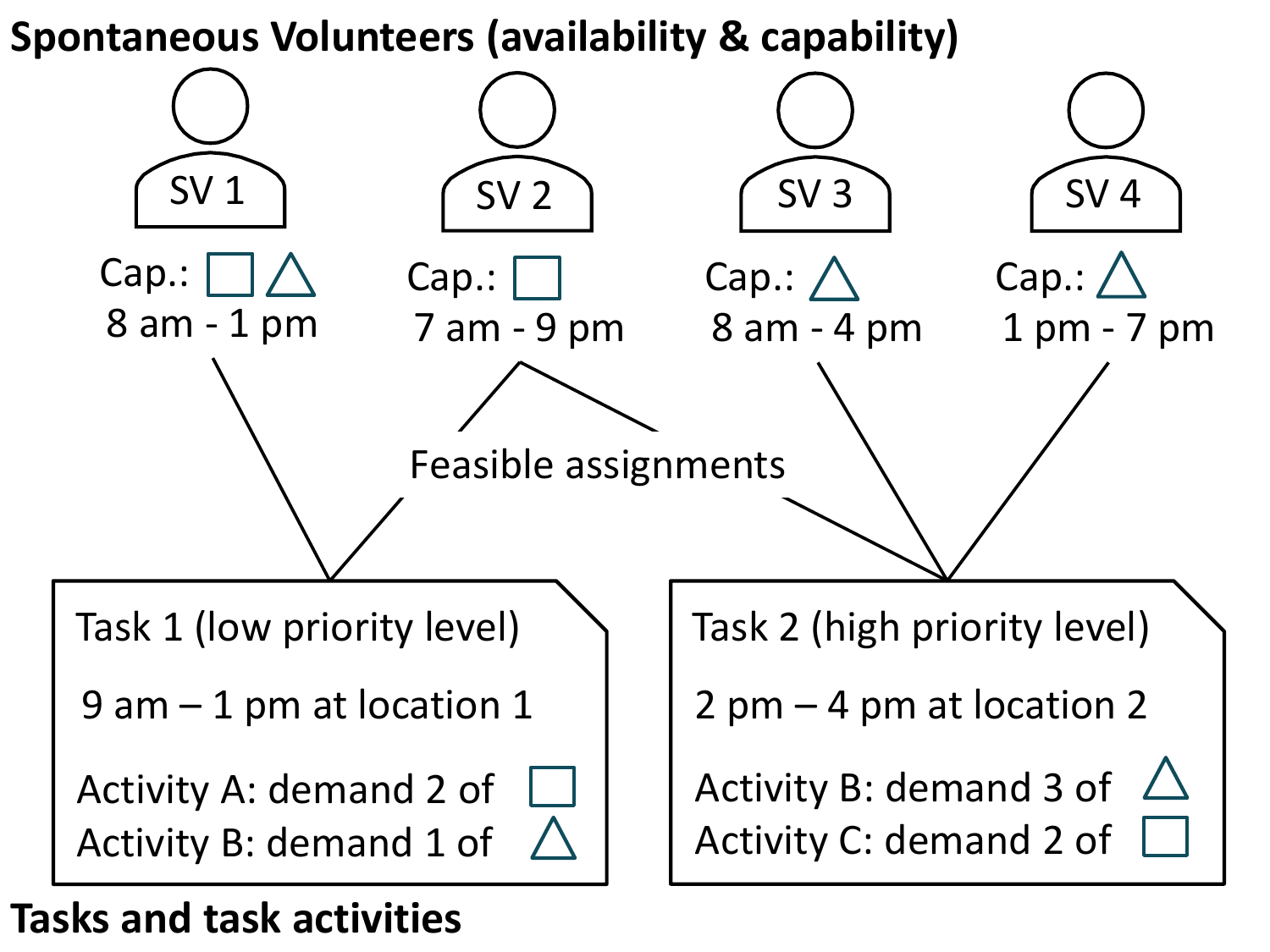}    	
            	\caption{Illustrative example of the spontaneous volunteer coordination problem, adapted from \cite{Sperling2022}.}	
            	\label{fig:sample_situation}
            \end{figure}

        Figure~\ref{fig:sample_situation} illustrates a simplified situation in which several SVs with different availability windows must be assigned to activities belonging to two tasks with different priority levels.
        
        Finally, disaster response operations evolve dynamically.
        New SVs may arrive, new tasks may emerge, and operational requirements may change over time.
        Consequently, SV coordination decisions are repeatedly made in a sequence of planning instances, where the current allocation serves as input for the subsequent decision problem.

    \subsection{Modeling requirements}
    \label{subsec:requirements}

        Based on practitioner interviews and prior work on SV coordination, a set of operational requirements must be considered when coordinating SVs. 
        These requirements capture both operational constraints and decision-making preferences of disaster managers. 
        A detailed discussion and formalization of these requirements can be found in \citet{Sperling2022}. 
        In the following, we briefly summarize the most relevant requirements that guide the design of the optimization model and the heuristic solution approach.
        
        \textbf{Task prioritization.}
        Tasks differ in their urgency and are therefore associated with priority levels determined by relief organizations. 
        In situations of SV shortage, SVs should be assigned preferentially to tasks with higher priority levels. 
        Within tasks of comparable priority, practitioners aim to maintain predefined workload ratios in order to balance relief efforts across tasks.
        
        \textbf{SV-task compatibility.}
        SVs may only be assigned to task activities if they possess the required capabilities and are available during the corresponding time periods. 
        Furthermore, each SV can work on at most one task activity at a time.
        
        \textbf{Operational constraints.}
        Assignments must respect several operational constraints. 
        These include travel or setup times between task locations, minimum working durations on task activities, and upper bounds on the total working time of SVs.
        
        \textbf{Assignment stability.}
        Since disaster response operations evolve dynamically, planning decisions are repeatedly updated over time. 
        However, previously assigned SVs are not reassigned in subsequent planning steps, resulting in a non-preemptive scheduling structure.
        
        \textbf{Workload balancing.}
        When several task activities require similar resources, relief organizations aim to balance workloads across these activities, particularly in situations of SV shortage.
        
        Based on these requirements, we formulate the optimization model in the next section and later derive a heuristic solution procedure that approximates the resulting decision structure.  

\section{SVCP model}
\label{sec:model} 

    The SVCP is formulated as a lexicographically ordered mixed-integer optimization model introduced in \cite{Sperling2022}. 
    In the following, we summarize the notation and model components relevant for the heuristic proposed in Section~\ref{sec:heuristic}.    

    \subsection{Notation}
    \label{subsec:notation}
        The model considers the assignment of SVs to task activities over a discretized planning horizon.
        Table~\ref{tab:notation} summarizes the indices, sets, parameters, and variables used throughout the model and the heuristic procedure.
        
        {\footnotesize
            \begin{longtable}{l p{0.75\linewidth}}            
                \caption{Indices, sets, parameters and variables of the optimization model and heuristic.}
                \label{tab:notation}\\
                
                \toprule
                Notation & Description/Definition \\
                \midrule
                \endfirsthead
                
                \toprule
                Notation & Description/Definition \\
                \midrule
                \endhead
                
                \midrule
                \multicolumn{2}{r}{Continued on next page} \\
                \midrule
                \endfoot
                
                \bottomrule
                \endlastfoot              
                
                \multicolumn{2}{c}{\textbf{Indices}}\\
                \cmidrule(lr){1-2}
                $t=1,\dots,T$ & Time slots\\
                $v=1,\ldots,V$ & Volunteers\\
                $c=1,\ldots,C$ & Capabilities\\ 
                $a=1,\ldots,A$ & Activities\\ 
                $p=1,\ldots,P$ & Priority levels\\
                $k=1,\ldots,K$ & Priority classes\\
                
                \addlinespace
                \multicolumn{2}{c}{\textbf{Sets}}\\
                \cmidrule(lr){1-2}
                $\mathcal{P}_{k}$ & Priority levels in priority class $k$\\
                $\hat{\mathcal{A}}_{k,t}$ & Task activities in time slot $t$ with priority levels of priority class $k$\\
                $\mathcal{A}_{p,t}$ & Task activities in time slot $t$ with priority level $p$\\  
             
                \addlinespace
                \multicolumn{2}{c}{\textbf{Parameters}}\\
                \cmidrule(lr){1-2}
                $\alpha_{k}$ & $=1$ if $|\mathcal{P}_{k}|>1$ ($0$ else)\\               
                $\sigma_{p,p+1}$ & Balancing factor for priorities $p$ and $p+1$ in the same priority class\\[2pt]
                
                $r_{a,t}$ & $=1$ if time slot $t$ is within the time frame of task activity $a$ ($0$ else)\\			
                $req_{a,c}$ & $=1$ if task activity $a$ requires capability $c$ ($0$ else) \\
                $p_{a}$ & Priority level of task activity $a$ - depends only on the underlying task\\
                $d_{a,t}$ & Demand/Supply ratio of volunteers for task activity $a$ in time slot $t$\\[2pt]
                
                $av_{v,t}$ & $=1$ if volunteer $v$ is available in time slot $t$ ($0$ else)\\			
                $cap_{v,c}$ & $=1$ if volunteer $v$ offers capability $c$ ($0$ else)\\				
                $o_{v,a,t}$ & $=1$ if volunteer $v$ has an assignment to task activity $a$ in time slot $t$ from the solution of the previous SVCP instance ($0$ else)\\[2pt]
                
                $n_a$ & Number of volunteers required for task activity $a$\\
                $n_{p,t}$ & Number of volunteers required for all task activities in $\mathcal{A}_{p,t}$ in time slot $t$\\[2pt]
                
                $\tau_{min}$ & a volunteer must work consecutively on the same task activity\\							
                $\bar{\tau}_{max}$ & Total maximum hours worked during the past $T$ slots\\
                $\tau_{v,a}$ & volunteer $v$ needs to reach task activity $a$ from his/her current position\\		
                $s_{a,a'}$ & Time required to travel from task activity $a$ to task activity $a'$\\	
                $w_t$ & Weight indicating severity of unmet demands in time slot $t$\\
                
                \addlinespace
                \multicolumn{2}{c}{\textbf{Variables}}\\
                \cmidrule(lr){1-2}
                $x_{v,a,t}$ & $=1$ if volunteer $v$ is assigned to task activity $a$ in time slot $t$ ($0$ else)\\	
                $L_{a,t}(\mathbf{x})$ & Workload of task activity $a$ in time slot $t$\\
                $\bar{L}_{p,t}(\mathbf{x})$ & Average workload of all task activities with priority level $p$ in time slot $t$\\
                $\Delta_{a,a',t}(\mathbf{x})$ & Imbalance between workloads $L_{a,t}$ and $L_{a',t}$ in time slot $t$\\ 
                $\Lambda_{p,p',t}(\mathbf{x})$ & Imbalance between average workloads $\bar{L}_{p,t}$ and $\bar{L}_{p',t}$ in time slot $t$\\    
    		\end{longtable}	
        }
       
        Priority levels are grouped into ordered priority classes $\mathcal P_k$ forming a partition of $\{1,\dots,P\}$ with $p<p'$ for all $p\in\mathcal P_k$, $p'\in\mathcal P_{k'}$, and $k<k'$, where $\hat{\mathcal A}_{k,t}$ denotes activities of priority class $k$ active in time slot $t$.
        Binary decision variables $x_{v,a,t}$ indicate whether volunteer $v$ is assigned to activity $a$ at time slot $t$.
        The workload of activity $a$ in time slot $t$ is defined as        
            \begin{equation}
            \label{eq:def-workloadtaskact}
                L_{a,t}(\mathbf{x}) = \frac{1}{n_a}\sum_{v=1}^{V} x_{v,a,t}.
            \end{equation}        
        The average workload of activities with priority level $p$ in time slot $t$ is        
            \begin{equation}
            \label{eq:def-workloadpriority}
                \bar{L}_{p,t}(\mathbf{x}) =
                \frac{1}{n_{p,t}}
                \sum_{v=1}^{V}
                \sum_{a\in\mathcal{A}_{p,t}} x_{v,a,t}.
            \end{equation}
        Imbalance between workloads is captured by variables $\Delta_{a,a',t}(\mathbf{x})$ and $\Lambda_{p,p',t}(\mathbf{x})$ defined as
            \begin{equation}
            \label{eq:def-lambda}
            \Lambda_{p,p',t}(\mathbf{x}) =
            \begin{cases}
            \left(\bar{L}_{p,t}(\mathbf{x})-\sigma_{p',p}^{-1}\bar{L}_{p',t}(\mathbf{x})\right)^+ & \text{if } p=p'+1,\\
            \left(\bar{L}_{p,t}(\mathbf{x})-\sigma_{p,p'}\bar{L}_{p',t}(\mathbf{x})\right)^+ & \text{if } p=p'-1 .
            \end{cases}
            \end{equation}        
        and        
            \begin{equation}
            \label{eq:def-delta}
            \Delta_{a,a',t}(\mathbf{x}) =
            \left(L_{a,t}(\mathbf{x})-L_{a',t}(\mathbf{x})\right)^+ .
            \end{equation}        
        where $(y)^+ := \max(0,y)$ denotes the positive part of $y$.
        To capture capability scarcity, we define the supply weight factor        
            \begin{equation}
            \label{eq:def-not-prefer-a-taskactivity}
                d_{a,t} =
                \min \left(
                1,\,
                \frac{n_a r_{a,t}}
                {req_{a,c}\sum_{v=1}^{V} av_{v,t} cap_{v,c}}
                \right).
            \end{equation}    
        The factor $d_{a,t}$ scales balancing penalties and is later used in the heuristic to prioritize scarce capabilities.        

    \subsection{Optimization model}
    \label{subsec:opt_model}
    
        The SVCP is formulated as a lexicographically ordered multi-objective optimization model following \cite{Sperling2022}. 
        The objectives reflect operational preferences reported by practitioners. 
        
    		\begin{align*}
    		    \label{OF:max-assignments-highest-priority-class}\tag{OF $1$}
    		        &\max \left(\sum_{v=1}^V\sum_{a\in\hat{\mathcal A}_{K,t}}\sum_{t=1}^T w_{t} x_{v,a,t} \right)\\
    		        &\qquad\vdots\tag{OF $\vdots$}\\
    		    \label{OF:max-assignments-lowest-priority-class}\tag{OF $K$}
    		        &\max \left(\sum_{v=1}^V\sum_{a\in\hat{\mathcal A}_{1,t}}\sum_{t=1}^T w_{t} x_{v,a,t} \right)\\
    		    \label{OF:min-imbalance-intraprio}\tag{OF $K+1$}
    		        &\min \left(\sum_{k=1}^K\alpha_{k}\sum_{\substack{ p,p'\in\mathcal{P}_{k}}}\sum_{t=1}^T \Lambda_{p,p',t}\left(\mathbf{x}\right) \right)\\
    	        \label{OF:min-imbalance-interprio}\tag{OF $K+2$}
    		        &\min \left(\sum_{p=1}^{P}\sum_{\substack{ a,a'\in\mathcal{A}_{p,t}}}\sum_{t=1}^{T}d_{a,t} d_{a',t}  \Delta_{a,a',t}\left(\mathbf{x}\right) \right)
    		\end{align*}       
            
        The objectives are evaluated in lexicographic order. 
        The first $K$ objectives maximize time-weighted assignments to higher priority classes ($w_t > w_{t+1} > 0$), favoring early task coverage. 
        Objective \eqref{OF:min-imbalance-intraprio} penalizes deviations from predefined workload ratios within priority classes, whereas objective \eqref{OF:min-imbalance-interprio} balances workloads across activities with the same priority level.
        The feasible assignment space is defined by constraints ensuring activity staffing limits, volunteer availability and capability compatibility, assignment stability across planning instances, travel and setup times, minimum assignment durations, and limits on total working time.
        Assignment variables are binary.

\section{Priority-driven constructive heuristic}
\label{sec:heuristic}

    \subsection{Heuristic design}
    \label{sec:heuristic_idea}
    
\begin{figure}[h]
                \centering
                \begin{tikzpicture}[
                node distance=2.0cm,
                box/.style={
                rectangle,
                draw,
                rounded corners,
                align=center,
                text width=0.35\linewidth,
                minimum height=1cm,
                font=\footnotesize
                },
                arrow/.style={->, thick}
                ]
                
                \node[box] (model1) {Priority objectives\\ \eqref{OF:max-assignments-highest-priority-class}--\eqref{OF:max-assignments-lowest-priority-class}};
                
                \node[box, right=2.5cm of model1] (rule1) {Priority rule\\Select highest priority class};
                
                \node[box, below of=model1] (model2) {Time weighting in priority objectives\\ \eqref{OF:max-assignments-highest-priority-class}--\eqref{OF:max-assignments-lowest-priority-class}};
                
                \node[box, right=2.5cm of model2] (rule2) {Time rule\\Select earliest time slot};
                
                \node[box, below of=model2] (model3) {Workload balancing objectives\\ \eqref{OF:min-imbalance-intraprio},\eqref{OF:min-imbalance-interprio}};
                
                \node[box, right=2.5cm of model3] (rule3) {Workload rule\\Select lowest weighted workload};
                
                \node[box, below of=model3] (model4) {Capability scarcity};
                
                \node[box, right=2.5cm of model4] (rule4) {Volunteer ordering\\Prioritize scarce capabilities};
                
                \draw[arrow] (model1) -- (rule1);
                \draw[arrow] (model2) -- (rule2);
                \draw[arrow] (model3) -- (rule3);
                \draw[arrow] (model4) -- (rule4);
                
                \end{tikzpicture}
                
            \caption{Relationship between the objective structure of the SVCP optimization model and the decision rules used in the priority-driven constructive heuristic.}
            \label{fig:heuristic_mapping}
            \end{figure}
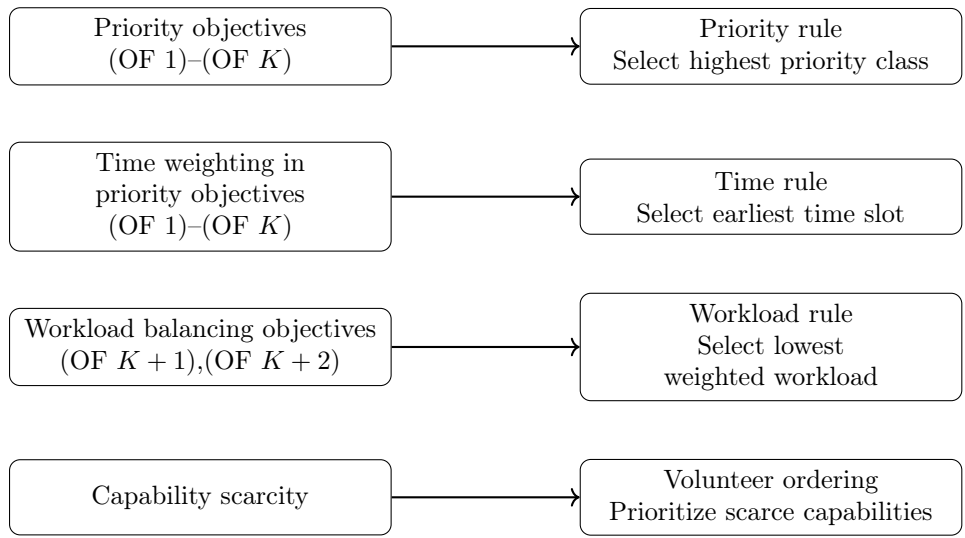

        The heuristic design exploits structural properties of the lexicographic objective hierarchy of the SVCP introduced in Section~\ref{sec:model}. 
        Figure~\ref{fig:heuristic_mapping} illustrates how the decision rules of the heuristic are derived from the objective structure of the optimization model. 
        In particular, the priority structure of tasks, the scarcity of volunteer capabilities, and the workload balancing requirements motivate a priority-driven constructive assignment procedure.

        The figure highlights that each heuristic decision rule corresponds to a specific component of the lexicographic objective hierarchy of the optimization model.
        
        \textit{Priority rule.}
        Only activity--time combinations from the highest currently active priority class are considered.
        This reflects the lexicographic structure of the priority-maximization objectives \eqref{OF:max-assignments-highest-priority-class}--\eqref{OF:max-assignments-lowest-priority-class}.
        
        \textit{Time rule.}
        Within this subset, the earliest time slot is selected in order to emphasize early task coverage, which is consistent with the time-weighted objective structure of the model.
        
        \textit{Workload rule.}
        Among activities active at the selected time slot, the activity with the lowest weighted workload is chosen. 
        If multiple activities satisfy the above criteria, ties are broken arbitrarily.        
        The weighted workload $\sigma_{p_a,p_a+1} L_{a,t^*}(\mathbf{x})$ reflects the workload-balancing structure of the optimization model and approximates objectives \eqref{OF:min-imbalance-intraprio} and \eqref{OF:min-imbalance-interprio}. 
        
        Before the constructive assignment procedure starts, volunteers are ordered according to the scarcity of their capabilities. 
        The scarcity score reflects the demand–supply ratio captured by the parameter $d_{a,t}$ and prioritizes volunteers whose capabilities are required for activities with limited supply. 
        Algorithm~\ref{alg:preprocessing} summarizes this preprocessing step.      
            
        Overall, the heuristic can be interpreted as a greedy approximation of the lexicographic objective hierarchy, where assignments are generated sequentially according to priority, temporal precedence, and workload balancing considerations.

    \subsection{Heuristic procedure}
    \label{sec:heuristic_procedure}
        Algorithm~\ref{alg:svcp-heuristic} summarizes the constructive assignment and scheduling procedure.
        During the constructive procedure, activity--time combinations are iteratively selected and feasible volunteers are assigned while respecting operational constraints. 
        The selection of activity--time combinations follows three decision rules derived from the lexicographic objective structure of the optimization model.
        
        Let
        \[
        \mathcal{A}^{\text{active}} =
        \{(a,t) \mid r_{a,t}=1,\; L_{a,t}(\mathbf{x}) < 1 \}
        \]
        denote the set of currently active activity--time combinations, i.e., task activities that are active in time slot $t$ and still require additional volunteers.
        For each priority class $k$, we define
        \[
        \mathcal{A}_{k}^{\text{active}}
        =
        \{(a,t) \in \mathcal{A}^{\text{active}} \mid p_a \in \mathcal{P}_k\},
        \]
        i.e., the subset of active combinations belonging to priority class $k$.
        At each iteration, the heuristic selects one element $(a,t) \in \mathcal{A}_{k^*}^{\text{active}}$ according to the following decision rules.      

        Priority classes are ordered such that larger values of $k$ correspond to higher priority classes.
        Formally, the selected activity--time combination $(a^*,t^*)$ satisfies
        \[
            k^* = \max\{k \mid \mathcal{A}_{k}^{\text{active}} \neq \emptyset\},
        \]
        \[
            t^* = \min\{t \mid (a,t) \in \mathcal{A}_{k^*}^{\text{active}} \},
        \]
        and
        \[
            a^* =\arg\min_{a:(a,t^*)\in\mathcal{A}_{k^*}^{\text{active}}}\sigma_{p_a,p_a+1} L_{a,t^*}(\mathbf{x}).
        \]
        
        The activity selection rule is based on the weighted workload
        \[\sigma_{p_a,p_a+1} L_{a,t^*}(\mathbf{x}),\]
        which scales the activity workload by the priority-balancing factor of the optimization model.
        This weighting reflects the desired workload ratios between adjacent priority levels within the same priority class (objective (OF K+1)).
        Within the same priority level, the rule reduces to minimizing the workload $L_{a,t^*}(\mathbf{x})$, which promotes balanced staffing across activities (objective (OF K+2)).
        
        For the selected activity--time combination $(a^*,t^*)$, feasible volunteers are identified based on capability compatibility, availability, travel times, and remaining working time. 
        Among feasible candidates, the volunteer who can start the assignment earliest is selected.
        Let $V^{cand}$ denote the set of volunteers for whom a feasible assignment interval exists.
        More precisely, the selected volunteer $v^*$ satisfies
        \[
        v^* = \arg\min_{(v,t_s,t_e)\in V^{cand}} t_s,
        \]
        where $t_s = 1+\tau_{v,a}$ and $t_e$ denote the earliest and latest feasible start and end times of the assignment interval returned by the feasibility evaluation procedure.
        This rule increases early task coverage and helps maintain flexibility for later assignments.

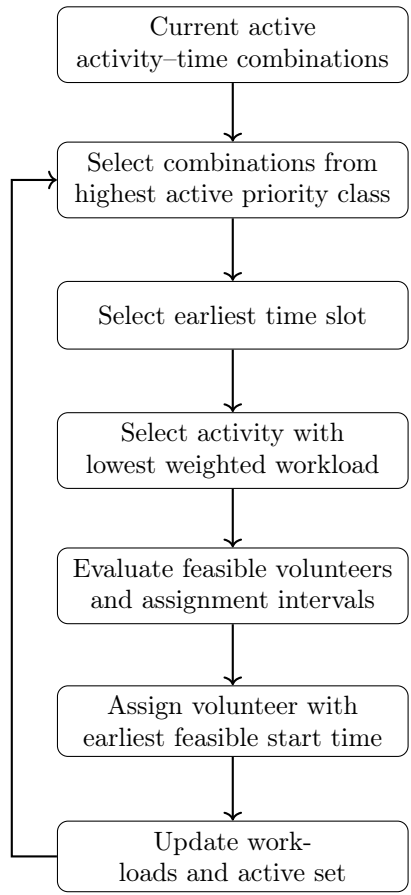
\begin{figure}
            \centering
                \begin{tikzpicture}[
                node distance=1.8cm,
                box/.style={
                rectangle,
                draw,
                rounded corners,
                align=center,
                text width=0.32\linewidth,
                minimum height=0.9cm,
                font=\footnotesize
                },
                arrow/.style={->, thick}
                ]
                
                \node[box] (active) {Current active\\ activity--time combinations};
                
                \node[box, below of=active] (priority) {Select combinations from highest active priority class};
                
                \node[box, below of=priority] (time) {Select earliest time slot};
                
                \node[box, below of=time] (workload) {Select activity with lowest weighted workload};
                
                \node[box, below of=workload] (feasible) {Evaluate feasible volunteers and assignment intervals};
                
                \node[box, below of=feasible] (assign) {Assign volunteer with earliest feasible start time};
                
                \node[box, below of=assign] (update) {Update workloads and active set};
                
                \draw[arrow] (active) -- (priority);
                \draw[arrow] (priority) -- (time);
                \draw[arrow] (time) -- (workload);
                \draw[arrow] (workload) -- (feasible);
                \draw[arrow] (feasible) -- (assign);
                \draw[arrow] (assign) -- (update);
                \draw[arrow] (update.west) -- ++(-0.6,0) |- (priority.west);
                
                \end{tikzpicture}
            
            \caption{Decision structure of the proposed priority-driven constructive heuristic. 
            At each iteration, the algorithm selects an activity--time combination according to the priority, time, and workload rules before evaluating feasible volunteer assignments.} 
            \label{fig:heuristic_overview}
            \end{figure}

        Figure~\ref{fig:heuristic_overview} summarizes the resulting decision workflow of the constructive heuristic.

            \begin{algorithm}[H]
                \caption{Preprocessing: Sort volunteers}
                \label{alg:preprocessing}
                
                \KwIn{
                    Volunteer set $V$; parameters $d_{a,t}$, $req_{a,c}$, $cap_{v,c}$
                    }        
                \KwOut{Reordered volunteer set $V$}        
                \BlankLine 
                
                \tcp{Compute scarcity score for each volunteer}        
                \ForEach{$v \in V$}{
                    $s_v \leftarrow \min_{a,t:\; req_{a,c}=1,\; cap_{v,c}=1} d_{a,t}$\;
                }                
                \BlankLine                       
                Sort volunteers in nondecreasing order of $s_v$\; 
                \Return reordered $V$        
            \end{algorithm}         
        
            \begin{algorithm}[H]
                \caption{Priority-driven constructive heuristic for the SVCP}
                \label{alg:svcp-heuristic}
                \DontPrintSemicolon        
                \KwIn{Model parameters and data}
                \KwOut{Volunteer assignment and schedule}        
                \BlankLine
                $V \leftarrow \textsc{Preprocessing}(\cdot)$\;
                \BlankLine
                                   
                Initialize $o_{v,a,t}$ for all $v,a,t$\;
                Compute workloads $L_{a,t}(\mathbf{x})$ and weighted workloads $\sigma_{p_a,p_a+1} L_{a,t}(\mathbf{x})$\;
                $\mathcal{A}^{\text{active}} \leftarrow \{(a,t)\mid r_{a,t}=1,\; L_{a,t}(\mathbf{x})<1\}$\;                         
                \BlankLine         
               
                \While{$\mathcal{A}^{\text{active}} \neq \emptyset$}{
                
                    $\mathcal{A}_{k^*}^{\text{active}} \leftarrow \textsc{SelectSubsetHighestPriorityClass}(\mathcal{A}^{\text{active}})$
                    \BlankLine                     
                     
                    $(a^*,t^*)\leftarrow \textsc{SelectBestActiveCombination}(\mathcal{A}_{k^*}^{\text{active}})$\;    
                    \BlankLine                    
                    
                    $V^{cand} \leftarrow$ \textsc{VolunteerFeasibilityEvaluation}$(\cdot)$\;  
                    \BlankLine      
                    \tcp{Assignment or removal}
                      \If{$V^{cand} \neq \emptyset$}{
                        Select $(v^*,t_s^*,t_e^*) \in V^{cand}$ with minimal $t_s^*$\;
                        \For{$t=t_s^*$ \KwTo $t_e^*$}{
                        $x_{v^*,a^*,t} \leftarrow 1$\;
                        }
                        Update workloads $L_{a,t}(\mathbf{x})$ and weighted workloads $\sigma_{p_a,p_a+1} L_{a,t}(\mathbf{x})$\;
                        Remove all $(a^*,t)$ from $\mathcal{A}^{\text{active}}$ with $L_{a^*,t}(\mathbf{x})=1$\;
                      }
                      \Else{
                        Remove $(a^*,t^*)$ from $\mathcal{A}^{\text{active}}$\;
                      }
                }
                \Return Volunteer assignment and scheduling
            \end{algorithm}
            
        At each iteration, the algorithm selects an activity--time combination from the highest priority class that still requires volunteers.
        Within this subset, the earliest time slot is considered first and the activity with the lowest weighted workload is chosen.

        Candidate volunteers are then evaluated with respect to capability compatibility, availability, travel constraints, and remaining working time.
        Among feasible candidates, the volunteer with the earliest feasible start time is assigned for the maximal feasible interval.
        If no feasible volunteer exists, the activity--time combination is removed from the active set.            

        Detailed pseudocode for these subroutines is provided in~\ref{sec:subroutines}.

    \subsection{Algorithmic properties}
    \label{sec:heuristic_properties}
    
        The priority-driven structure of the heuristic induces several useful algorithmic properties.

        \textbf{Priority preservation.}
        The restriction to activity--time combinations belonging to the highest currently active priority class preserves the lexicographic structure of the objective hierarchy. 
        Assignments that contribute to higher-priority objectives are therefore always generated before assignments associated with lower-priority classes.
        
        \textbf{Early task coverage.}
        The preference for earlier time slots promotes early task coverage. 
        Since the objective function assigns larger weights to earlier time slots, the earliest-time selection rule implicitly aligns the constructive decisions with the temporal weighting structure of the model.
        
        \textbf{Scarce capability protection.}
        The preprocessing step that orders volunteers according to capability scarcity helps protect scarce capabilities during the construction process. 
        Volunteers whose capabilities are required for activities with limited supply are considered earlier, reducing the risk that these volunteers are allocated to tasks that could also be performed by more common capabilities.
        
        \textbf{Workload balancing.}
        The workload-based activity selection rule gradually reduces staffing imbalances between activities. 
        Activities with relatively low staffing ratios are prioritized, which approximates the workload balancing objectives of the optimization model during the constructive process.
    
    \subsection{Complexity analysis}
    \label{sec:heuristic_complexity}
    
        The computational complexity of the heuristic is dominated by the volunteer feasibility evaluation step.
        For each selected activity–time combination, the feasibility of volunteers is evaluated to identify feasible candidates.
        In the worst case, this requires evaluating all volunteers in $V$.
        Since each activity--time combination $(a,t)$ can enter the set $\mathcal{A}^{\text{active}}$ at most once and is removed once it is processed or becomes fully assigned, the total number of $|\mathcal{A}^{active}|$ is bounded by $|A||T|$.

        Consequently, the worst-case time complexity of the heuristic is $O(|A||T||V|)$.
        The complexity therefore scales linearly in the number of volunteers and the number of activity--time combinations.
        In practice, the number of feasibility checks is typically smaller because many activity--time combinations become satisfied early during the constructive process.
        This supports the applicability of the heuristic in time-critical disaster response settings where coordination decisions must be updated frequently.

\section{Computational study}
\label{sec:computational_study} 

    The computational study evaluates the performance of the proposed priority-driven constructive heuristic and compares it with solutions obtained by solving the exact optimization model using the off-the-shelf solver Gurobi.
    
    The computational analysis is structured around three research questions.
        \begin{itemize}
            \item[\textbf{Q1}] How closely does the proposed heuristic approximate the solutions obtained by exact optimization with respect to the lexicographic objective hierarchy?
        
            \item[\textbf{Q2}] How does the computational runtime of the heuristic compare to exact optimization across dynamically evolving disaster response scenarios?
        
            \item[\textbf{Q3}] What runtime–quality trade-off is achieved by the heuristic and under which conditions heuristic solutions provide an effective alternative to exact optimization?
        \end{itemize}    

    Because the objectives are evaluated lexicographically, solution quality is assessed by comparing the objective values hierarchically according to the lexicographic objective structure. 
    The heuristic wall-clock time includes the time required for initialization, preprocessing (Algorithm~\ref{alg:preprocessing}), and execution of the priority-driven constructive heuristic (Algorithm~\ref{alg:svcp-heuristic}).
    The off-the-shelf solver wall-clock time consists of the time required for initialization, model construction, and optimization.

    To ensure comparability with previous work, we adopt the scenario generation procedure introduced by \citet{Sperling2022}. 
    The experimental setup is based on empirical data collected during the 2013 flood in Halle (Saale), Germany, and reflects operational characteristics of spontaneous volunteer coordination observed in practice.
    The set of operational tasks derived from this event is reported in~\ref{subsec:tasks}.
    
    \subsection{Experimental setup}
	\label{subsec:experimental_design}
    
        The dynamic scenario framework follows \citet{Sperling2022}.
        Volunteer coordination decisions are updated every 30 minutes, reflecting the operational requirements reported by practitioners.
        Each update generates a new optimization instance while taking into account assignments determined in previous instances, resulting in 20 sequential instances per scenario.
        The planning horizon of each instance spans $T=48$ time slots representing a period of 24 hours.
        Consequently, each scenario covers an overall time span of $33.5$ hours.

        To account for stochastic variability in volunteer arrivals, each scenario is simulated ten times using different random seeds.
        Each simulation run generates a distinct sequence of volunteer arrivals and availability patterns while keeping the scenario parameters fixed.
        Overall, the computational study therefore contains $20 \cdot 16 \cdot 10 = 3\,200$ optimization instances.
        This large number of instances ensures that the experimental results are statistically robust and capture a wide range of operational conditions.
        The activity types and their associated capability requirements used in the experiments are listed in~\ref{subsec:activities_capabilities}.

        Constant parameter values used in all scenarios are summarized in Table~\ref{tab:parameters}.
        The configuration follows the empirical setting in~\citet{Sperling2022}.
            \begin{table}[h]
            \centering
            \caption{Constant parameter values used in the computational study
            (adopted from \citet{Sperling2022}).}
            \label{tab:parameters}
                \begin{tabular}{l l p{0.45\linewidth}}
                \toprule
                Category & Parameter & Value \\
                \midrule                
                    Planning horizon
                    & $T$
                    & $48$ time slots ($24$ hours) \\
                    
                    Capabilities
                    & $C$
                    & $6$ \\
                    
                    Priorities
                    & $P$
                    & $3$ priority levels \\
                    
                    Priority classes
                    & $K$
                    & $2$ classes \\
                    
                    Priority classes
                    & $\mathcal{P}_1,\mathcal{P}_2$
                    & $\mathcal{P}_1=\{1,2\}$,\;
                    $\mathcal{P}_2=\{3\}$ \\
                    
                    Priority balancing
                    & $\sigma_{1,2}$
                    & $\frac{1}{3}$ \\
                    
                    Assignment duration
                    & $\tau_{min}$
                    & $4$ time slots ($2$ hours) \\
                    
                    Maximum working time
                    & $\bar{\tau}_{max}$
                    & $16$ time slots ($8$ hours) \\
                    
                    Initial travel time
                    & $\tau_{v,a}$
                    & $2$ time slots ($1$ hour) \\
                    
                    Time weighting
                    & $w_t$
                    & $1-\frac{t-1}{T}$ \\                
                \bottomrule
                \end{tabular}
            \end{table}

        In addition, to calculate the travel time between different activity locations $s_{a,a'}$, we set the speed to 10 km/h, following the parameterization used in \citet{Sperling2022}.
            
    \subsection{Scenario and instance generation}
	\label{subsec:instance_generation}
        The experimental scenarios are generated using the stochastic simulation procedure proposed by \citet{Sperling2022}.
        The generation process captures the dynamic arrival of tasks and SVs during disaster response operations.
        SV arrivals follow a Poisson-based arrival process controlled by the parameter $\lambda$, while the maximum number of volunteers per scenario is bounded by the scenario configuration.
        As a result, the total number of available volunteers gradually increases over the 20 instances of a scenario until the specified maximum number of volunteers is reached.

        Scenario variations are defined along four dimensions: the number of newly arriving tasks, the intensity of volunteer arrivals, the maximum number of available volunteers, and the probability that volunteers possess specific capabilities.        
            \begin{table}[h]
            \small
                \centering
                \caption{Experimental factors used in the scenario generation (based on \citet{Sperling2022}).}
                \label{tab:variation-parameter}
                \begin{tabular}{{p{0.3\linewidth} p{0.15\linewidth} p{0.45\linewidth}}}
                \toprule
                Factor & Levels & Description \\
                \midrule
                Added tasks per instance 
                & $\{1,2\}$ 
                & Number of new tasks that appear in each instance. \\
                
                Volunteer arrival parameter $\lambda$ 
                & $\{7,11\}$ 
                & Parameter of the Poisson distribution controlling the arrival intensity of volunteers. \\
                
                Maximum number of volunteers 
                & $\{5\,000,10\,000\}$ 
                & Upper bound on the total number of volunteers available in a scenario. \\
                
                Capability probability 
                & $\{0.3,0.5\}$ 
                & Probability that a volunteer possesses a particular capability. \\
                
                \bottomrule
                \end{tabular}
            \end{table}
        Combining the parameter levels across the four experimental factors results in a design with $2^4=16$ distinct scenarios.
        
        This simulation procedure reproduces the characteristic surge and decline of SV arrivals observed during the 2013 flood response in Halle (Saale), Germany, as documented in \citet{Sperling2022}.
        This arrival pattern reflects the typical convergence behavior of spontaneous volunteers observed in real disaster response operations.
        Across all scenarios, the resulting instances contain up to 5\,000 or 10\,000 volunteers and up to 27 tasks. 
        These tasks correspond to up to $A = 87$ activities and a maximum demand of 3\,030 volunteers. 
        Considering the planning horizon of $T = 48$ time slots, this results in up to $A\cdot T = 4\,176$ activity–time combinations.

        The computational experiments are based on a full-factorial scenario design with four experimental factors: the maximum number of volunteers, the number of newly arriving tasks per instance, the probability that a volunteer possesses a required capability, and the arrival intensity of volunteers controlled by the parameter $\lambda$.

        Each factor is evaluated at two levels, resulting in $2^4=16$ distinct scenarios summarized in Table~\ref{tab:scenario_configurations}.
            \begin{table}[h]
            \centering
            \caption{Scenario configurations used in the computational study following \citet{Sperling2022}.}
            \label{tab:scenario_configurations}
            \begin{tabular}{p{0.15\linewidth} p{0.15\linewidth} p{0.1\linewidth} p{0.2\linewidth} p{0.1\linewidth}}
            \toprule
            Scenario & Volunteers & Added tasks & Capability probability & Arrival rate $\lambda$ \\
            \midrule
            1  & 5\,000  & 1 & 0.3 & 7  \\
            2  & 5\,000  & 1 & 0.3 & 11 \\
            3  & 5\,000  & 1 & 0.5 & 7  \\
            4  & 5\,000  & 1 & 0.5 & 11 \\
            5  & 5\,000  & 2 & 0.3 & 7  \\
            6  & 5\,000  & 2 & 0.3 & 11 \\
            7  & 5\,000 & 2 & 0.5 & 7  \\
            8  & 5\,000 & 2 & 0.5 & 11 \\
            9  & 10\,000 & 1 & 0.3 & 7  \\
            10 & 10\,000 & 1 & 0.3 & 11 \\
            11 & 10\,000 & 1 & 0.5 & 7  \\
            12 & 10\,000 & 1 & 0.5 & 11 \\
            13 & 10\,000 & 2 & 0.3 & 7  \\
            14 & 10\,000 & 2 & 0.3 & 11 \\
            15 & 10\,000 & 2 & 0.5 & 7  \\
            16 & 10\,000 & 2 & 0.5 & 11 \\
            \bottomrule
            \end{tabular}
            \end{table}
        
    \subsection{Solution methods}
	\label{subsec:solution_methods}    
        Two solution approaches are evaluated in the computational study.
        First, the proposed priority-driven constructive heuristic described in Section~\ref{sec:heuristic} is applied to each instance.
        Second, we solve the mixed-integer optimization model introduced by \citet{Sperling2022} using the off-the-shelf solver Gurobi. 
        Since the model is formulated as a lexicographic multi-objective problem, each objective is optimized sequentially. 
        A solver time limit of 5 minutes is imposed for each objective, resulting in a maximum solving time of 20 minutes per instance.
        In addition, an overall wall-clock time limit of 30 minutes is enforced for each instance in order to reflect the operational requirement that coordination decisions must be updated every 30 minutes during disaster response operations.
        
        The heuristic is not subject to this time limit because its runtime remains well below the operational decision interval in all experiments.
        The reported runtime of the heuristic includes initialization, preprocessing, and the execution of the main algorithm.        

        It is important to note that the time limit imposed in Gurobi applies only to the optimization phase and does not include the time required for model construction. 
        Consequently, the total wall-clock time observed in our experiments includes both model generation and solver time.

        In contrast, the computational results reported in \citet{Sperling2022} focus on solver runtime only. 
        Our experiments indicate that a substantial portion of the total wall-clock time can be spent on model generation rather than on the optimization itself. 
        Since the SVCP requires the repeated construction of large mixed-integer optimization models in each planning instance, model generation can become a significant computational overhead in rolling-horizon decision environments.        

    \subsection{Computational environment}
	\label{subsec:computational_environment}
        Experiments were executed on the Noctua 2 high-performance computing system at the Paderborn Center for Parallel Computing (PC²), equipped with AMD EPYC 7763 processors.

        Both the heuristic algorithm and the optimization model were implemented in Python (version 3.12.3).
        The exact optimization model was solved using the off-the-shelf solver Gurobi (version 12.0.2).

        Experiments were submitted through the SLURM workload manager.
        The heuristic algorithm was executed using a single CPU thread, as its design does not rely on parallelization. 
        In contrast, the exact optimization was configured to use 16 parallel threads in order to exploit the parallel branch-and-bound capabilities of Gurobi.
        
        All experiments were repeated using identical hardware configurations to ensure reproducibility of the runtime measurements.

    \subsection{Results}
    \label{subsec:results}    
        This section evaluates the performance of the proposed heuristic in comparison with the exact optimization model.
        The analysis focuses on three aspects: solution quality, computational runtime, and the resulting runtime–quality trade-off across the simulated disaster response scenarios.
        
        The computational experiments are based on the scenario design summarized in Table~\ref{tab:scenario_configurations}.
        For each of the 16 scenarios, 20 sequential rolling-horizon instances are solved using both solution approaches.
        Each scenario is simulated ten times using different random seeds in order to capture stochastic variability in volunteer arrivals.
        Unless stated otherwise, the reported results represent the median performance across these stochastic replications.
        
        The heuristic solutions and the solutions obtained by the exact optimization are evaluated using the lexicographic objective hierarchy defined in Section~\ref{sec:model}.
        The relative gap is computed with respect to the best feasible solution returned by the exact optimization at termination.
        If the solver reaches the imposed time limit, the best incumbent solution is used as reference.
        
        Instances for which the solver reached the time limit and additional descriptive statistics are reported in~\ref{sec:add_results}.
        The heuristic parameters remain fixed across all scenarios and were not tuned for individual instances.

    \subsubsection{Solution quality}
    \label{subsec:solution_quality}    
\begin{figure}[htbp]
            \centering
                \begin{subfigure}[b]{0.49\textwidth}
                \centering
                \includegraphics[page=1,width=\textwidth]{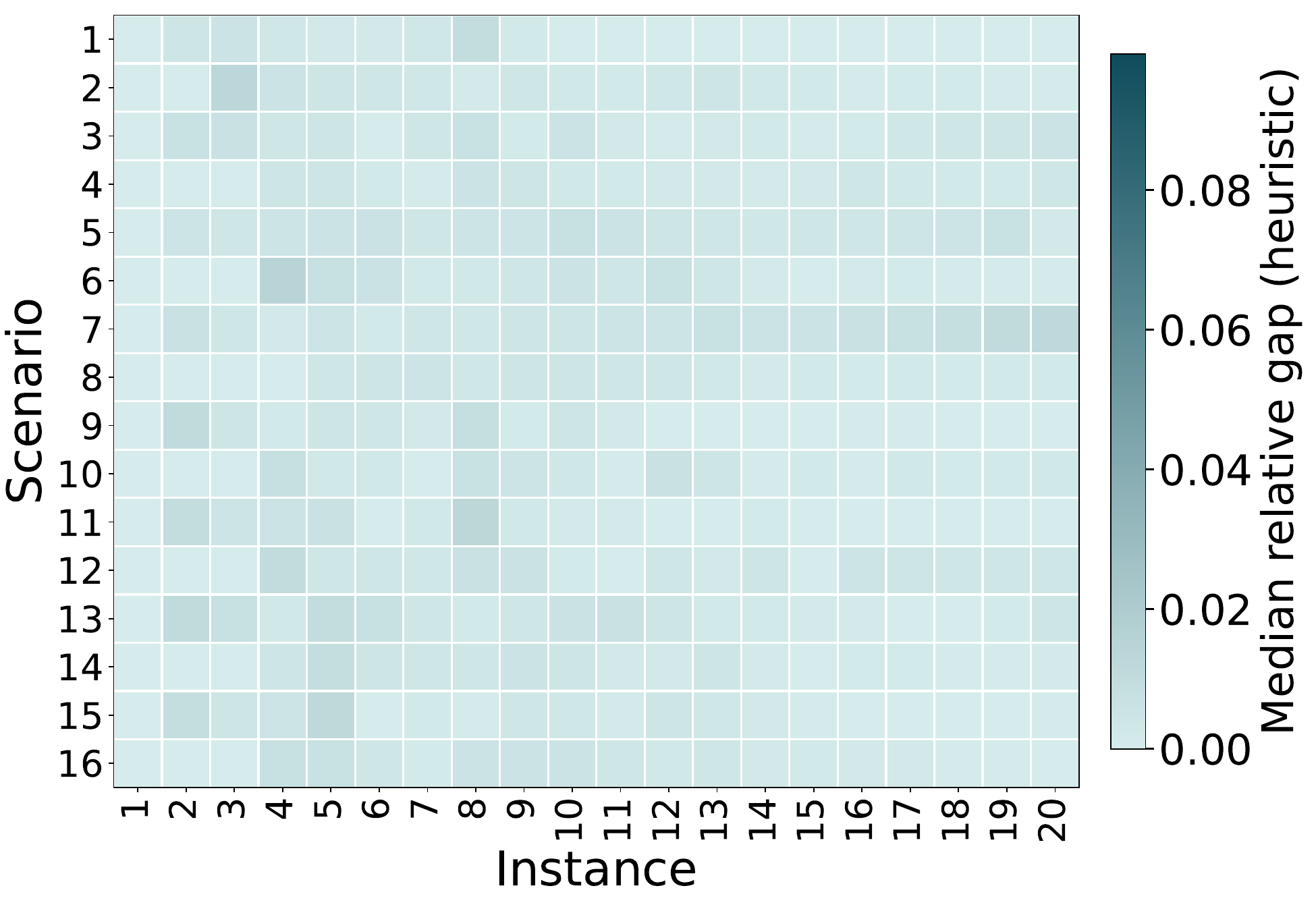}
                \caption{Objective 1}
                \label{fig:gap_obj1}
                \end{subfigure}
                \hfill
                \begin{subfigure}[b]{0.49\textwidth}
                \centering
                \includegraphics[page=2,width=\textwidth]{figures/F1_instance_heatmaps_median_grouped_scale.pdf}
                \caption{Objective 2}
                \label{fig:gap_obj2}
                \end{subfigure}
                
                \vspace{0.4cm}
                
                \begin{subfigure}[b]{0.49\textwidth}
                \centering
                \includegraphics[page=3,width=\textwidth]{figures/F1_instance_heatmaps_median_grouped_scale.pdf}
                \caption{Objective 3}
                \label{fig:gap_obj3}
                \end{subfigure}
                \hfill
                \begin{subfigure}[b]{0.49\textwidth}
                \centering
                \includegraphics[page=4,width=\textwidth]{figures/F1_instance_heatmaps_median_grouped_scale.pdf}
                \caption{Objective 4}
                \label{fig:gap_obj4}
                \end{subfigure}
                \vspace{-0.4cm}
                \begin{subfigure}[b]{0.8\textwidth}
                \centering
                \includegraphics[page=5,width=\textwidth]{figures/F1_instance_heatmaps_median_grouped_scale.pdf}                
                \end{subfigure}
            
            \caption{Median relative gaps of the heuristic across scenarios and rolling-horizon instances for the four objectives.
            Rows correspond to the 16 scenario configurations defined in Table~\ref{tab:scenario_configurations}, while columns represent the sequential rolling-horizon instances.
            The color intensity indicates the magnitude of the median relative gap computed across the ten stochastic replications of each scenario–instance combination.
            Objectives 1–4 correspond to priority coverage, time-weighted task fulfillment, inter-activity workload balancing, and intra-priority workload balancing, respectively.}
            \label{fig:gap_all}
            \end{figure}

        We first analyze the solution quality achieved by the proposed priority-driven constructive heuristic.
        Figure~\ref{fig:gap_all} reports the median relative gaps of the heuristic across all scenarios and rolling-horizon instances for the four objectives.
        
        \paragraph{Primary objectives}        
        For the first two objectives, which capture priority coverage and time-weighted task fulfillment, the heuristic achieves solutions that are very close to those obtained by the exact optimization.
        Across almost all scenarios and instances, the median gaps remain below one percent for Objective~1 and below three percent for Objective~2.
        
        These results indicate that the priority-driven selection rules embedded in the heuristic successfully reproduce the lexicographic priority structure of the optimization model.
        In particular, the restriction to the highest active priority class ensures that assignments contributing to higher-priority objectives are generated before lower-priority tasks are considered.
        
        \paragraph{Workload balancing objectives}        
        Larger deviations can be observed for Objective~3, which captures the workload ratio balancing between activities.
        The largest gaps occur in the earliest instances of the rolling-horizon process, where the number of available volunteers is still limited.
        In these early stages, the feasible assignment space is strongly restricted by capability compatibility and minimum assignment durations.
        As a result, the constructive heuristic prioritizes high-priority coverage and early assignments, while workload balancing can only be approximated.
        
        As additional volunteers arrive in later instances, the gaps decrease substantially, and the heuristic solutions converge toward those obtained by the exact optimization.
        This effect occurs because the increasing number of volunteers expands the feasible assignment space, allowing the heuristic to better satisfy workload balancing objectives.
        
        Occasionally large relative gaps appear for Objective~3.
        These values occur when the optimal objective value is close to zero, which inflates ratio-based gap measures even for small absolute deviations.
        Such effects are inherent to relative gap metrics for minimization objectives and therefore need to be interpreted with caution.
        
        Objective~4, which captures workload balancing within the same priority level, shows moderate deviations across scenarios.
        In most instances, the gaps remain below 20\%, and they decrease gradually as the volunteer pool grows over time.
        
        \paragraph{Impact of rolling-horizon dynamics and scenario parameters}        
        The gap patterns are strongly influenced by the rolling-horizon dynamics of the coordination process.
        Early instances typically exhibit larger gaps because only a limited number of volunteers are available, and the feasible assignment space is highly constrained.
        In later instances, the volunteer pool grows substantially, which increases assignment flexibility and allows the heuristic to better approximate the balancing objectives.
        
        The solution quality also varies systematically across the scenario configurations reported in Table~\ref{tab:scenario_configurations}. 
        Scenarios with larger volunteer pools and higher task arrival rates tend to exhibit larger gaps in the early rolling-horizon instances. 
        In these scenarios, the assignment space becomes substantially larger, which increases the difficulty of balancing workload objectives while simultaneously respecting capability compatibility and minimum assignment durations. 
        However, as the volunteer pool grows over time, the heuristic rapidly converges toward the solutions obtained by the exact optimization, even in these challenging scenarios.
        
        The magnitude of the gaps also depends on the scenario configuration.
        Scenarios with large volunteer pools and higher task arrival rates (e.g., Scenarios~13–16 in Table~\ref{tab:scenario_configurations}) generate substantially larger assignment spaces and therefore represent the most challenging instances for both solution approaches.
        Nevertheless, even in these demanding scenarios the heuristic maintains stable performance across the majority of instances.
        
        \paragraph{Solver time limits}        
        The results further indicate that the exact optimization frequently reaches the imposed time limit, particularly in scenarios with large volunteer pools.
        In such cases the reported gaps are computed relative to the best incumbent solution obtained by the solver.
        Consequently, the heuristic occasionally produces solutions of comparable quality even when the exact optimization cannot fully optimize the corresponding instance within the available time budget.
        
        Overall, the results demonstrate that the proposed heuristic closely approximates the lexicographic objective hierarchy of the optimization model while maintaining robust solution quality across a wide range of dynamically evolving disaster response scenarios.

    \subsubsection{Runtime comparison}
    \label{subsec:runtime}
    
\begin{figure}[htbp]
        \centering
            \includegraphics[width=0.85\textwidth]{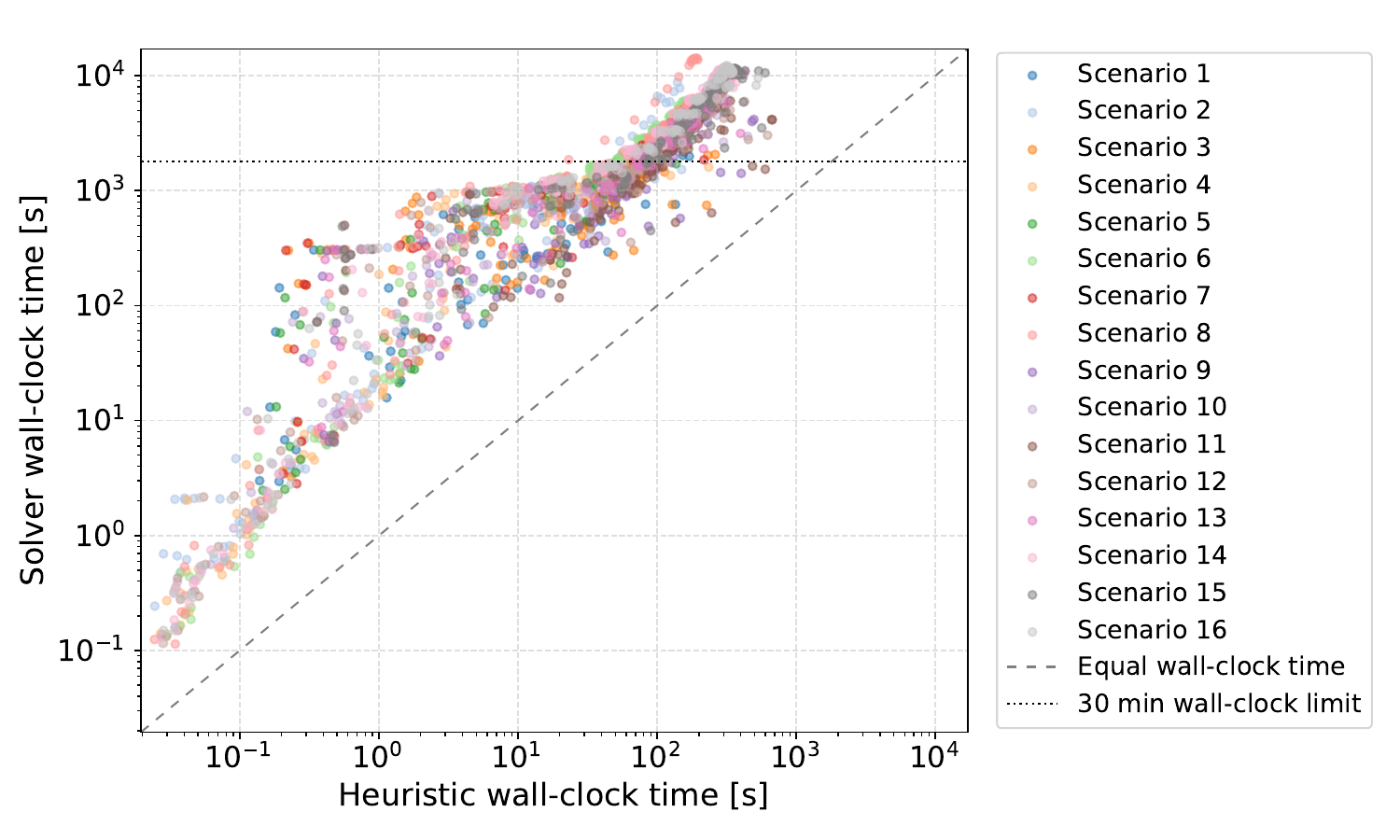}
        \caption{Runtime comparison between the exact optimization and the proposed heuristic.}
        \label{fig:runtime_scatter}
        \end{figure}

        Figure~\ref{fig:runtime_scatter} illustrates the runtime comparison across all scenario–instance combinations.
        Each point represents the wall-clock runtime of a single scenario–instance–run combination.
        The diagonal line indicates equal wall-clock time, while the horizontal line marks the operational decision interval of 30 minutes imposed by the rolling-horizon coordination process.
        
        The results show a clear separation between the two solution approaches.
        Across the vast majority of instances, the solver requires substantially more computation time than the heuristic.
        Most observations lie far above the diagonal line, indicating that the heuristic is typically one to two orders of magnitude faster than the exact optimization.
        
        This difference becomes particularly evident when considering the operational time limit.
        In total, the exact optimization exceeds the 30-minute wall-clock time limit in 61.8\% of all instances, whereas the heuristic never violates this limit.
        Across all instances, the heuristic achieves a median runtime speedup of approximately $28\times$ compared to the exact optimization.

        The runtime differences also vary systematically across the scenario configurations reported in Table~\ref{tab:scenario_configurations}.
        Scenarios with large volunteer pools and higher task arrival rates, in particular Scenarios~13–16, generate substantially larger assignment spaces and therefore represent the most computationally demanding instances.
        In these scenarios the number of feasible volunteer–activity assignments grows rapidly during the rolling-horizon process, which leads to large mixed-integer models and frequent violations of the imposed solver time limit.
        In contrast, scenarios with smaller volunteer pools (Scenarios~1–8) typically remain computationally manageable for the exact optimization during the early instances of the rolling horizon.
        
        A closer examination of the instance dynamics reveals a consistent pattern across scenarios.
        The runtime of the exact optimization approach increases steadily across the sequential rolling-horizon instances.
        This effect is caused by the growing number of available volunteers and feasible assignments over time, which leads to increasingly large mixed-integer models in later planning stages.
        These results indicate that instance size is a key driver of the difficulty of the SVCP, affecting both exact optimization and heuristic solution quality.
        
        Detailed runtime breakdowns provided in~\ref{sec:add_results} show that a substantial portion of the total wall-clock time is spent on model construction rather than on the optimization phase itself.
        Since the SVCP must be repeatedly rebuilt and solved for each planning instance, model generation becomes a major computational bottleneck in rolling-horizon environments.
        
        In contrast, the proposed heuristic operates directly on the problem data and therefore avoids repeated model generation.
        As a result, its runtime remains stable across instances and scenarios and consistently stays well below the operational decision interval.
        
        Overall, the runtime results highlight a fundamental scalability limitation of the exact optimization approach in dynamic disaster response settings.
        While the exact optimization can produce high-quality solutions for smaller instances, the repeated construction and solution of large mixed-integer models becomes computationally expensive as the system evolves.
        The proposed heuristic, in contrast, maintains consistently low wall-clock times while producing solutions of comparable quality, making it well suited for real-time decision support in rolling-horizon volunteer coordination environments.
        The observed runtime behavior confirms that the heuristic scales linearly with the size of the volunteer pool and the number of activity–time combinations, which is consistent with the theoretical complexity analysis presented in Section~\ref{sec:heuristic}.
        
        Additional runtime statistics are summarized in Tables~\ref{tab:runtime_statistics} and \ref{tab:runtime_exceedance_scenarios} in the appendix.

    \subsubsection{Runtime-quality tradeoff}
    \label{subsec:scalability}
        
        The final analysis investigates the trade-off between computational runtime and solution quality achieved by the proposed heuristic.
        Because the first two objectives represent the highest-priority goals in the lexicographic objective hierarchy, their solution quality is preserved by the heuristic in almost all instances.
        Trade-offs therefore primarily arise for the workload-balancing objectives.
        
\begin{figure}[htbp]
            \centering
            \includegraphics[page=1,width=0.9\textwidth]{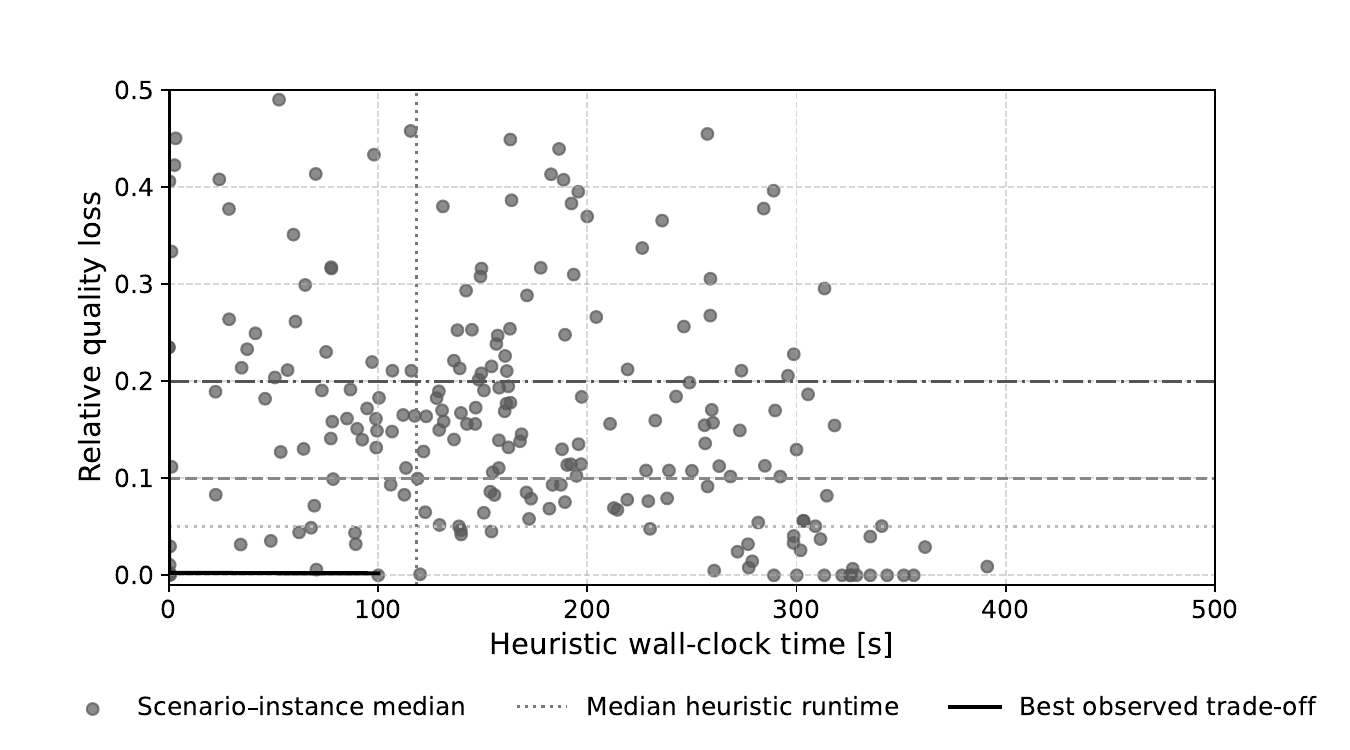}
            \caption{Runtime–quality trade-off of the heuristic for Objective~3.
            Each point represents the median performance across stochastic
            replications for a scenario–instance combination.
            The curve indicates the best observed runtime–quality trade-off,
            while the vertical line marks the median heuristic runtime.}
            \label{fig:tradeoff}
            \end{figure}

        Figure~\ref{fig:tradeoff} illustrates the relationship between heuristic runtime and the relative quality loss for Objective~3.
        Each point represents the median performance across the stochastic replications of a scenario–instance combination.
        
        Table~\ref{tab:tradeoff_summary} summarizes the distribution of heuristic wall-clock times and the corresponding quality losses across all scenario–instance combinations.
        While the maximum observed relative gap appears large, such values occur only in rare cases where the optimal objective value is close to zero.
        In these situations, relative gap measures may exaggerate deviations even though the corresponding absolute differences in workload balancing remain small.
        
            \begin{table}[h]
            \centering
            \caption{Summary statistics of the runtime–quality trade-off for the heuristic (Objective~3).}
            \label{tab:tradeoff_summary}
            \begin{tabular}{lcc}
            \toprule
            Statistic & Wall-clock time [s] & Relative quality loss \\
            \midrule
            Median & 118.23 & 0.21 \\
            25th percentile & 22.77 & 0.10 \\
            75th percentile & 192.25 & 0.68 \\
            Maximum & 391.15 & 31.46 \\
            \bottomrule
            \end{tabular}
            \end{table}
        
            \begin{flushleft}
            \footnotesize
            Note: Very large relative gaps may occur when the optimal objective
            value is close to zero, which inflates ratio-based gap measures even
            for small absolute deviations.
            \end{flushleft}
        
        The results show that the heuristic consistently produces solutions within short computation times while maintaining high solution quality.
        Across most scenario–instance combinations, the heuristic runtime remains below approximately four minutes, with a median runtime of 118 seconds.
        This indicates that the heuristic typically produces solutions well within the operational decision interval of 30 minutes considered in the rolling-horizon coordination process.
        
        At the same time, the observed quality losses for the workload-balancing objective remain moderate in most cases.
        The majority of instances exhibit relative gaps below approximately 20\%, indicating that the heuristic captures the main balancing structure of the optimization model despite its constructive nature.
        
        The Pareto envelope in Figure~\ref{fig:tradeoff} highlights the best observed runtime–quality combinations across all scenarios.
        It shows that near-zero quality losses can already be achieved at moderate runtimes, whereas additional runtime improvements provide only limited reductions in the remaining quality loss.
        
        Overall, the results confirm that the proposed heuristic achieves a favorable runtime–quality balance.
        While exact optimization may still produce slightly better workload-balancing decisions in some instances, the heuristic delivers solutions of comparable quality within very short runtimes.
        This makes the approach particularly suitable for real-time decision support in rolling-horizon volunteer coordination during disaster response operations.

        These observations suggest that heuristic solution approaches are particularly advantageous in large-scale and dynamically evolving disaster response settings where exact optimization becomes computationally prohibitive.

\section{Conclusion}
\label{sec:conclusion}

    This paper studied the coordination of spontaneous volunteers in large-scale disaster response operations. 
    Building on the multi-objective optimization model introduced by Sperling and Schryen (2022), we developed a priority-driven constructive heuristic for the spontaneous volunteer coordination problem (SVCP). 
    The heuristic exploits structural properties of the optimization model, including the lexicographic priority hierarchy, capability scarcity among volunteers, and workload balancing requirements.

    A large-scale computational study based on empirically grounded disaster response scenarios evaluated the performance of the proposed approach. 
    The experiments considered rolling-horizon coordination settings with up to 10\,000 volunteers, 27 tasks, and more than 4\,000 activity–time combinations per instance. 
    Across the 3\,200 simulated instances, the heuristic closely approximates the solutions obtained by exact optimization for the primary objectives of priority coverage and time-weighted task fulfillment, with median relative gaps below 1\% and 3\%, respectively. 
    For the workload-balancing objectives, the heuristic exhibits larger deviations in early rolling-horizon instances but converges toward the exact solutions as additional volunteers become available. 
    At the same time, the proposed heuristic achieves substantial computational advantages. 
    Across all experiments, the heuristic achieves a median runtime speedup of approximately $28\times$ compared to exact optimization, while the exact solver exceeds the operational decision limit of 30 minutes in more than 60\% of all instances. 
    In contrast, the heuristic consistently produces solutions within a few minutes, thereby satisfying the real-time decision requirements of disaster response coordination.

    These findings demonstrate that structure-exploiting heuristics can provide effective decision support for large-scale spontaneous volunteer coordination. 
    In particular, the proposed approach enables real-time decision making in rolling-horizon coordination environments where repeated exact optimization becomes computationally prohibitive.
    
    The results of this study should be interpreted in light of several limitations. 
    The computational experiments are based on empirically grounded scenarios derived from a single disaster response event and therefore reflect a specific operational setting. 
    Although the scenario design captures important characteristics of spontaneous volunteer coordination, additional empirical studies across different disaster contexts would further strengthen the generalizability of the findings.
    
    Several directions for future research emerge from this work. 
    First, future studies could investigate the applicability of the proposed coordination approach across different disaster types and organizational contexts, where task structures, volunteer arrival patterns, and coordination mechanisms may differ substantially.

    Second, the current model assumes a centralized coordination structure with a single volunteer coordinator. 
    In practice, volunteer coordination is often distributed across multiple organizations and coordination centers. 
    Extending the model to such multi-coordinator environments therefore represents an important avenue for future research.

    Third, hybrid solution approaches that combine structure-exploiting heuristics with exact optimization could be investigated. 
    For example, heuristic solutions could be used to initialize or guide exact solvers in large instances, potentially combining the computational efficiency of heuristics with the solution quality guarantees of optimization methods.

    Finally, future work could explore extensions of the model that explicitly incorporate uncertainty in volunteer arrivals and task demands, for example through stochastic optimization or online decision-making frameworks. 
    Such extensions may further improve decision support capabilities in highly dynamic disaster response environments.





\section*{Acknowledgements}
The authors gratefully acknowledge the funding of this project by computing time provided by the Paderborn Center for Parallel Computing (PC2).

\bibliographystyle{elsarticle-num-names} 
\bibliography{dom_vols}


\appendix
\section{Heuristic subroutines}
\label{sec:subroutines}
    
    The priority-driven constructive heuristic described in Section~\ref{sec:heuristic} consists of three main decision components:(i) the identification of activity--time combinations belonging to the currently highest priority class, (ii) the selection of the most critical activity--time pair based on temporal precedence and workload considerations, and (iii) the evaluation of feasible volunteer assignments. 
    The following algorithms implement the subroutines of the priority-driven constructive heuristic for the SVCP, as described in Algorithm~\ref{alg:svcp-heuristic}.

    \paragraph{Priority-based activity selection}
        Algorithm~\ref{alg:select-subset} identifies the subset of currently active activity--time combinations that belong to the highest priority class. 
        Since the objective structure of the optimization model follows a lexicographic hierarchy of priority classes, the heuristic restricts its search to combinations associated with the highest priority level that still requires additional volunteers. 
        This filtering step reduces the search space and ensures that assignments are generated in accordance with the priority structure of the model.
        
        Algorithm~\ref{alg:select-best-combination} selects a single activity--time pair from the subset identified in the previous step. 
        The selection rule first prioritizes earlier time slots to emphasize early task coverage and then chooses the activity with the smallest weighted workload. 
        The weighted workload incorporates the workload balancing parameters of the optimization model and therefore approximates the lexicographic balancing objectives of the original formulation.
    
             \begin{algorithm}[H]
                \caption{Select Subset with Highest Priority Class}
                \label{alg:select-subset}
                
                \KwIn{
                Active combinations $\mathcal{A}^{\text{active}}$, priority classes $\{1,\dots,K\}$ ($K$ is highest)
                }
                
                \KwOut{
                Subset $\mathcal{A}_{k^*}^{\text{active}} \subseteq \mathcal{A}^{\text{active}}$ containing only combinations of the highest present priority class
                }
                
                \BlankLine
                
                \tcp{Determine the highest priority class present in $\mathcal{A}^{\text{active}}$}
                $k^* \leftarrow \max \{\, p_a \mid (a,t) \in \mathcal{A}^{\text{active}} \,\}$\;
                
                \BlankLine
                
                \tcp{Select all combinations belonging to this priority class}
                $\mathcal{A}_{k^*}^{\text{active}} \leftarrow 
                \{\, (a,t) \in \mathcal{A}^{\text{active}} \mid p_a \in P_{k^*} \,\}$\;
                
                \BlankLine
                
                \Return $\mathcal{A}_{k^*}^{\text{active}}$\;
                
            \end{algorithm}

             \begin{algorithm}[H]
             \caption{Select Best Active Combination}
             \label{alg:select-best-combination}
                    \KwIn{
                    Active combinations with highest present priority class $\mathcal{A'}$ , wighted workload $\sigma_{p_a,p_a+1}\cdot L_{a,t}(\mathbf{x})$  
                    }
                    \KwOut{Selected combination $(a^*,t^*)$}
                    \BlankLine
                
                    \tcp{Select earliest time slot}
                    $t^* \leftarrow \min\{\, t \mid (a,t)\in \mathcal{A'} \,\}$\;
                    
                    \BlankLine
                    \tcp{Select activity with minimal weighted workload at $t^*$}
                    $a^* \leftarrow 
                    \arg\min_{a:\,(a,t^*)\in\mathcal{A'}}
                    \{\,\sigma_{p_a,p_a+1}\cdot L_{a,t^*}(\mathbf{x})\,\}$\;
                    
                    \BlankLine
                    
                    \Return $(a^*,t^*)$\;
                \end{algorithm}

    \paragraph{Volunteer feasibility evaluation}
        Algorithm~\ref{alg:volunteer-feasibility} evaluates which volunteers can feasibly be assigned to the selected activity--time pair. 
        For each volunteer, the procedure verifies capability compatibility, remaining working time, and other operational constraints. 
        If these conditions are satisfied, the maximal feasible assignment interval for the volunteer is determined using Algorithm~\ref{alg:max-feasible-interval}. 
        All volunteers for whom a feasible assignment interval exists are collected in a candidate set.
        
        Algorithm~\ref{alg:max-feasible-interval} determines the maximal contiguous time interval during which a volunteer can perform the selected task. 
        Starting from the reference time slot, the procedure iteratively extends the interval backward and forward while respecting availability, travel-time, and workload constraints. 
        If the resulting interval satisfies the minimum assignment duration requirement, the interval is returned as a feasible assignment period; otherwise, no feasible interval is reported.
            
                \begin{algorithm}[H]
                \caption{Volunteer Feasibility Evaluation}
                \label{alg:volunteer-feasibility}
                    \KwIn{
                    Volunteer set $V$; task--time pair $(a^*,t^*)$; 
                    assignment variables $x_{v,a,t}$; parameters 
                    $av_{v,t}$, $cap_{v,c}$, $req_{a,c}$,
                    $\tau_{v,a}$, $s_{a,a'}$, $\tau_{min}$, $\bar{\tau}_{max}$;
                    workloads $L_{a,t}(\mathbf{x})$}
                    
                    \KwOut{Candidate set $V^{cand}$}                    
                    \BlankLine
                    
                    $V^{cand} \leftarrow \emptyset$\;
                    \BlankLine
                    \ForEach{$v \in V$}{
                        \If{$\sum_{a=1}^{A}\sum_{t=1}^{T} x_{v,a,t} \ge \bar{\tau}_{max}$}{
                            \textbf{go to next $v$}\;
                        }
                        \If{$\sum_{c=1}^{C} req_{a^*,c}\cdot cap_{v,c} = 0$}{
                            \textbf{go to next $v$}\;
                        } 
                        
                        \BlankLine
                        \tcp{Select maximal interval for $v$ (Algorithm~\ref{alg:max-feasible-interval})}
                        $(t_s,t_e) \leftarrow$ 
                        \textsc{MaximalFeasibleInterval}$(\cdot)$\;                
                        \BlankLine
                        
                        \If{ $(t_s,t_e) \neq \emptyset$ }{
                            Add $(v,t_s,t_e)$ to $V^{cand}$\;
                        }  
                    } 
                    
                    \Return $V^{cand}$\;
                \end{algorithm}       
                                    
                \begin{algorithm}[H]
                \caption{Maximal Feasible Interval Detection}
                \label{alg:max-feasible-interval}
                    \KwIn{
                    Volunteer $v$; task $a^*$; reference time $t^*$; 
                    current assignment $\mathbf{x}$; parameters 
                    $av_{v,t}$, $\tau_{v,a}$, $s_{a,a'}$, $\tau_{min}$;
                    workloads $L_{a,t}(\mathbf{x})$
                    }
                    \KwOut{Maximal feasible interval $[t_s,t_e]$ or $\emptyset$}
                    \BlankLine
                    \tcp{Initialize interval boundaries}
                    $t_s \leftarrow t^*$,\;
                    $t_e \leftarrow t^*$\;
                    \BlankLine
                    \tcp{Extend interval backward}
                    \While{$av_{v,t_s-1}=1$
                    \textbf{and} travel/setup constraints induced by $\tau_{v,a^*}$ and $s_{a,a'}$ are satisfied
                    \textbf{and} $L_{a^*,t_s-1}(\mathbf{x}) < 1$}{
                      $t_s \leftarrow t_s - 1$\;
                    }
                    \BlankLine
                    \tcp{Extend interval forward}
                    \While{$av_{v,t_e+1}=1$
                    \textbf{and} travel/setup constraints induced by $\tau_{v,a^*}$ and $s_{a,a'}$ are satisfied
                    \textbf{and} $L_{a^*,t_e+1}(\mathbf{x}) < 1$}{
                      $t_e \leftarrow t_e + 1$\;
                    }
                    
                    \If{$t_e - t_s + 1 < \tau_{min}$}{
                      \Return $\emptyset$\;
                    }
                    
                    \Return $[t_s,t_e]$\;
                \end{algorithm}
                
    Together, these components allow the heuristic to approximate the lexicographic decision logic of the optimization model while maintaining low computational complexity.

\section{Additional computational results}
\label{sec:add_results}

    Additional figures and tables are provided in this appendix to support the reproducibility and transparency of the computational evaluation.
    In particular, the figures and tables reported here provide a detailed scenario-level view of heuristic solution quality across all simulated instances.

    \subsection{Additional quality analysis}
    \label{subsec:add_gap_instances}
        
        Tables~\ref{tab:gap_obj1_1}--~\ref{tab:gap_obj4_2} summarize heuristic solution quality at the instance and scenario level.
            \begin{table}[h]
            \centering
            \caption{Median relative gap per instance (Objective 1, Scenarios 1--10)}
            \label{tab:gap_obj1_1}
            \setlength{\tabcolsep}{3pt}
            \renewcommand{\arraystretch}{0.95}
            
                \begin{tabular}{l*{10}{S}}
                    \toprule
                    Scenario & {1} & {2} & {3} & {4} & {5} & {6} & {7} & {8} & {9} & {10} \\
                    \midrule
                    1  & 0.000 & 0.004 & 0.006 & 0.003 & 0.002 & 0.002 & 0.003 & 0.010 & 0.002 & 0.000 \\
                    2  & 0.000 & 0.000 & 0.013 & 0.005 & 0.005 & 0.003 & 0.003 & 0.001 & 0.004 & 0.003 \\
                    3  & 0.000 & 0.007 & 0.006 & 0.004 & 0.005 & 0.000 & 0.003 & 0.007 & 0.001 & 0.005 \\
                    4  & 0.000 & 0.000 & 0.000 & 0.004 & 0.005 & 0.002 & 0.001 & 0.006 & 0.005 & 0.003 \\
                    5  & 0.000 & 0.005 & 0.003 & 0.005 & 0.005 & 0.006 & 0.004 & 0.005 & 0.005 & 0.007 \\
                    6  & 0.000 & 0.000 & 0.000 & 0.015 & 0.007 & 0.006 & 0.002 & 0.003 & 0.003 & 0.005 \\
                    7  & 0.000 & 0.006 & 0.004 & 0.002 & 0.005 & 0.002 & 0.003 & 0.003 & 0.004 & 0.004 \\
                    8  & 0.000 & 0.000 & 0.000 & 0.000 & 0.003 & 0.004 & 0.005 & 0.003 & 0.004 & 0.004 \\
                    9  & 0.000 & 0.010 & 0.005 & 0.002 & 0.004 & 0.004 & 0.002 & 0.008 & 0.001 & 0.004 \\
                    10 & 0.000 & 0.000 & 0.000 & 0.008 & 0.002 & 0.003 & 0.000 & 0.007 & 0.005 & 0.003 \\
                    11 & 0.000 & 0.010 & 0.005 & 0.005 & 0.007 & 0.000 & 0.002 & 0.012 & 0.002 & 0.001 \\
                    12 & 0.000 & 0.000 & 0.000 & 0.010 & 0.003 & 0.004 & 0.003 & 0.007 & 0.006 & 0.001 \\
                    13 & 0.000 & 0.010 & 0.007 & 0.003 & 0.009 & 0.007 & 0.004 & 0.002 & 0.003 & 0.006 \\
                    14 & 0.000 & 0.000 & 0.000 & 0.004 & 0.009 & 0.004 & 0.004 & 0.004 & 0.005 & 0.005 \\
                    15 & 0.000 & 0.009 & 0.004 & 0.005 & 0.011 & 0.000 & 0.002 & 0.001 & 0.004 & 0.003 \\
                    16 & 0.000 & 0.000 & 0.000 & 0.007 & 0.007 & 0.003 & 0.001 & 0.005 & 0.005 & 0.005 \\
                    \bottomrule
                \end{tabular}
            \end{table}
            
            \begin{table}
            \centering
            \caption{Median relative gap per instance (Objective 1, Scenarios 11--20)}
            \label{tab:gap_obj1_2}
            \setlength{\tabcolsep}{3pt}
            \renewcommand{\arraystretch}{0.95}
            
                \begin{tabular}{l*{10}{S}}
                    \toprule
                    Scenario & {11} & {12} & {13} & {14} & {15} & {16} & {17} & {18} & {19} & {20} \\
                    \midrule
                    1 & 0.000 & 0.000 & 0.000 & 0.000 & 0.000 & 0.000 & 0.000 & 0.000 & 0.000 & 0.000 \\
                    2 & 0.002 & 0.003 & 0.004 & 0.002 & 0.002 & 0.001 & 0.001 & 0.001 & 0.001 & 0.001 \\
                    3 & 0.002 & 0.000 & 0.002 & 0.002 & 0.001 & 0.001 & 0.003 & 0.003 & 0.004 & 0.005 \\
                    4 & 0.002 & 0.002 & 0.002 & 0.001 & 0.002 & 0.004 & 0.002 & 0.002 & 0.002 & 0.004 \\
                    5 & 0.006 & 0.004 & 0.003 & 0.003 & 0.003 & 0.003 & 0.004 & 0.005 & 0.007 & 0.002 \\
                    6 & 0.003 & 0.007 & 0.003 & 0.001 & 0.000 & 0.001 & 0.001 & 0.001 & 0.001 & 0.000 \\
                    7 & 0.005 & 0.005 & 0.007 & 0.006 & 0.006 & 0.006 & 0.007 & 0.008 & 0.010 & 0.012 \\
                    8 & 0.004 & 0.003 & 0.003 & 0.001 & 0.002 & 0.001 & 0.002 & 0.001 & 0.002 & 0.002 \\
                    9 & 0.002 & 0.000 & 0.000 & 0.000 & 0.000 & 0.001 & 0.001 & 0.000 & 0.000 & 0.000 \\
                    10 & 0.000 & 0.006 & 0.004 & 0.000 & 0.002 & 0.001 & 0.002 & 0.001 & 0.002 & 0.003 \\
                    11 & 0.001 & 0.000 & 0.000 & 0.001 & 0.000 & 0.000 & 0.000 & 0.000 & 0.000 & 0.000 \\
                    12 & 0.000 & 0.004 & 0.002 & 0.004 & 0.000 & 0.005 & 0.004 & 0.004 & 0.004 & 0.004 \\
                    13 & 0.006 & 0.004 & 0.002 & 0.002 & 0.001 & 0.001 & 0.000 & 0.000 & 0.001 & 0.004 \\
                    14 & 0.002 & 0.002 & 0.004 & 0.001 & 0.000 & 0.001 & 0.001 & 0.000 & 0.001 & 0.001 \\
                    15 & 0.001 & 0.004 & 0.003 & 0.002 & 0.001 & 0.000 & 0.000 & 0.000 & 0.000 & 0.000 \\
                    16 & 0.004 & 0.003 & 0.004 & 0.002 & 0.002 & 0.002 & 0.002 & 0.000 & 0.001 & 0.000 \\
                    \bottomrule
                \end{tabular}
            \end{table}
           
            \begin{table}[h]
            \captionof{table}{Median relative gap per instance (Objective 2, , Scenarios 1--10)}
            \label{tab:gap_obj2_1}
            \setlength{\tabcolsep}{3pt}
            \renewcommand{\arraystretch}{0.95}  
                \begin{tabular}{l*{10}{S}}
                    \toprule
                    Scenario & {1} & {2} & {3} & {4} & {5} & {6} & {7} & {8} & {9} & {10} \\
                    \midrule        
                    
                    1 & 0.00 & 0.03 & {0.02\textsuperscript{*}} & 0.02 & 0.01 & 0.00 & 0.00 & 0.00 & 0.01 & 0.01 \\
                    2 & 0.00 & 0.00 & 0.00 & 0.01 & {0.01\textsuperscript{*}} & 0.01 & {0.01\textsuperscript{*}} & 0.01 & 0.00 & {0.00\textsuperscript{*}} \\
                    3 & 0.00 & {0.03\textsuperscript{*}} & {0.03\textsuperscript{*}} & {0.02\textsuperscript{*}} & 0.01 & 0.01 & {0.01\textsuperscript{*}} & 0.00 & 0.01 & 0.01\\
                    4 & 0.00 & 0.00 & 0.01 & 0.02 & {0.03\textsuperscript{*}} & {0.05\textsuperscript{*}} & {0.03\textsuperscript{*}} & 0.01 & {0.01\textsuperscript{*}} & {0.01\textsuperscript{*}} \\
                    5 & 0.00 & 0.03 & {0.03\textsuperscript{*}} & {0.02\textsuperscript{*}} & 0.01 & 0.02 & 0.01 & 0.01 & 0.01 & 0.01 \\
                    6 & 0.00 & 0.00 & 0.00 & 0.00 & {0.01\textsuperscript{*}} & {0.01\textsuperscript{*}} & {0.04\textsuperscript{*}} & {0.03\textsuperscript{*}} & {0.01\textsuperscript{*}} & {0.01\textsuperscript{*}} \\
                    7 & 0.00 & {0.05\textsuperscript{*}} & {0.07\textsuperscript{*}} & {0.02\textsuperscript{*}} & 0.03 & {0.00\textsuperscript{*}} & 0.01 & 0.00 & 0.01 & 0.01 \\
                    8 & 0.00 & 0.00 & 0.00 & {0.00\textsuperscript{*}} & {0.03\textsuperscript{*}} & 0.03 & {0.02\textsuperscript{*}} & {0.01\textsuperscript{*}} & {0.04\textsuperscript{*}} & {0.01\textsuperscript{*}} \\
                    9 & 0.00 & 0.03 & 0.02 & 0.01 & {0.01\textsuperscript{*}} & 0.01 & 0.01 & {0.00\textsuperscript{*}} & 0.01 & 0.01 \\
                    10 & 0.00 & 0.00 & 0.00 & 0.01 & {0.02\textsuperscript{*}} & {0.01\textsuperscript{*}} & 0.00 & {0.00\textsuperscript{*}} & {0.01\textsuperscript{*}} & {0.01\textsuperscript{*}}\\
                    11 & 0.01 & {0.03\textsuperscript{*}} & 0.01 & {0.01\textsuperscript{*}} & 0.01 & 0.01 & 0.01 & 0.00 & 0.01 & 0.01\\
                    12 & 0.00 & 0.00 & 0.00 & 0.00 & {0.02\textsuperscript{*}} & {0.02\textsuperscript{*}} & 0.01 & {0.01\textsuperscript{*}} & {0.01\textsuperscript{*}} & {0.01\textsuperscript{*}} \\
                    13 & 0.00 & 0.03 & {0.01\textsuperscript{*}} & 0.02 & 0.02 & 0.01 & {0.01\textsuperscript{*}} & 0.01 & 0.01 & 0.01 \\
                    14 & 0.00 & 0.00 & 0.00 & 0.02 & {0.03\textsuperscript{*}} & {0.03\textsuperscript{*}} & {0.01\textsuperscript{*}} & {0.02\textsuperscript{*}} & {0.01\textsuperscript{*}} & {0.01\textsuperscript{*}} \\
                    15 & 0.01 & {0.04\textsuperscript{*}} & 0.03 & 0.03 & {0.01\textsuperscript{*}} & {0.01\textsuperscript{*}} & {0.01\textsuperscript{*}} & {0.01\textsuperscript{*}} & 0.01 & 0.02 \\
                    16 & 0.00 & 0.00 & 0.00 & {0.00\textsuperscript{*}} & {0.03\textsuperscript{*}} & {0.04\textsuperscript{*}} & {0.05\textsuperscript{*}} & {0.02\textsuperscript{*}} & {0.00\textsuperscript{*}} & {0.01\textsuperscript{*}} \\
                    \bottomrule
                    \multicolumn{11}{l}{{*} Solving time limit reached (bound used).}\\
                    \end{tabular}                
            \end{table}
            
            \begin{table}[h]
            \captionof{table}{Median relative gap per instance (Objective 2, , Scenarios 11--20)}
            \label{tab:gap_obj2_2}
            \setlength{\tabcolsep}{3pt}
            \renewcommand{\arraystretch}{0.95}  
                \begin{tabular}{l*{10}{S}}
                    \toprule
                    Scenario & {11} & {12} & {13} & {14} & {15} & {16} & {17} & {18} & {19} & {20} \\
                    \midrule            
                    
                    1 & 0.01 & 0.01 & 0.00 & 0.01 & 0.01 & 0.00 & 0.01 & 0.01 & 0.01 & 0.01 \\
                    2 & 0.00 & 0.00 & 0.00 & 0.00 & 0.01 & 0.00 & 0.00 & 0.01 & 0.00 & 0.00 \\
                    3 & 0.01 & 0.01 & 0.01 & 0.01 & 0.01 & 0.00 & 0.00 & 0.00 & 0.00 & 0.00 \\
                    4 & {0.01\textsuperscript{*}} & {0.00\textsuperscript{*}} & {0.01\textsuperscript{*}} & {0.00\textsuperscript{*}} & 0.01 & 0.01 & 0.00 & 0.00 & 0.00 & 0.00 \\
                    5 & 0.01 & 0.01 & 0.01 & 0.01 & 0.01 & 0.01 & 0.01 & 0.01 & 0.00 & 0.00 \\
                    6 & {0.01\textsuperscript{*}} & {0.00\textsuperscript{*}} & {0.00\textsuperscript{*}} & 0.01 & 0.00 & 0.00 & 0.00 & 0.00 & 0.00 & 0.00 \\
                    7 & 0.01 & 0.00 & 0.00 & 0.00 & 0.00 & 0.00 & 0.00 & 0.00 & 0.00 & 0.00 \\
                    8 & {0.00\textsuperscript{*}} & {0.00\textsuperscript{*}} & {0.00\textsuperscript{*}} & {0.00\textsuperscript{*}} & {0.00\textsuperscript{*}} & 0.00 & 0.00 & 0.00 & 0.00 & 0.00 \\
                    9 & 0.01 & 0.01 & 0.01 & 0.01 & 0.01 & 0.02 & 0.01 & 0.01 & 0.01 & 0.02 \\
                    10 & {0.01\textsuperscript{*}} & {0.00\textsuperscript{*}} & 0.01 & 0.01 & 0.01 & 0.01 & 0.01 & 0.01 & 0.01 & 0.01 \\
                    11 & 0.01 & 0.02 & 0.02 & 0.02 & 0.02 & 0.02 & 0.02 & 0.02 & 0.03 & 0.03 \\
                    12 & 0.01 & 0.01 & 0.01 & {0.01\textsuperscript{*}} & 0.01 & 0.01 & 0.00 & 0.01 & 0.01 & 0.01 \\
                    13 & 0.01 & 0.00 & 0.00 & 0.00 & 0.00 & 0.00 & 0.00 & 0.00 & 0.00 & 0.00 \\
                    14 & {0.01\textsuperscript{*}} & {0.01\textsuperscript{*}} & {0.01\textsuperscript{*}} & {0.01\textsuperscript{*}} & 0.00 & 0.00 & 0.00 & 0.00 & 0.00 & 0.00 \\
                    15 & 0.01 & 0.01 & 0.01 & 0.02 & 0.02 & 0.02 & 0.02 & 0.03 & 0.04 & 0.02 \\
                    16 & {0.01\textsuperscript{*}} & {0.01\textsuperscript{*}} & {0.01\textsuperscript{*}} & {0.00\textsuperscript{*}} & {0.00\textsuperscript{*}} & 0.00 & {0.00\textsuperscript{*}} & 0.00 & 0.00 & 0.00 \\
                    \bottomrule
                    \multicolumn{11}{l}{{*} Solving time limit reached (bound used).}\\
                    \end{tabular}                
            \end{table}

            \begin{table}[h]
            \captionof{table}{Median relative gap per instance (Objective 3, , Scenarios 1--10)}
            \label{tab:gap_obj3_1}
            \setlength{\tabcolsep}{3pt}
            \renewcommand{\arraystretch}{0.95}  
                \begin{tabular}{l*{10}{S}}
                    \toprule
                    Scenario & {1} & {2} & {3} & {4} & {5} & {6} & {7} & {8} & {9} & {10} \\
                    \midrule        
                    1 & 1.29 & 0.63 & {1.12\textsuperscript{*}} & {0.87\textsuperscript{*}} & 0.19 & 0.04 & {0.04\textsuperscript{*}} & {0.05\textsuperscript{*}} & {0.05\textsuperscript{*}} & 0.15 \\
                    2 & 0.95 & 0.01 & 1.68 & 0.53 & {0.58\textsuperscript{*}} & {0.85\textsuperscript{*}} & {2.88\textsuperscript{*}} & {0.09\textsuperscript{*}} & {0.22\textsuperscript{*}} & {0.19\textsuperscript{*}} \\
                    3 & 3.06 & {1.80\textsuperscript{*}} & {1.15\textsuperscript{*}} & {1.20\textsuperscript{*}} & {0.53\textsuperscript{*}} & 0.25 & {0.21\textsuperscript{*}} & 0.08 & 0.10 & 0.17 \\
                    4 & 0.95 & 0.01 & 0.41 & 1.26 & {3.96\textsuperscript{*}} & {2.03\textsuperscript{*}} & {1.05\textsuperscript{*}} & {1.44\textsuperscript{*}} & {0.58\textsuperscript{*}} & {0.36\textsuperscript{*}} \\
                    5 & 1.29 & 0.59 & {1.68\textsuperscript{*}} & {2.99\textsuperscript{*}} & {0.55\textsuperscript{*}} & {0.50\textsuperscript{*}} & 0.32 & 0.44 & 0.46 & 0.39 \\
                    6 & 0.95 & 0.01 & 0.02 & 0.34 & {2.06\textsuperscript{*}} & {2.80\textsuperscript{*}} & {9.28\textsuperscript{*}} & {6.17\textsuperscript{*}} & {0.76\textsuperscript{*}} & {0.83\textsuperscript{*}} \\
                    7 & 3.06 & 0.99 & {10.71\textsuperscript{*}} & {3.74\textsuperscript{*}} & 2.56 & 0.42 & 1.02 & 0.65 & 0.67 & {0.31\textsuperscript{*}} \\
                    8 & 0.95 & 4.52 & 0.01 & 0.12 & {0.94\textsuperscript{*}} & {31.46\textsuperscript{*}} & {19.07\textsuperscript{*}} & {14.43\textsuperscript{*}} & {1.38\textsuperscript{*}} & {1.52\textsuperscript{*}} \\
                    9 & 3.70 & {0.43\textsuperscript{*}} & {0.74\textsuperscript{*}} & {0.41\textsuperscript{*}} & 0.13 & {0.14\textsuperscript{*}} & 0.04 & {0.01\textsuperscript{*}} & {0.14\textsuperscript{*}} & 0.09 \\
                    10 & 0.24 & 0.01 & 0.01 & 0.64 & {1.79\textsuperscript{*}} & {3.76\textsuperscript{*}} & {0.38\textsuperscript{*}} & {1.03\textsuperscript{*}} & {0.24\textsuperscript{*}} & {0.18\textsuperscript{*}} \\
                    11 & 4.90 & {1.50\textsuperscript{*}} & {0.54\textsuperscript{*}} & 0.27 & {0.22\textsuperscript{*}} & {0.01\textsuperscript{*}} & 0.00 & {0.23\textsuperscript{*}} & {0.09\textsuperscript{*}} & {0.08\textsuperscript{*}} \\
                    12 & 1.07 & 0.87 & 1.39 & {1.47\textsuperscript{*}} & {3.72\textsuperscript{*}} & {4.34\textsuperscript{*}} & {0.24\textsuperscript{*}} & {0.30\textsuperscript{*}} & {0.32\textsuperscript{*}} & {0.67\textsuperscript{*}} \\
                    13 & 3.70 & 0.46 & {0.68\textsuperscript{*}} & {0.73\textsuperscript{*}} & 0.90 & {0.99\textsuperscript{*}} & {1.30\textsuperscript{*}} & {0.31\textsuperscript{*}} & {0.34\textsuperscript{*}} & {0.46\textsuperscript{*}} \\
                    14 & 0.24 & 0.01 & 0.03 & {0.82\textsuperscript{*}} & {3.59\textsuperscript{*}} & {12.61\textsuperscript{*}} & {4.52\textsuperscript{*}} & {0.93\textsuperscript{*}} & {0.52\textsuperscript{*}} & {0.80\textsuperscript{*}} \\
                    15 & 4.90 & {1.39\textsuperscript{*}} & {0.66\textsuperscript{*}} & {0.70\textsuperscript{*}} & {1.90\textsuperscript{*}} & {1.79\textsuperscript{*}} & {0.62\textsuperscript{*}} & {0.63\textsuperscript{*}} & {0.37\textsuperscript{*}} & {0.31\textsuperscript{*}} \\
                    16 & 1.07 & 0.58 & 0.01 & {1.20\textsuperscript{*}} & {25.42\textsuperscript{*}} & {18.09\textsuperscript{*}} & {5.10\textsuperscript{*}} & {2.96\textsuperscript{*}} & {1.90\textsuperscript{*}} & {1.99\textsuperscript{*}} \\
                    \bottomrule
                    \multicolumn{11}{l}{{*} Solving time limit reached (bound used).}\\
                    \end{tabular}                
            \end{table}
            
            \begin{table}[h]
            \captionof{table}{Median relative gap per instance (Objective 3, , Scenarios 11--20)}
            \label{tab:gap_obj3_2}
            \setlength{\tabcolsep}{3pt}
            \renewcommand{\arraystretch}{0.95}  
                \begin{tabular}{l*{10}{S}}
                    \toprule
                    Scenario & {11} & {12} & {13} & {14} & {15} & {16} & {17} & {18} & {19} & {20} \\
                    \midrule
                    1 & 0.20 & 0.14 & 0.14 & 0.15 & 0.12 & 0.07 & 0.06 & 0.06 & 0.05 & 0.05 \\
                    2 & {0.27\textsuperscript{*}} & {0.16\textsuperscript{*}} & {0.16\textsuperscript{*}} & 0.15 & 0.22 & 0.22 & 0.17 & 0.17 & 0.17 & 0.16 \\
                    3 & 0.17 & 0.10 & 0.09 & 0.13 & 0.19 & 0.22 & 0.30 & 0.26 & 0.21 & 0.12 \\
                    4 & {0.20\textsuperscript{*}} & {0.05\textsuperscript{*}} & {0.19\textsuperscript{*}} & {0.17\textsuperscript{*}} & {0.10\textsuperscript{*}} & 0.16 & {0.26\textsuperscript{*}} & 0.07 & 0.05 & 0.11 \\
                    5 & 0.18 & {0.20\textsuperscript{*}} & 0.15 & 0.14 & 0.14 & 0.14 & 0.17 & 0.22 & 0.39 & 0.67 \\
                    6 & {0.22\textsuperscript{*}} & {0.19\textsuperscript{*}} & {0.16\textsuperscript{*}} & {0.20\textsuperscript{*}} & 0.24 & {0.23\textsuperscript{*}} & {0.18\textsuperscript{*}} & 0.20 & 0.22 & 0.25 \\
                    7 & {0.45\textsuperscript{*}} & 0.07 & 0.08 & 0.10 & 0.12 & 0.14 & 0.19 & 0.25 & 0.61 & 1.23 \\
                    8 & {0.68\textsuperscript{*}} & {0.21\textsuperscript{*}} & {0.26\textsuperscript{*}} & {0.29\textsuperscript{*}} & {0.32\textsuperscript{*}} & {0.42\textsuperscript{*}} & 0.44 & 0.41 & {0.39\textsuperscript{*}} & 0.40 \\
                    9 & 0.06 & 0.10 & 0.11 & 0.07 & 0.08 & 0.05 & 0.08 & 0.11 & 0.16 & 0.22 \\
                    10 & {0.17\textsuperscript{*}} & {0.32\textsuperscript{*}} & {0.09\textsuperscript{*}} & {0.12\textsuperscript{*}} & {0.16\textsuperscript{*}} & {0.08\textsuperscript{*}} & 0.11 & 0.14 & 0.10 & 0.12 \\
                    11 & 0.12 & 0.27 & 0.22 & 0.16 & 0.19 & 0.20 & 0.18 & 0.15 & 0.17 & 0.19 \\
                    12 & {0.15\textsuperscript{*}} & {0.18\textsuperscript{*}} & {0.13\textsuperscript{*}} & {0.07\textsuperscript{*}} & {0.11\textsuperscript{*}} & 0.16 & 0.01 & 0.03 & 0.01 & 0.12 \\
                    13 & 0.38 & {0.16\textsuperscript{*}} & 0.04 & 0.04 & 0.05 & 0.06 & 0.06 & 0.09 & 0.06 & 0.03 \\
                    14 & {0.75\textsuperscript{*}} & {0.26\textsuperscript{*}} & {0.27\textsuperscript{*}} & {0.11\textsuperscript{*}} & {0.04\textsuperscript{*}} & {0.06\textsuperscript{*}} & 0.11 & 0.21 & 0.23 & 0.13 \\
                    15 & {0.40\textsuperscript{*}} & 0.01 & 0.00 & 0.00 & 0.00 & 0.00 & 0.00 & 0.00 & 0.00 & 0.03 \\
                    16 & {0.37\textsuperscript{*}} & {0.30\textsuperscript{*}} & {0.02\textsuperscript{*}} & {0.00\textsuperscript{*}} & {0.00\textsuperscript{*}} & {0.00\textsuperscript{*}} & {0.00\textsuperscript{*}} & {0.01\textsuperscript{*}} & 0.04 & 0.06 \\
                    \bottomrule
                    \multicolumn{11}{l}{{*} Solving time limit reached (bound used).}\\
                    \end{tabular} 
            \end{table}

            \begin{table}[h]
            \captionof{table}{Median relative gap per instance (Objective 4, , Scenarios 1--10)}
            \label{tab:gap_obj4_1}
            \setlength{\tabcolsep}{3pt}
            \renewcommand{\arraystretch}{0.95}  
                \begin{tabular}{l*{10}{S}}
                    \toprule
                    Scenario & {1} & {2} & {3} & {4} & {5} & {6} & {7} & {8} & {9} & {10} \\
                    \midrule        
                    1 & {0.10\textsuperscript{*}} & {0.09\textsuperscript{*}} & {0.03\textsuperscript{*}} & {0.01\textsuperscript{*}} & {0.02\textsuperscript{*}} & {0.01\textsuperscript{*}} & {0.00\textsuperscript{*}} & {0.04\textsuperscript{*}} & {0.00\textsuperscript{*}} & {0.07\textsuperscript{*}} \\
                    2 & 10.13 & 1.32 & 0.26 & {0.13\textsuperscript{*}} & {0.14\textsuperscript{*}} & {0.08\textsuperscript{*}} & {0.05\textsuperscript{*}} & {0.01\textsuperscript{*}} & {0.03\textsuperscript{*}} & {0.06\textsuperscript{*}}\\
                    3 & {0.67\textsuperscript{*}} & {0.43\textsuperscript{*}} & {0.19\textsuperscript{*}} & {0.06\textsuperscript{*}} & {0.02\textsuperscript{*}} & 0.05 & {0.00\textsuperscript{*}} & {0.02\textsuperscript{*}} & {0.02\textsuperscript{*}} & {0.09\textsuperscript{*}}\\
                    4 & 10.58 & 2.74 & {0.71\textsuperscript{*}} & {0.32\textsuperscript{*}} & {0.44\textsuperscript{*}} & {0.13\textsuperscript{*}} & {0.08\textsuperscript{*}} & {0.06\textsuperscript{*}} & {0.13\textsuperscript{*}} & {0.07\textsuperscript{*}} \\
                    5 & {0.10\textsuperscript{*}} & {0.05\textsuperscript{*}} & {0.10\textsuperscript{*}} & {0.09\textsuperscript{*}} & 0.06 & {0.22\textsuperscript{*}} & 0.14 & {0.12\textsuperscript{*}} & 0.19 & {0.15\textsuperscript{*}} \\
                    6 & 10.13 & 1.88 & 0.39 & 0.22 & {0.12\textsuperscript{*}} & {0.12\textsuperscript{*}} & {0.16\textsuperscript{*}} & {0.17\textsuperscript{*}} & {0.11\textsuperscript{*}} & {0.25\textsuperscript{*}} \\
                    7 & {0.57\textsuperscript{*}} & {0.12\textsuperscript{*}} & {0.30\textsuperscript{*}} & {0.18\textsuperscript{*}} & 0.20 & 0.25 & 0.30 & {0.15\textsuperscript{*}} & 0.17 & {0.10\textsuperscript{*}} \\
                    8 & 10.58 & 2.58 & {0.62\textsuperscript{*}} & {0.49\textsuperscript{*}} & {0.50\textsuperscript{*}} & {0.22\textsuperscript{*}} & {0.10\textsuperscript{*}} & {0.12\textsuperscript{*}} & {0.23\textsuperscript{*}} & {0.25\textsuperscript{*}} \\
                    9 & {0.16\textsuperscript{*}} & {0.01\textsuperscript{*}} & {0.00\textsuperscript{*}} & {0.00\textsuperscript{*}} & {0.00\textsuperscript{*}} & {0.00\textsuperscript{*}} & 0.00 & {0.00\textsuperscript{*}} & {0.00\textsuperscript{*}} & 0.00 \\
                    10 & 0.11 & 0.16 & {0.09\textsuperscript{*}} & 0.09 & {0.11\textsuperscript{*}} & {0.09\textsuperscript{*}} & {0.00\textsuperscript{*}} & {0.07\textsuperscript{*}} & {0.05\textsuperscript{*}} & {0.04\textsuperscript{*}} \\
                    11 & {0.50\textsuperscript{*}} & {0.23\textsuperscript{*}} & {0.00\textsuperscript{*}} & 0.00 & {0.01\textsuperscript{*}} & {0.00\textsuperscript{*}} & {0.00\textsuperscript{*}} & 0.02 & {0.00\textsuperscript{*}} & 0.00 \\
                    12 & 2.80 & 0.79 & {0.24\textsuperscript{*}} & {0.19\textsuperscript{*}} & {0.23\textsuperscript{*}} & {0.22\textsuperscript{*}} & {0.20\textsuperscript{*}} & {0.11\textsuperscript{*}} & {0.20\textsuperscript{*}} & {0.19\textsuperscript{*}} \\
                    13 & {0.16\textsuperscript{*}} & {0.01\textsuperscript{*}} & {0.08\textsuperscript{*}} & {0.07\textsuperscript{*}} & {0.16\textsuperscript{*}} & {0.13\textsuperscript{*}} & {0.29\textsuperscript{*}} & {0.12\textsuperscript{*}} & {0.20\textsuperscript{*}} & {0.09\textsuperscript{*}} \\
                    14 & 0.11 & 0.24 & {0.15\textsuperscript{*}} & {0.13\textsuperscript{*}} & {0.17\textsuperscript{*}} & {0.13\textsuperscript{*}} & {0.04\textsuperscript{*}} & {0.21\textsuperscript{*}} & {0.27\textsuperscript{*}} & {0.21\textsuperscript{*}} \\
                    15 & {0.50\textsuperscript{*}} & {0.16\textsuperscript{*}} & {0.05\textsuperscript{*}} & {0.23\textsuperscript{*}} & {0.32\textsuperscript{*}} & {0.31\textsuperscript{*}} & {0.12\textsuperscript{*}} & {0.21\textsuperscript{*}} & {0.24\textsuperscript{*}} & {0.15\textsuperscript{*}} \\
                    16 & 2.80 & 1.67 & {0.38\textsuperscript{*}} & {0.28\textsuperscript{*}} & {0.29\textsuperscript{*}} & {0.10\textsuperscript{*}} & {0.29\textsuperscript{*}} & {0.33\textsuperscript{*}} & {0.51\textsuperscript{*}} & {0.29\textsuperscript{*}} \\   
                    \bottomrule
                    \multicolumn{11}{l}{{*} Solving time limit reached (bound used).}\\
                    \end{tabular}                
            \end{table}
            
            \begin{table}[h]
            \captionof{table}{Median relative gap per instance (Objective 4, , Scenarios 11--20)}
            \label{tab:gap_obj4_2}
            \setlength{\tabcolsep}{3pt}
            \renewcommand{\arraystretch}{0.95}  
                \begin{tabular}{l*{10}{S}}
                    \toprule
                    Scenario & {11} & {12} & {13} & {14} & {15} & {16} & {17} & {18} & {19} & {20} \\
                    \midrule
                    1 & 0.12 & 0.14 & 0.09 & 0.07 & 0.12 & 0.06 & 0.06 & 0.05 & 0.05 & 0.07 \\
                    2 & {0.02\textsuperscript{*}} & {0.02\textsuperscript{*}} & {0.07\textsuperscript{*}} & {0.01\textsuperscript{*}} & {0.07\textsuperscript{*}} & {0.08\textsuperscript{*}} & 0.06 & {0.07\textsuperscript{*}} & 0.11 & 0.11 \\
                    3 & {0.01\textsuperscript{*}} & 0.09 & 0.14 & 0.12 & 0.22 & 0.20 & 0.16 & 0.15 & 0.15 & 0.14 \\
                    4 & {0.15\textsuperscript{*}} & {0.07\textsuperscript{*}} & {0.13\textsuperscript{*}} & {0.09\textsuperscript{*}} & {0.09\textsuperscript{*}} & 0.07 & {0.12\textsuperscript{*}} & 0.20 & 0.13 & 0.18 \\
                    5 & {0.13\textsuperscript{*}} & {0.08\textsuperscript{*}} & 0.15 & 0.12 & 0.06 & 0.12 & 0.12 & 0.16 & 0.16 & 0.19 \\
                    6 & {0.23\textsuperscript{*}} & {0.20\textsuperscript{*}} & {0.13\textsuperscript{*}} & {0.12\textsuperscript{*}} & {0.02\textsuperscript{*}} & {0.06\textsuperscript{*}} & {0.06\textsuperscript{*}} & {0.09\textsuperscript{*}} & {0.11\textsuperscript{*}} & 0.14 \\
                    7 & {0.16\textsuperscript{*}} & 0.17 & 0.16 & 0.24 & 0.23 & 0.15 & 0.15 & 0.11 & 0.12 & 0.10 \\
                    8 & {0.18\textsuperscript{*}} & {0.13\textsuperscript{*}} & {0.08\textsuperscript{*}} & {0.11\textsuperscript{*}} & {0.16\textsuperscript{*}} & {0.18\textsuperscript{*}} & {0.22\textsuperscript{*}} & {0.21\textsuperscript{*}} & {0.31\textsuperscript{*}} & 0.21 \\
                    9 & 0.02 & 0.08 & 0.10 & 0.23 & 0.15 & 0.04 & 0.06 & 0.06 & 0.14 & 0.14 \\
                    10 & {0.03\textsuperscript{*}} & {0.05\textsuperscript{*}} & {0.03\textsuperscript{*}} & {0.00\textsuperscript{*}} & {0.04\textsuperscript{*}} & {0.03\textsuperscript{*}} & {0.06\textsuperscript{*}} & 0.15 & 0.13 & 0.10 \\
                    11 & 0.00 & 0.08 & 0.16 & 0.11 & 0.19 & 0.18 & 0.21 & 0.18 & 0.14 & 0.13 \\
                    12 & {0.01\textsuperscript{*}} & {0.06\textsuperscript{*}} & {0.07\textsuperscript{*}} & {0.09\textsuperscript{*}} & {0.07\textsuperscript{*}} & {0.12\textsuperscript{*}} & {0.10\textsuperscript{*}} & 0.16 & 0.13 & 0.20 \\
                    13 & {0.10\textsuperscript{*}} & {0.08\textsuperscript{*}} & 0.10 & 0.11 & 0.17 & 0.19 & 0.14 & 0.16 & 0.14 & 0.15 \\
                    14 & {0.10\textsuperscript{*}} & {0.13\textsuperscript{*}} & {0.17\textsuperscript{*}} & {0.18\textsuperscript{*}} & {0.05\textsuperscript{*}} & {0.10\textsuperscript{*}} & {0.09\textsuperscript{*}} & {0.07\textsuperscript{*}} & {0.10\textsuperscript{*}} & {0.07\textsuperscript{*}} \\
                    15 & {0.15\textsuperscript{*}} & {0.05\textsuperscript{*}} & 0.07 & 0.09 & 0.12 & 0.13 & 0.16 & 0.16 & 0.13 & 0.13 \\
                    16 & {0.25\textsuperscript{*}} & {0.19\textsuperscript{*}} & {0.18\textsuperscript{*}} & {0.03\textsuperscript{*}} & {0.06\textsuperscript{*}} & {0.03\textsuperscript{*}} & {0.07\textsuperscript{*}} & {0.10\textsuperscript{*}} & {0.20\textsuperscript{*}} & 0.24 \\   
                    \bottomrule
                    \multicolumn{11}{l}{{*} Solving time limit reached (bound used).}\\
                    \end{tabular} 
            \end{table}

    \subsection{Runtime analysis across scenarios and instances}
    \label{subsec:add_runtime_instances}

\begin{figure}[p]
        \centering
        
            \begin{subfigure}{0.48\textwidth}
            \includegraphics[width=\linewidth]{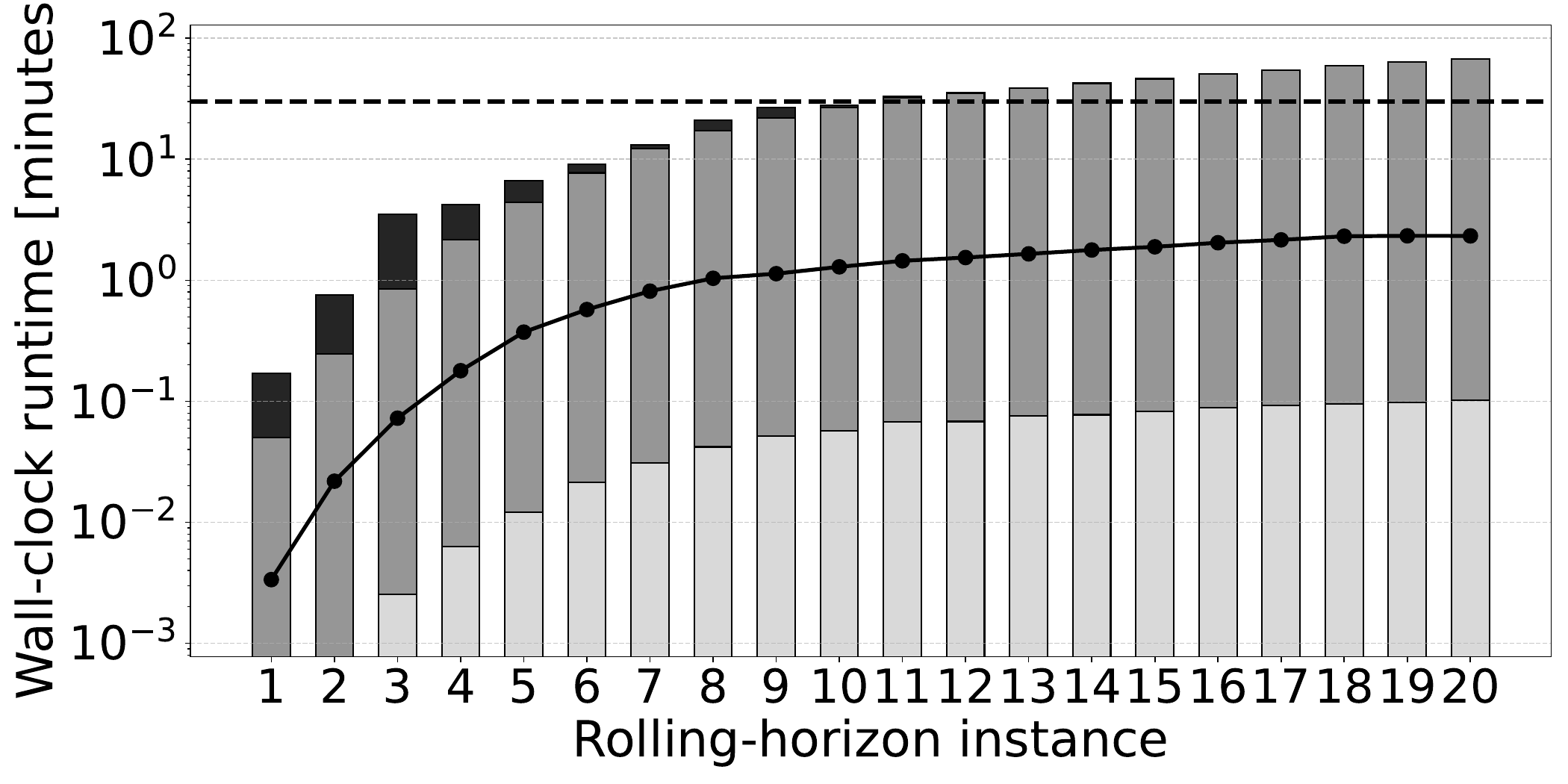}
            \caption{Scenario 1}
            \end{subfigure}
            \hfill
            \begin{subfigure}{0.48\textwidth}
            \includegraphics[width=\linewidth]{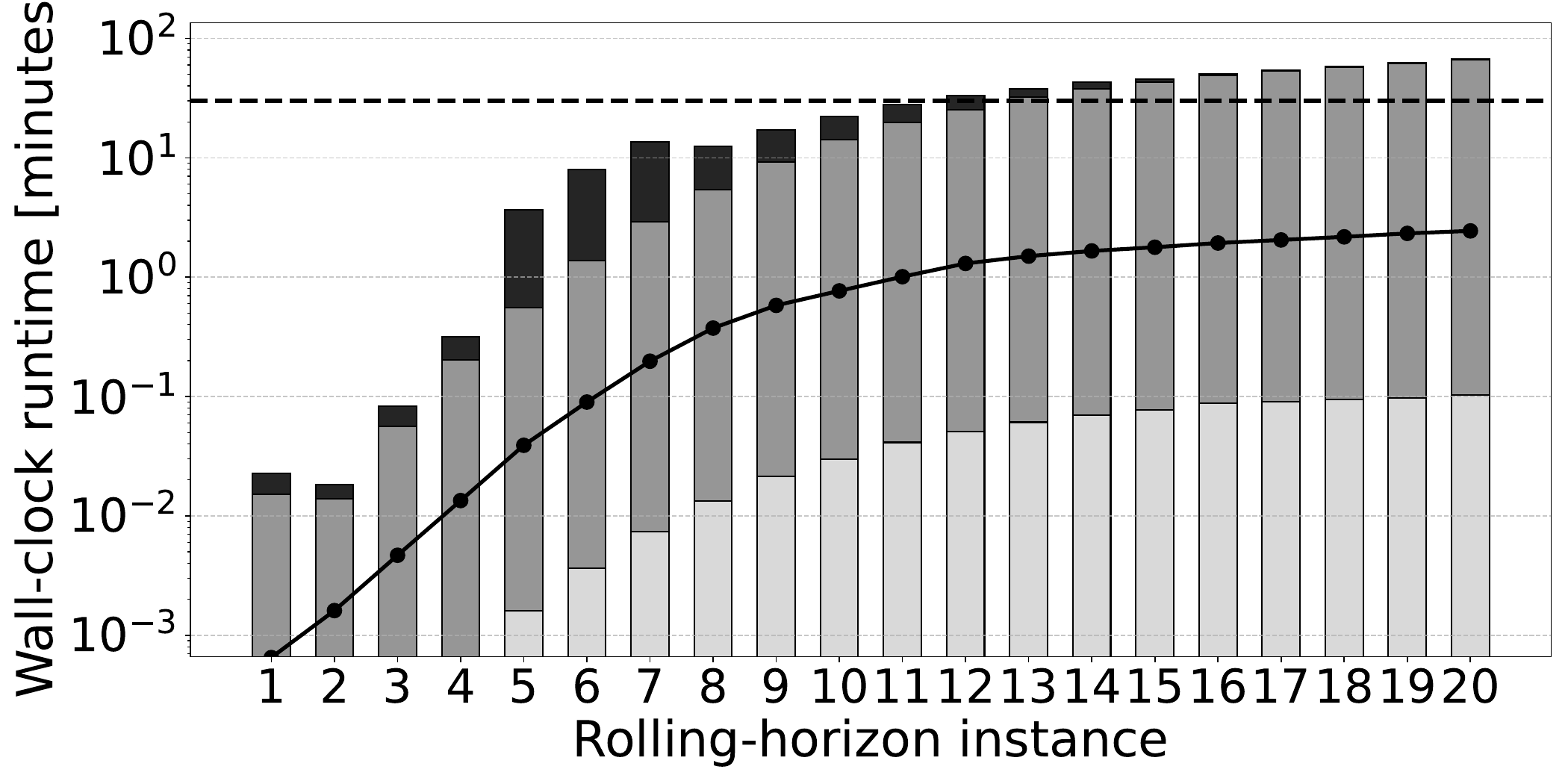}
            \caption{Scenario 2}
            \end{subfigure}
            
            \vspace{0.4cm}
            
            \begin{subfigure}{0.48\textwidth}
            \includegraphics[width=\linewidth]{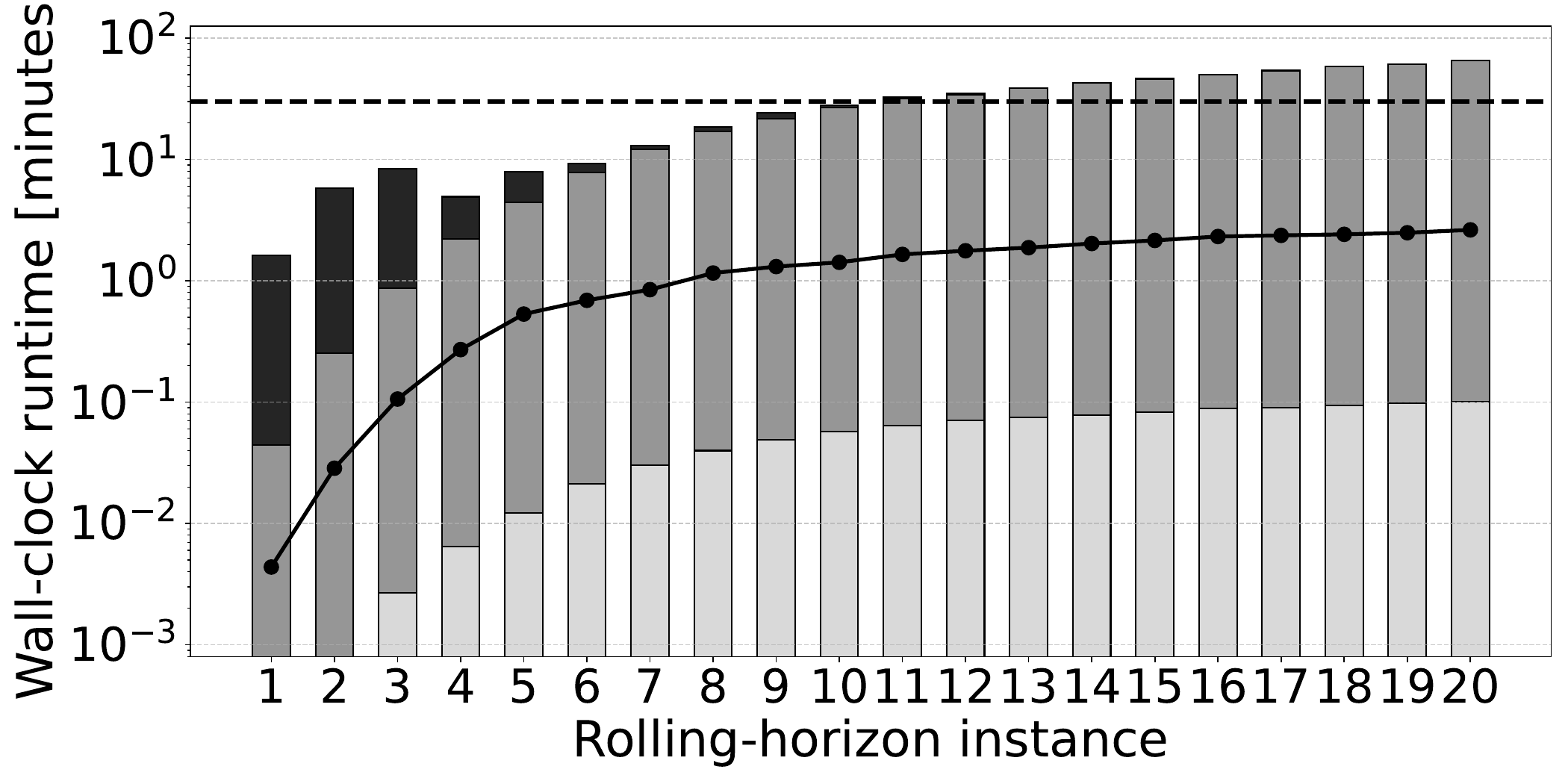}
            \caption{Scenario 3}
            \end{subfigure}
            \hfill
            \begin{subfigure}{0.48\textwidth}
            \includegraphics[width=\linewidth]{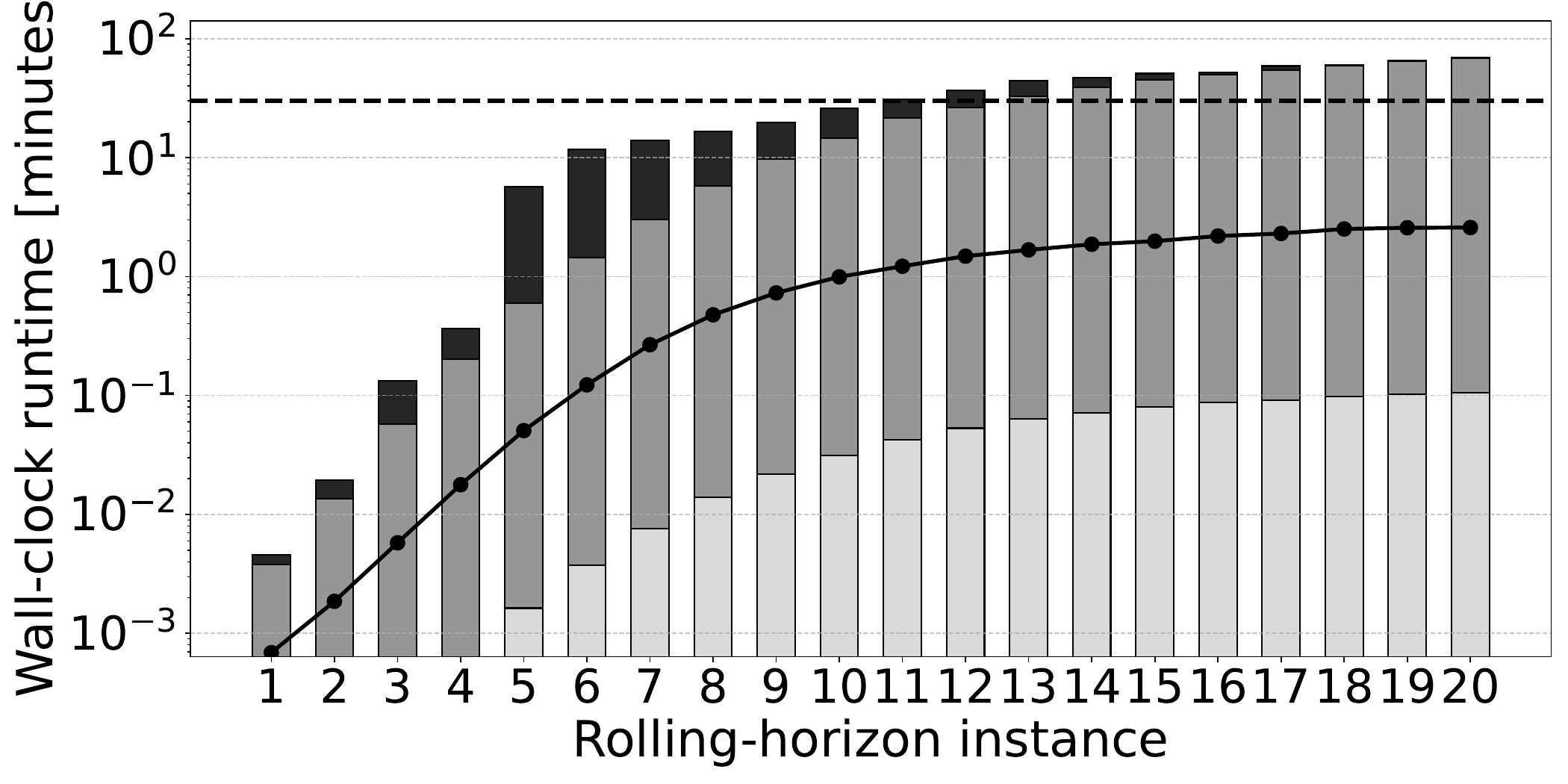}
            \caption{Scenario 4}
            \end{subfigure}
            
            \begin{subfigure}{0.7\textwidth}
            \includegraphics[width=\linewidth]{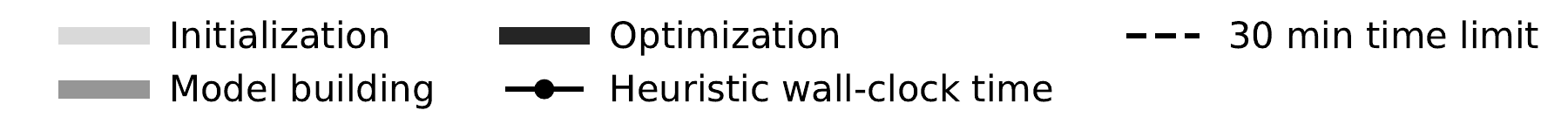}        
            \end{subfigure}
        
        \caption{Instance-level runtime breakdown for scenarios 1–4.
        Each plot reports the median runtime across the 10 stochastic replications
        for the 20 sequential instances.}
        \label{fig:runtime_instances_1_4}
        
        \end{figure}

\begin{figure}[p]
        \centering
        
            \begin{subfigure}{0.48\textwidth}
            \includegraphics[width=\linewidth]{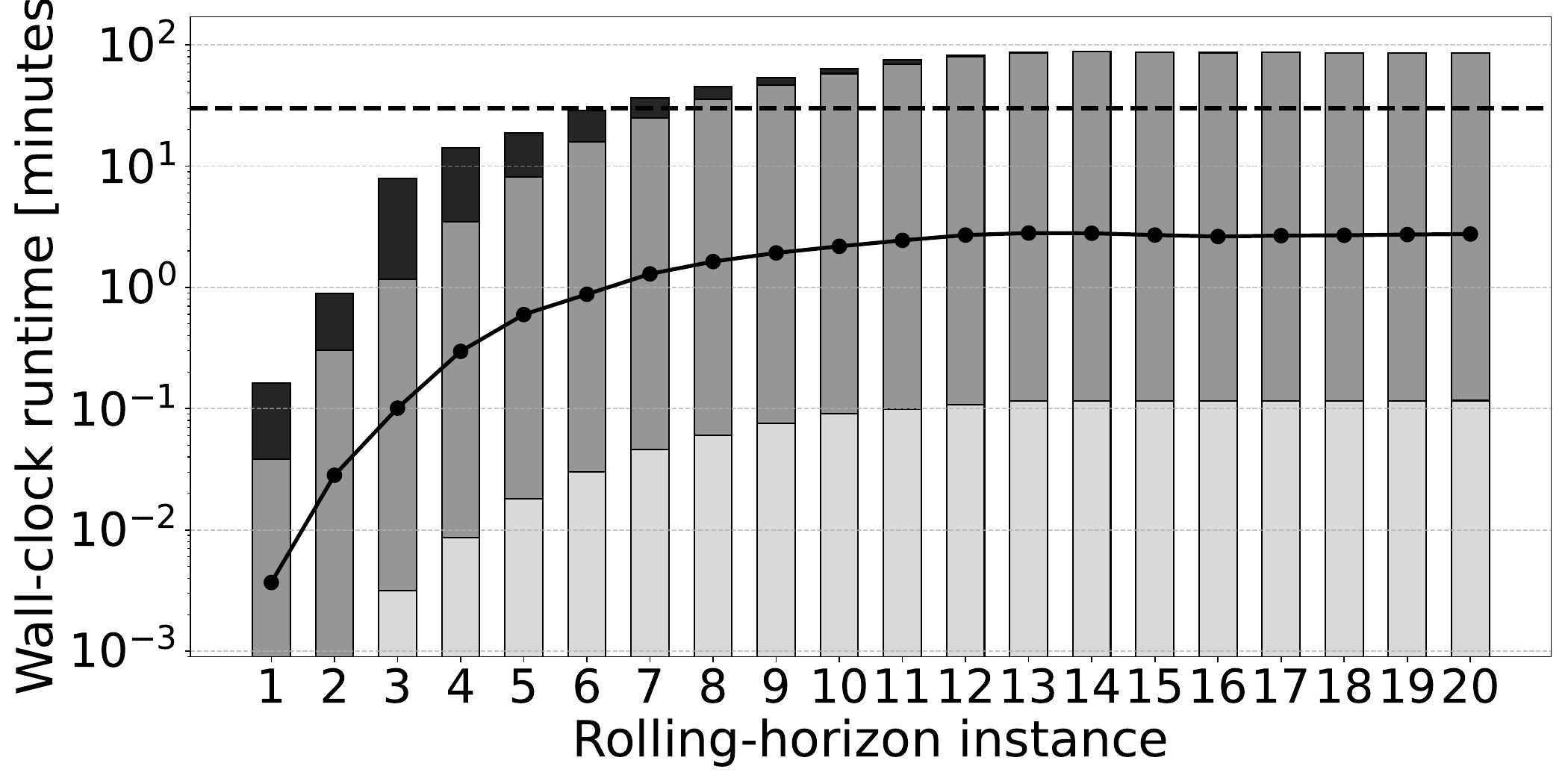}
            \caption{Scenario 5}
            \end{subfigure}
            \hfill
            \begin{subfigure}{0.48\textwidth}
            \includegraphics[width=\linewidth]{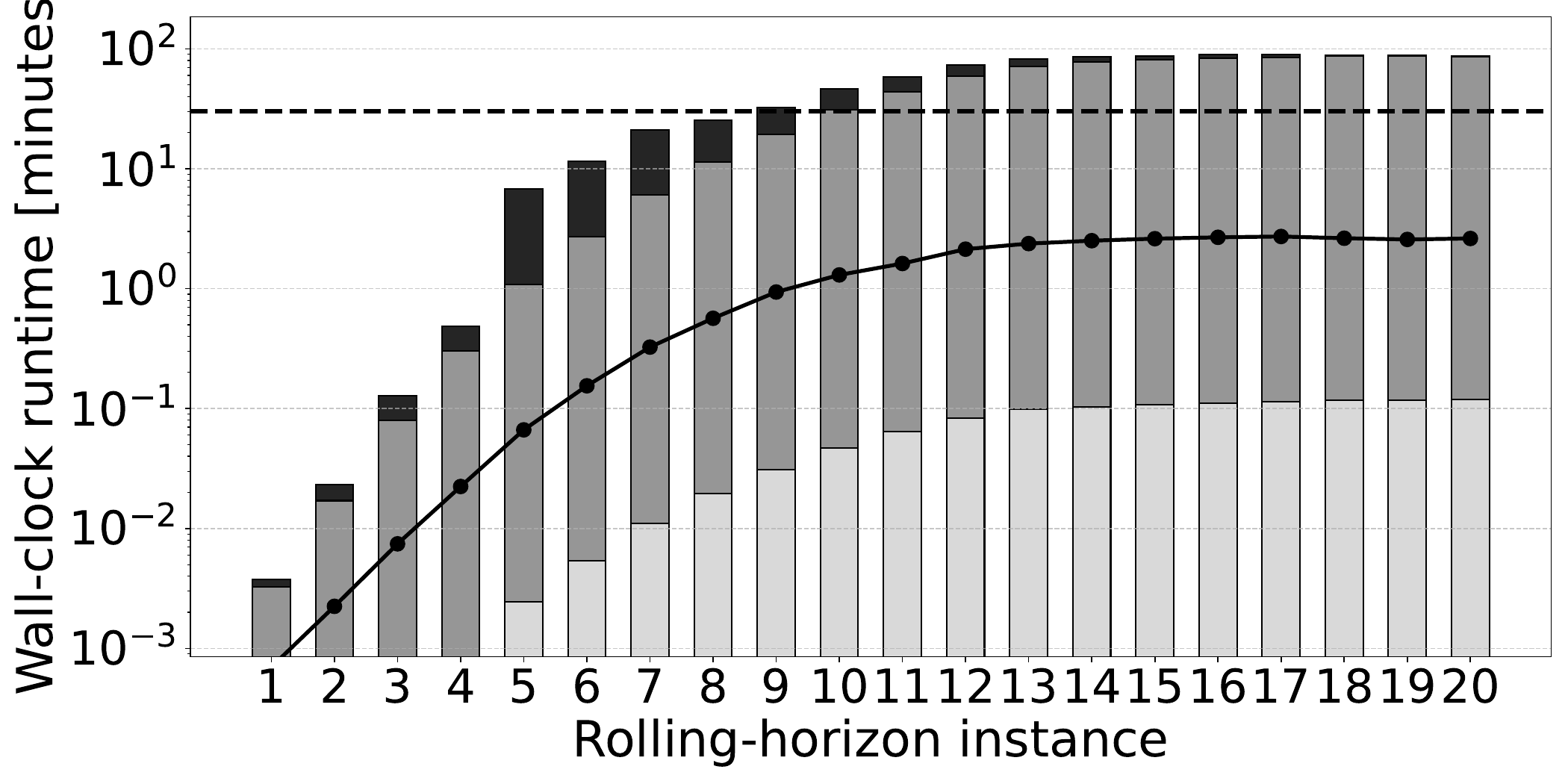}
            \caption{Scenario 6}
            \end{subfigure}
            
            \vspace{0.4cm}
            
            \begin{subfigure}{0.48\textwidth}
            \includegraphics[width=\linewidth]{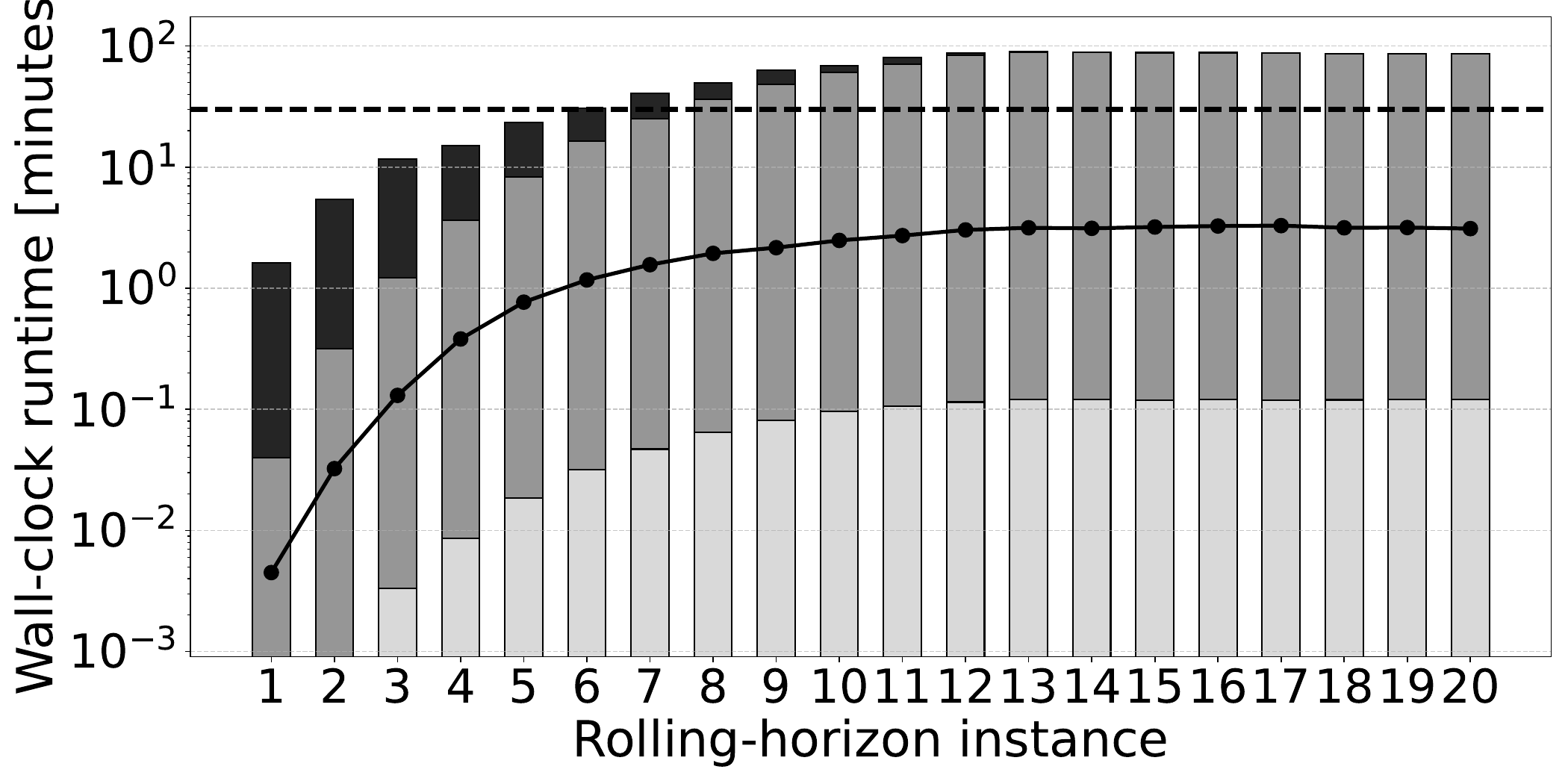}
            \caption{Scenario 7}
            \end{subfigure}
            \hfill
            \begin{subfigure}{0.48\textwidth}
            \includegraphics[width=\linewidth]{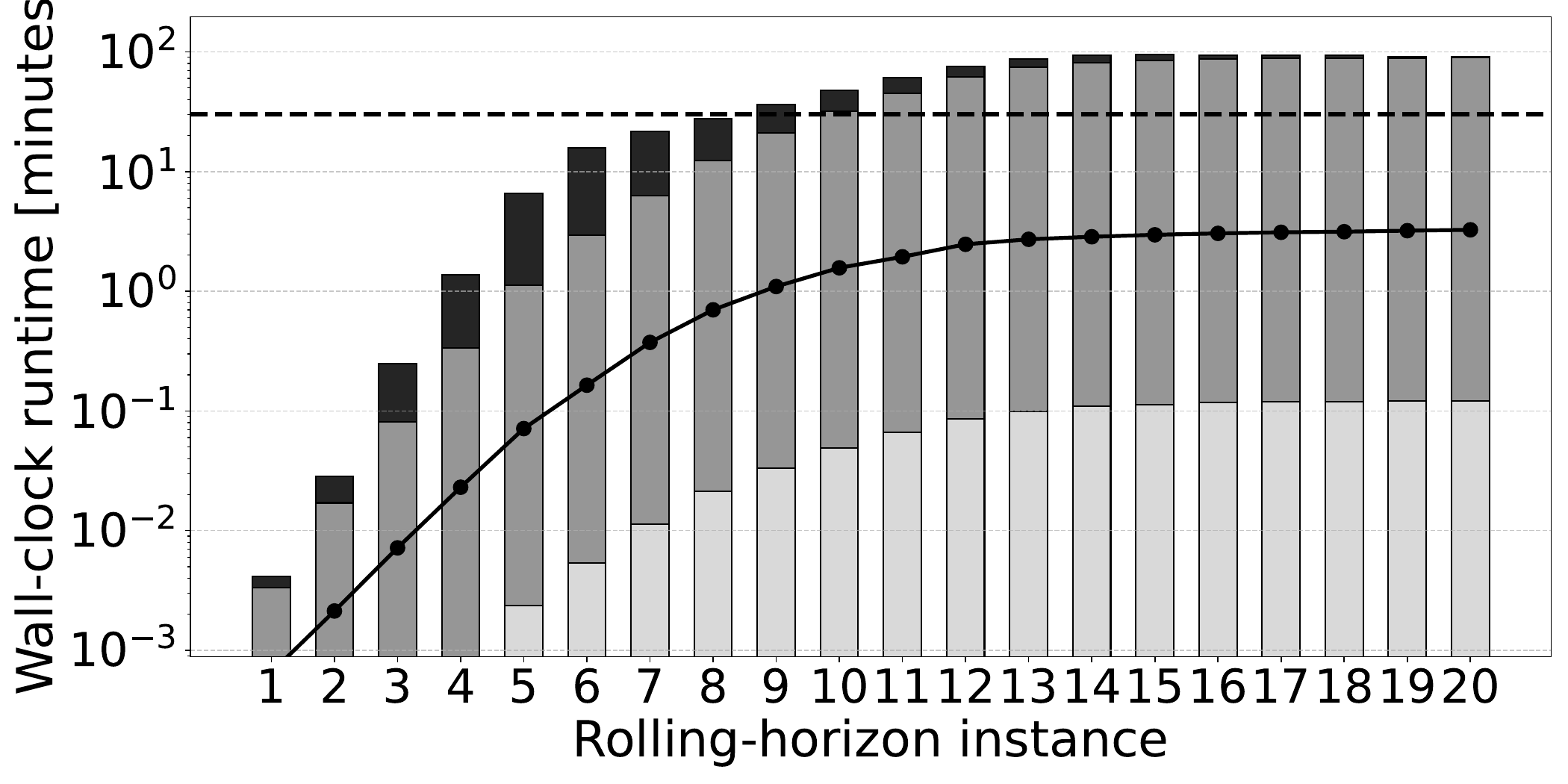}
            \caption{Scenario 8}
            \end{subfigure}

            \begin{subfigure}{0.7\textwidth}
            \includegraphics[width=\linewidth]{figures/runtime/legend.pdf}        
            \end{subfigure}
        
        \caption{Instance-level runtime breakdown for scenarios 5–8.
        Each plot reports the median runtime across the 10 stochastic replications
        for the 20 sequential instances.}
        \label{fig:runtime_instances_5_8}        
        \end{figure}

\begin{figure}[p]
        \centering
        
            \begin{subfigure}{0.48\textwidth}
            \includegraphics[width=\linewidth]{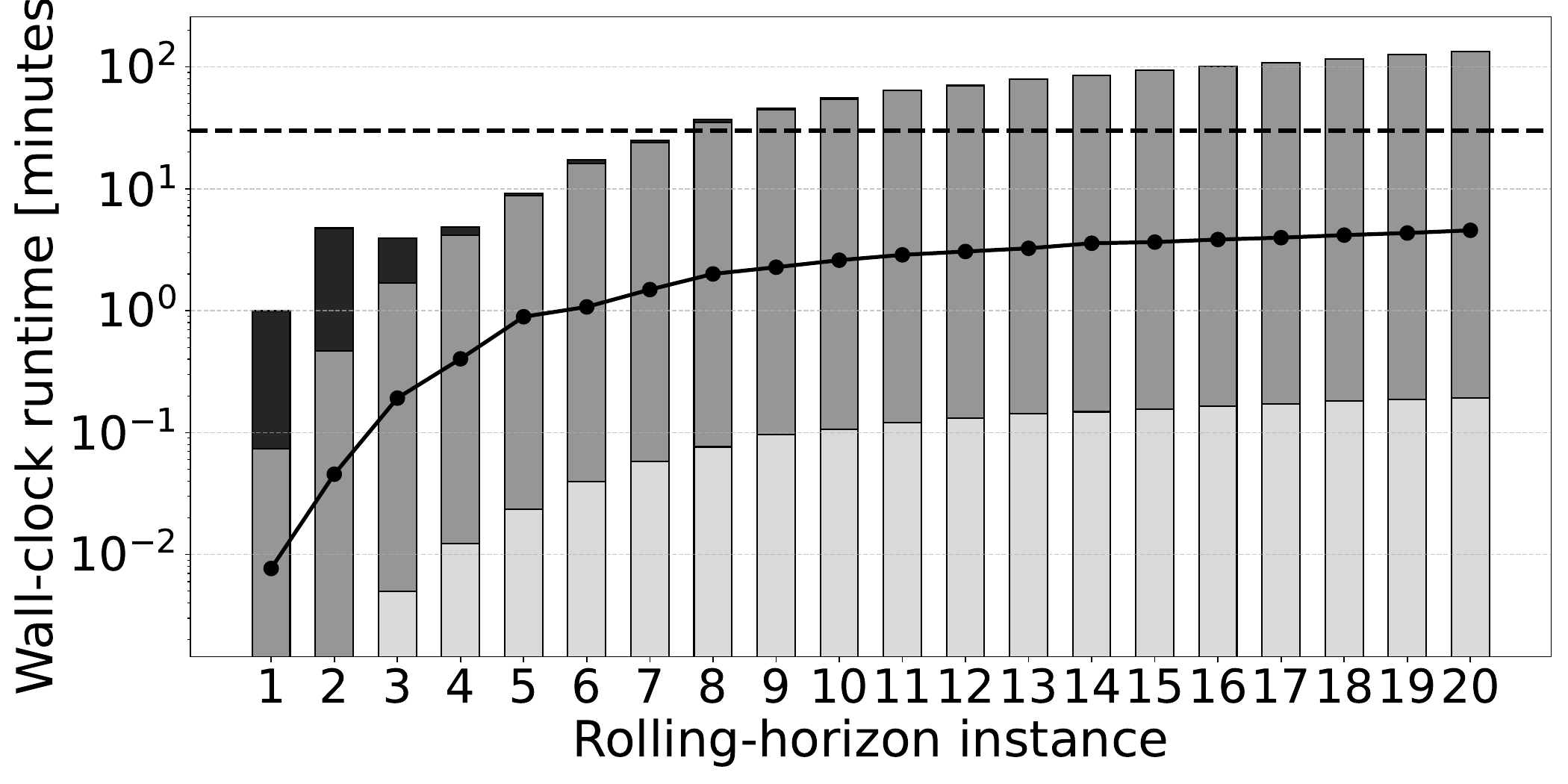}
            \caption{Scenario 9}
            \end{subfigure}
            \hfill
            \begin{subfigure}{0.48\textwidth}
            \includegraphics[width=\linewidth]{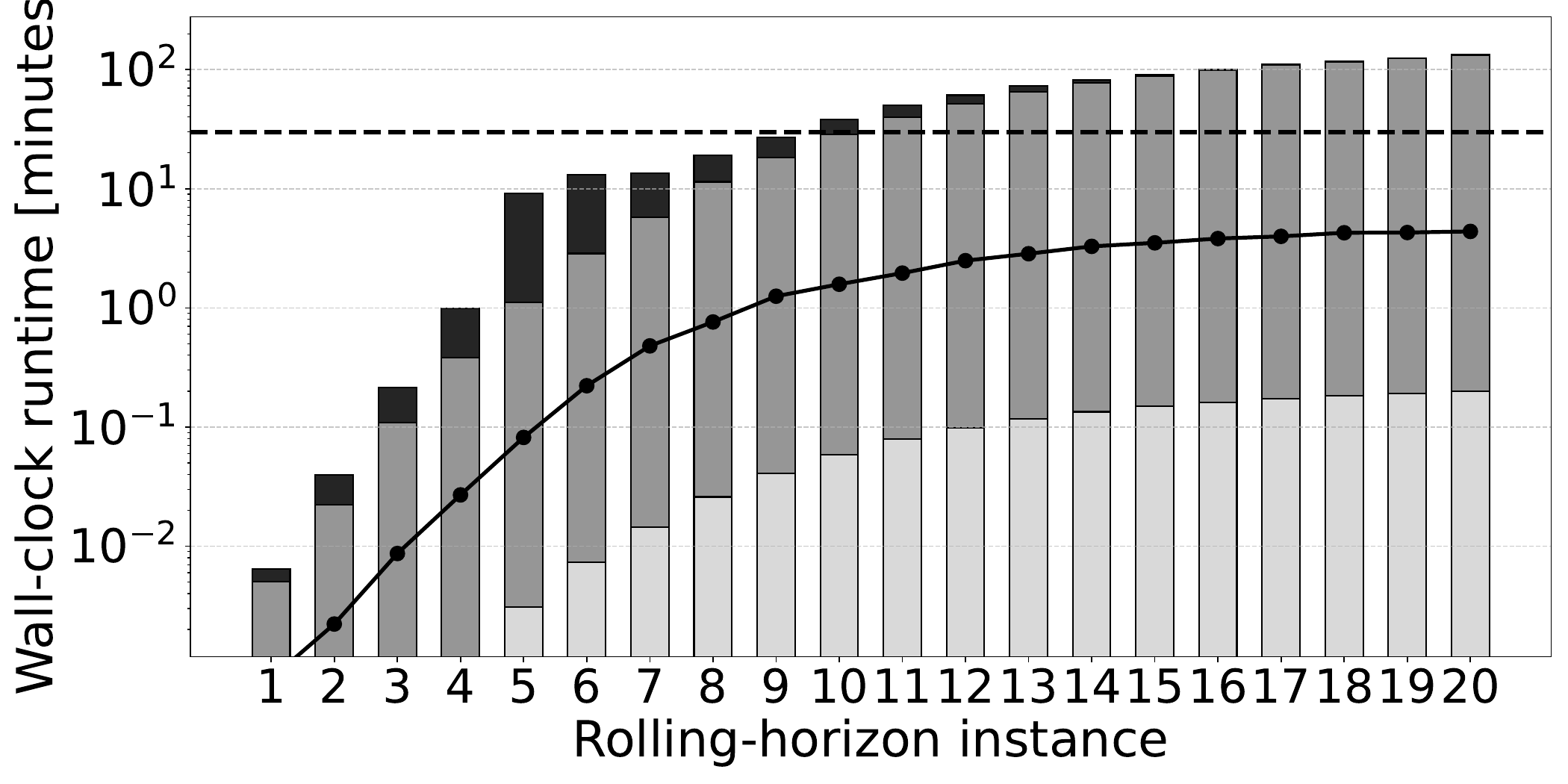}
            \caption{Scenario 10}
            \end{subfigure}
            
            \vspace{0.4cm}
            
            \begin{subfigure}{0.48\textwidth}
            \includegraphics[width=\linewidth]{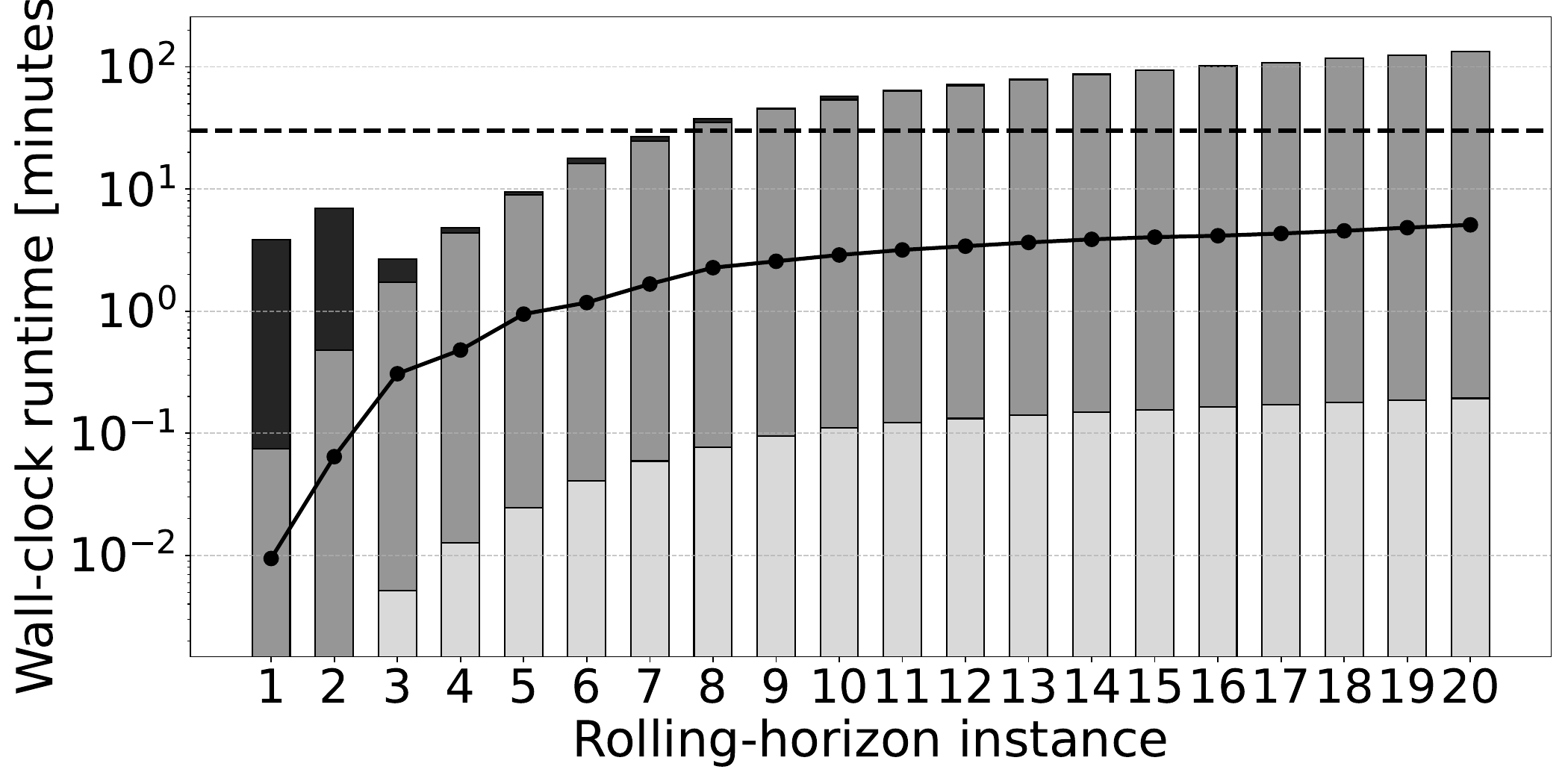}
            \caption{Scenario 11}
            \end{subfigure}
            \hfill
            \begin{subfigure}{0.48\textwidth}
            \includegraphics[width=\linewidth]{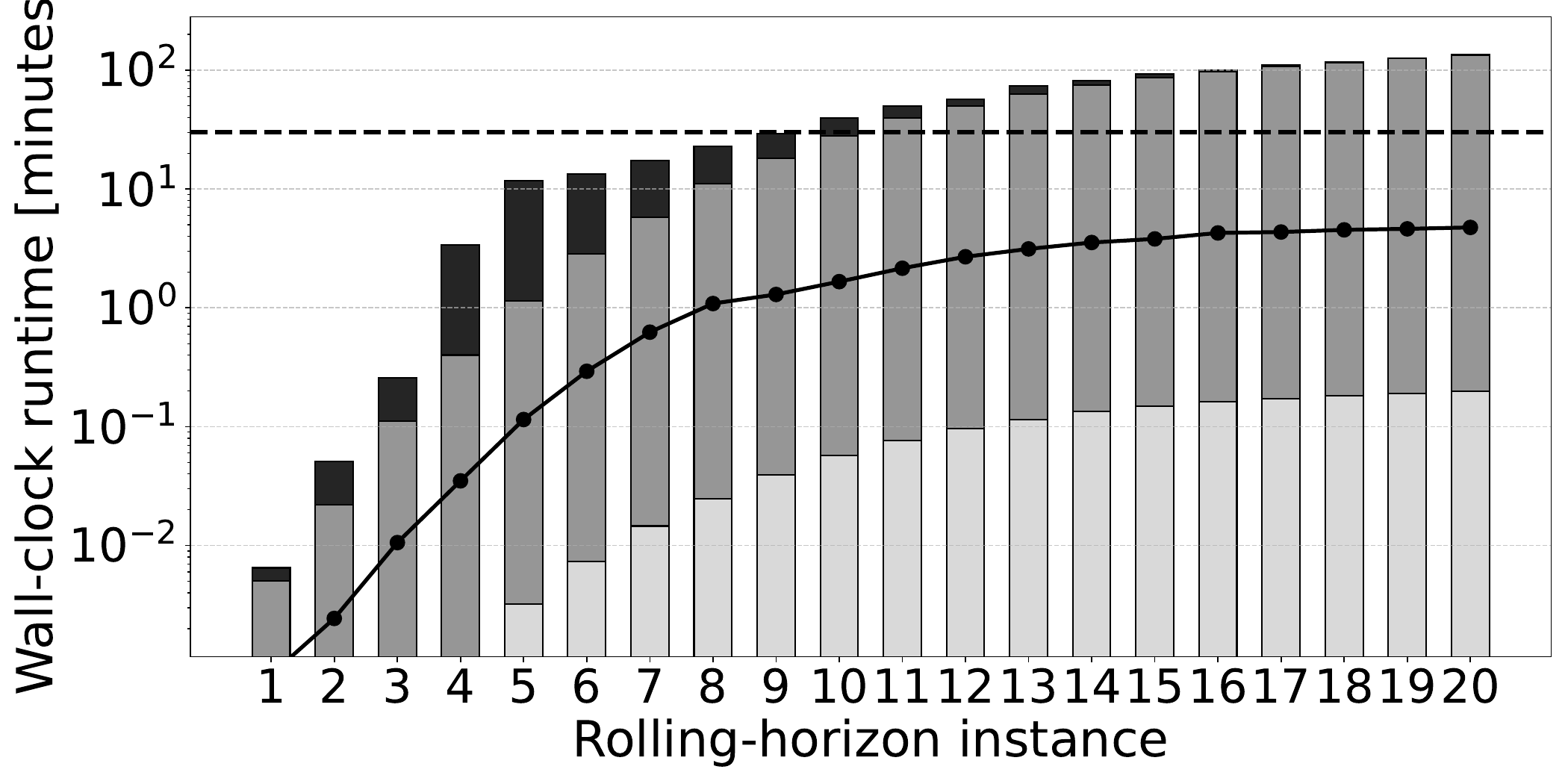}
            \caption{Scenario 12}
            \end{subfigure}

            \begin{subfigure}{0.7\textwidth}
            \includegraphics[width=\linewidth]{figures/runtime/legend.pdf}        
            \end{subfigure}
            
        \caption{Instance-level runtime breakdown for scenarios 9–12.
        Each plot reports the median runtime across the 10 stochastic replications
        for the 20 sequential instances.}
        \label{fig:runtime_instances_9_12}        
        \end{figure}

\begin{figure}[p]
        \centering
        
            \begin{subfigure}{0.48\textwidth}
            \includegraphics[width=\linewidth]{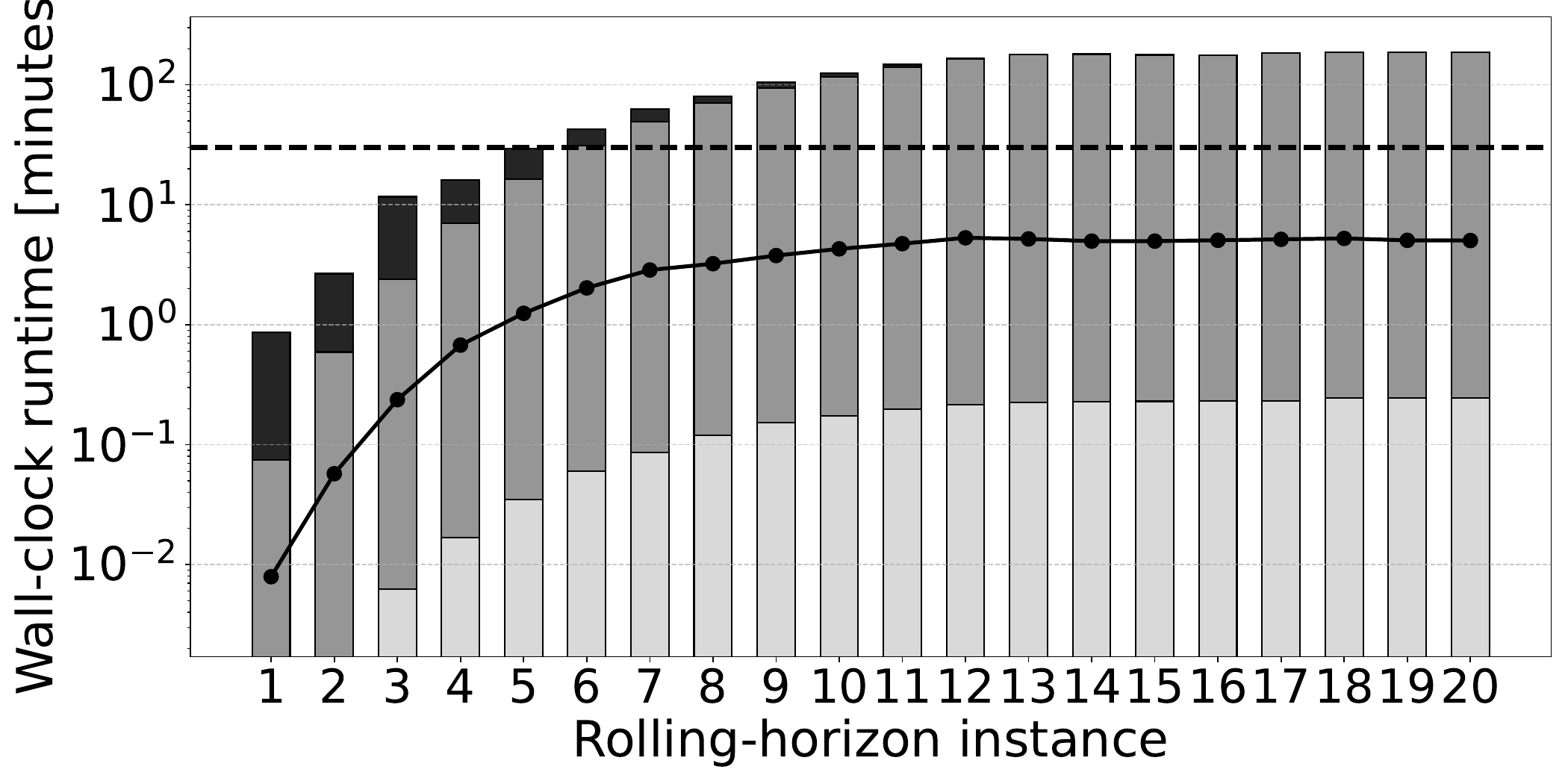}
            \caption{Scenario 13}
            \end{subfigure}
            \hfill
            \begin{subfigure}{0.48\textwidth}
            \includegraphics[width=\linewidth]{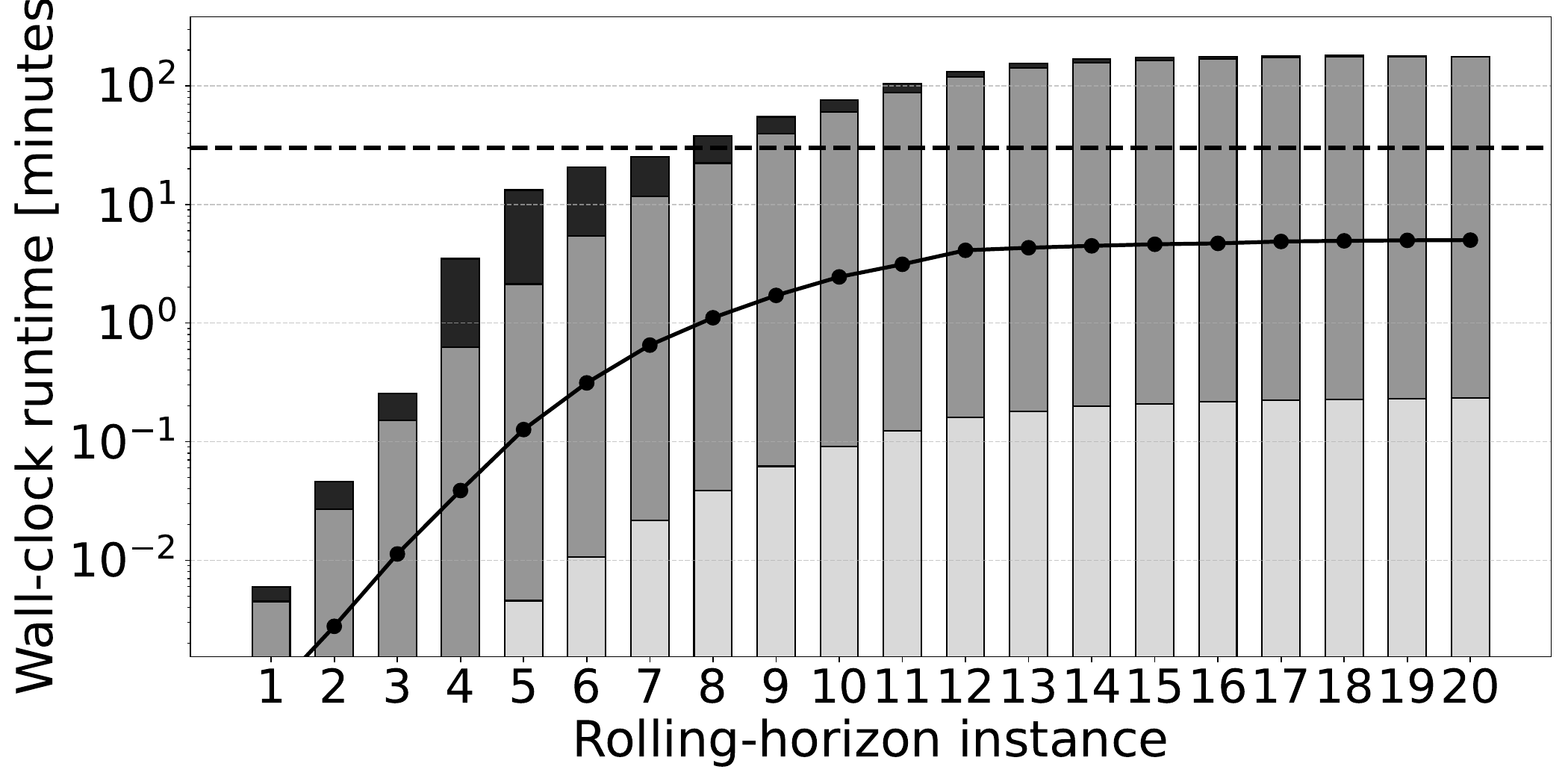}
            \caption{Scenario 14}
            \end{subfigure}
            
            \vspace{0.4cm}
            
            \begin{subfigure}{0.48\textwidth}
            \includegraphics[width=\linewidth]{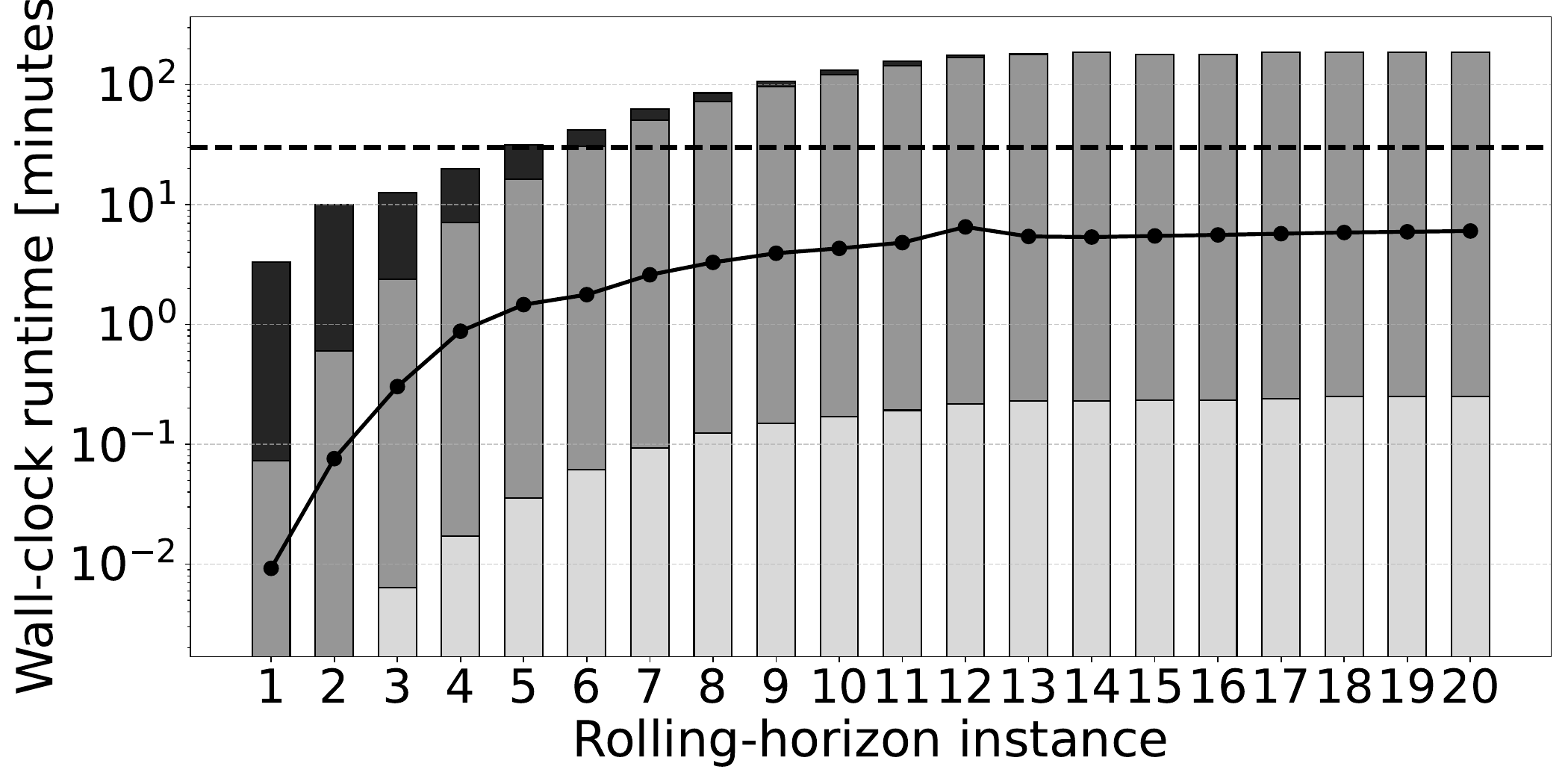}
            \caption{Scenario 15}
            \end{subfigure}
            \hfill
            \begin{subfigure}{0.48\textwidth}
            \includegraphics[width=\linewidth]{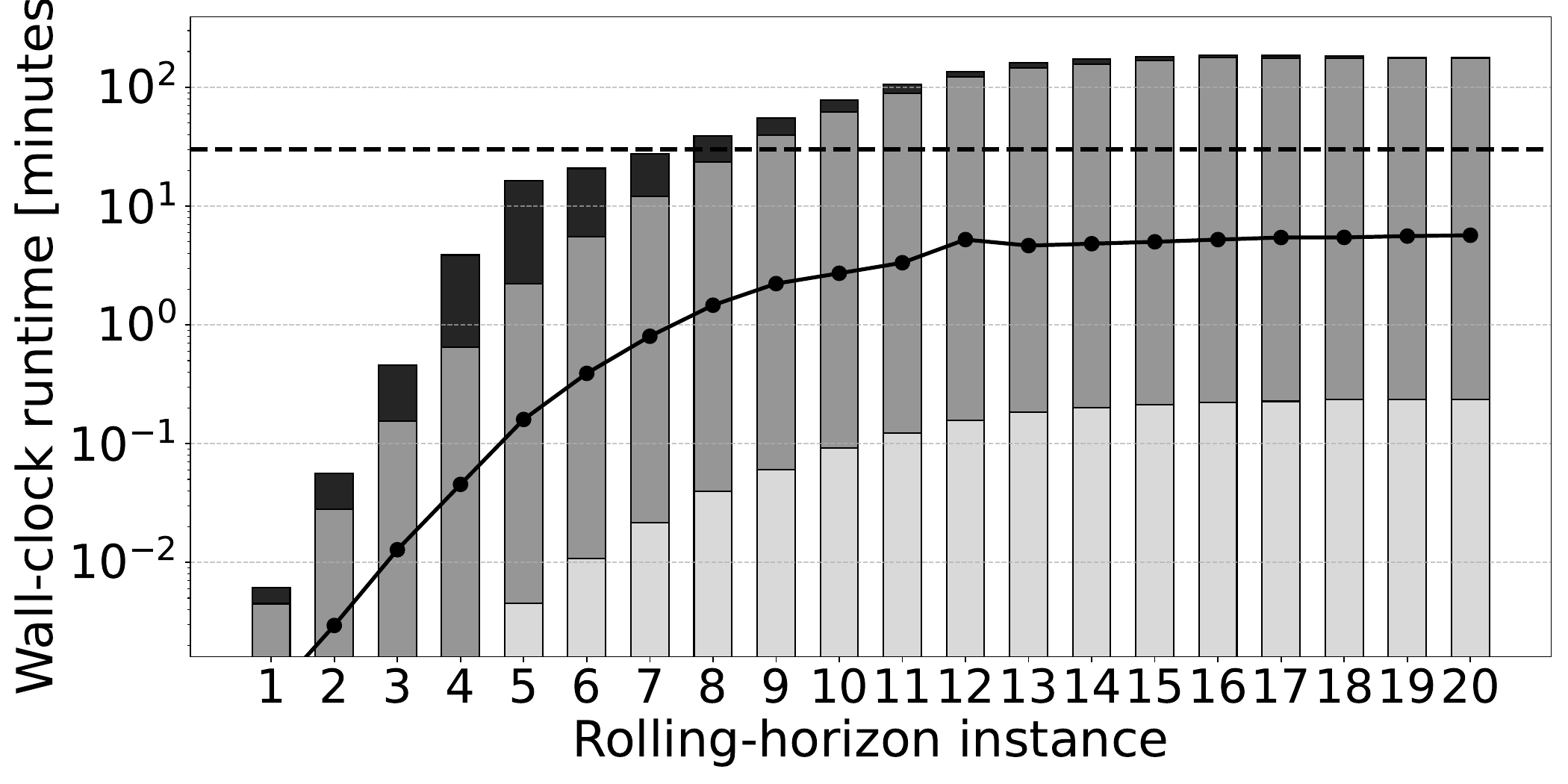}
            \caption{Scenario 16}
            \end{subfigure}

            \begin{subfigure}{0.7\textwidth}
            \includegraphics[width=\linewidth]{figures/runtime/legend.pdf}        
            \end{subfigure}
        
        \caption{Instance-level runtime breakdown for scenarios 13–16.
        Each plot reports the median runtime across the 10 stochastic replications
        for the 20 sequential instances.}
        \label{fig:runtime_instances_13_16}        
        \end{figure}

       This appendix reports detailed instance-level runtime breakdowns for all scenarios.
        For each scenario, the plots show the median runtime across the ten stochastic replications for the 20 sequential rolling-horizon instances (see Figures~\ref{fig:runtime_instances_1_4}--\ref{fig:runtime_instances_13_16}).
        The y-axis is displayed on a logarithmic scale in order to clearly visualize the large runtime differences between the heuristic and the exact solver.

    Summary runtime statistics are reported in Tables~\ref{tab:runtime_statistics} and \ref{tab:runtime_exceedance_scenarios}.
    
        \begin{table}[p]
        \centering
        \caption{Runtime statistics of the heuristic compared to the exact solver.}
        \label{tab:runtime_statistics}
            \begin{tabular}{lr}
            \toprule
            Statistic & Value \\
            \midrule
            Total instances & 3\;200 \\
            Solver time-limit violations (30 min) & 61.8\% \\
            Heuristic time-limit violations & 0.0\% \\
            Median runtime speedup & $27.9\times$ \\
            Mean runtime speedup & $44.9\times$ \\
            Speedup (min) & $2.6\times$ \\
            Speedup (max) & $1408.9\times$ \\
            \bottomrule
            \end{tabular}
        \end{table}

        \begin{table}[p]
        \centering
        \caption{Share of instances exceeding the 30-minute runtime limit per scenario.}
        \label{tab:runtime_exceedance_scenarios}
            \begin{tabular}{c r}
            \toprule
            Scenario & Time-limit violations \\
            \midrule
            1 & 52.5\% \\
            2 & 46.5\% \\
            3 & 52.0\% \\
            4 & 47.0\% \\
            5 & 71.0\% \\
            6 & 59.5\% \\
            7 & 72.0\% \\
            8 & 60.5\% \\
            9 & 65.0\% \\
            10 & 56.5\% \\
            11 & 66.5\% \\
            12 & 55.5\% \\
            13 & 77.0\% \\
            14 & 65.0\% \\
            15 & 77.5\% \\
            16 & 65.0\% \\
            \bottomrule
            \end{tabular}
        \end{table}

    \subsection{Additional runtime–quality trade-off results}
    \label{subsec:appendix_tradeoff}

\begin{figure}[htbp]
        \centering
        \includegraphics[page=2,width=0.75\textwidth]{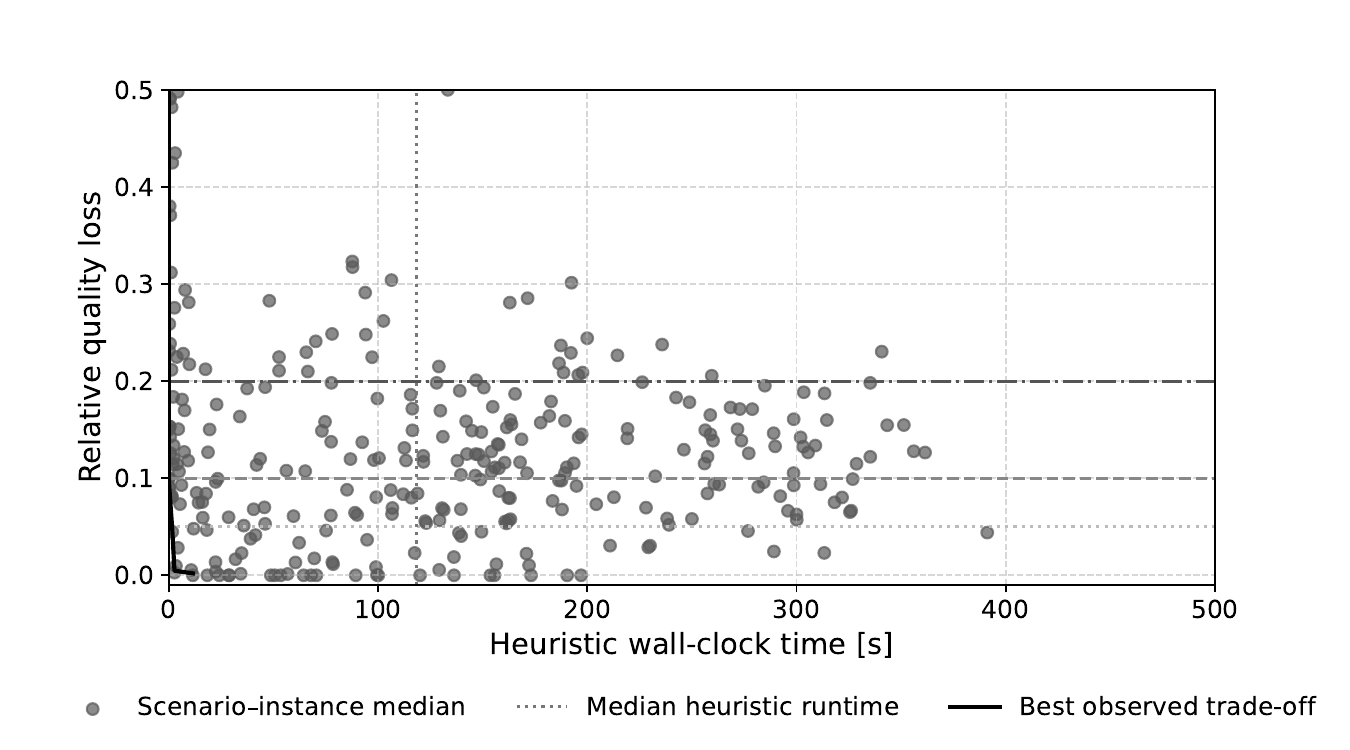}
        \caption{Runtime–quality trade-off of the heuristic for Objective~4.
        Each point represents the median performance across stochastic
        replications for a scenario–instance combination.}
        \label{fig:tradeoff_obj4}
        \end{figure}

    For completeness, Figure~\ref{fig:tradeoff_obj4} reports the runtime–quality trade-off for Objective~4, which represents workload balancing within the same priority class.

    The overall pattern is similar to the results observed for Objective~3 reported in Section~\ref{subsec:scalability}.
    In particular, the heuristic achieves low quality losses across most scenario–instance combinations while maintaining short computation times.
    Because the qualitative behavior closely resembles the results for Objective~3, the corresponding figure is reported here for completeness.
\section{Additional data for the computational study}
\label{sec:add_senario_instance_info}

    \subsection{Activity and Capability Categories}
    \label{subsec:activities_capabilities}
        The activities and capability categories originate from the KUBAS spontaneous volunteer coordination system.
        In this system, operational tasks are decomposed into activities that require specific volunteer capabilities, which are used for the assignment of volunteers to tasks \citep{Sperling2022}.
        Volunteers specify their offers of assistance based on location, availability, and a selection of generic capability types implemented in the system \citep{Betke2024HICSS}.
    
        Table~\ref{tab:activities} lists the activities together with their descriptions and the capability type required to perform them. 
            \begin{table}[h]
            \centering
            \caption{Activities and required capabilities}
            \label{tab:activities}
                \begin{tabularx}{\linewidth}{clcX}
                \toprule
                ID & Activity Description & Cap. ID & Capability Description \\
                \midrule
                1  & Carrying water containers & 1 & Heavy physical work \\
                2  & Filling water containers & 2 & Medium physical work \\
                3  & On-site documentation & 5 & Light documentation work \\
                4  & Cleaning and disinfection & 2 & Medium physical work \\
                5  & Information dissemination & 3 & Light physical work \\
                6  & Domestic assistance & 3 & Light physical work \\
                7  & Filling sandbags & 2 & Medium physical work \\
                8  & Carrying sandbags & 1 & Heavy physical work \\
                9  & Transport assignments & 6 & Caregiving tasks including vehicle operation \\
                10 & Local logistics & 3 & Light physical work \\
                11 & Dike inspection & 3 & Light physical work \\
                12 & Meal distribution & 4 & Caregiving tasks \\
                13 & Meal preparation & 2 & Medium physical work \\
                14 & Meal delivery & 2 & Medium physical work \\
                15 & Caregiving & 4 & Caregiving tasks \\
                \bottomrule
                \end{tabularx}
            \end{table}

    \subsection{Tasks from the 2013 Halle Flood}
    \label{subsec:tasks}
\begin{figure}[h]    
        \centering          
            \includegraphics[width=0.8\textwidth]{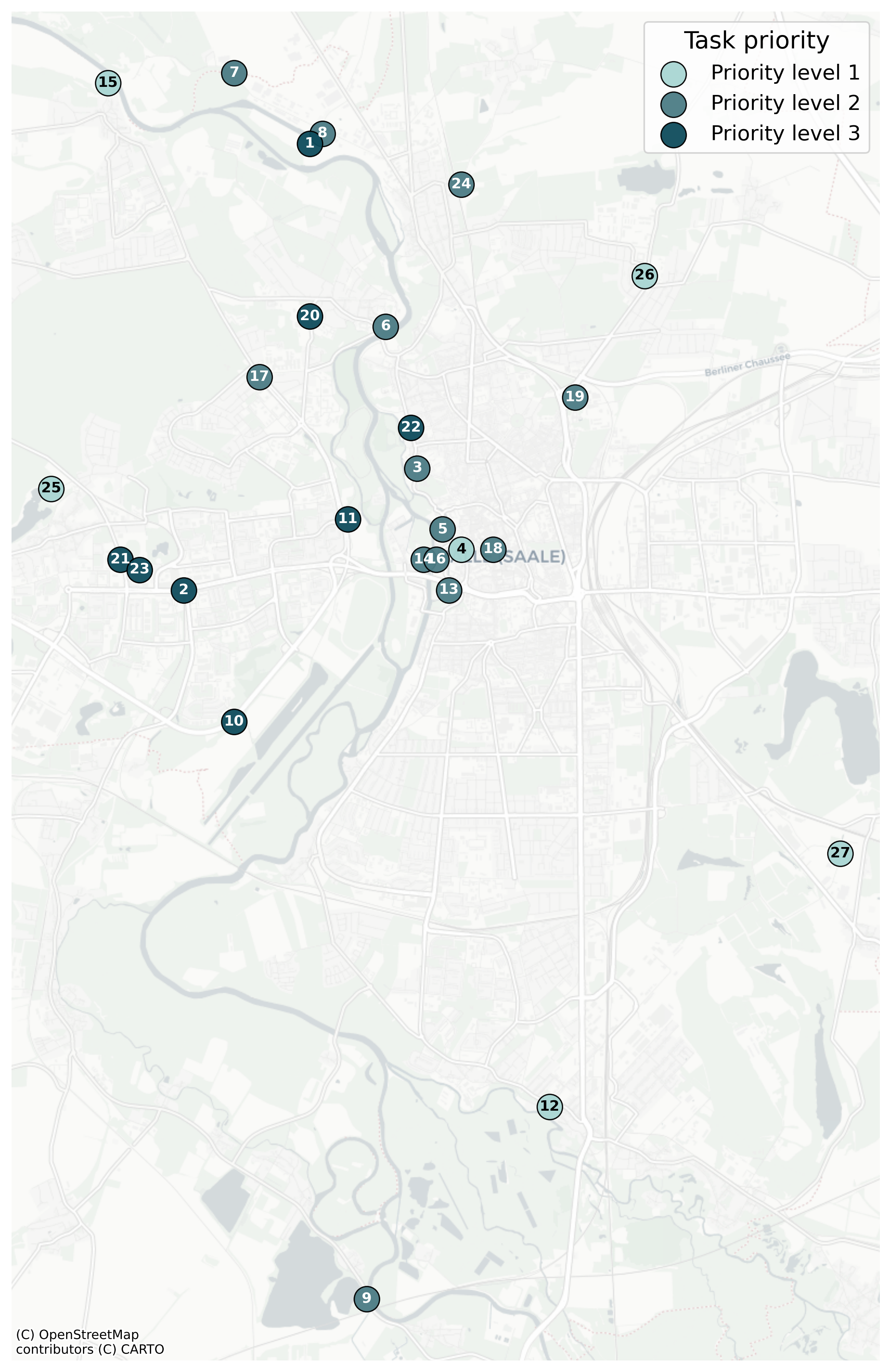}
        \caption{Spatial distribution of tasks and their priority levels. Colors indicate task priority (light = low, medium = medium, dark = high). 
        Numbers correspond to the task identifiers. Adapted from \citep{Sperling2022}.}
        \label{fig:task_map}
        \end{figure}

        The flood scenario 2013 in Halle, Germany, consists of a set of spatially distributed tasks that require different types of activities and varying numbers of volunteers. 
        Each task is associated with a geographic location and a priority level indicating its operational importance. 
        Figure~\ref{fig:task_map} illustrates the spatial distribution of all tasks within the considered operational area. 
        Tasks are represented as points on the map and are color-coded according to their priority level. 
        The numbers shown at each location correspond to the respective task identifiers.
    
        Each task may require multiple activities that must be performed by volunteers with suitable capabilities. 
        For every activity, a specific demand for volunteers is defined. 
        Table~\ref{tab:task_activities} summarizes the required activities for each task together with the corresponding number of volunteers needed to perform them.
    
        The scenario therefore captures both the spatial distribution of operational demands and the heterogeneous activity requirements associated with individual tasks. 
        This combination allows the evaluation of coordination strategies under realistic conditions involving spatial constraints and varying workload demands.
        
            \begin{longtable}{lll}   
            \caption{Overview of activities required for each task and the associated volunteer demand, adapted from \citep{Sperling2022}.}
            \label{tab:task_activities}\\
                \toprule
                Task id & Activity & Demand of Volunteers \\
                \midrule
                \endfirsthead
        
                \toprule
                Task id & Activity & Demand of Volunteers \\
                \midrule
                \endhead
        
                \midrule
                \multicolumn{3}{r}{Continued on next page} \\
                \midrule
                \endfoot
                
                \rowcolor{gray!10}
                1 & On-site documentation & 1 \\
                \rowcolor{gray!10}
                 & Carrying sandbags & 25 \\
                \rowcolor{gray!10}
                 & Meal distribution & 3 \\
                
                2 & On-site documentation & 2 \\
                 & Carrying sandbags & 96 \\
                 & Meal distribution & 10 \\
                
                \rowcolor{gray!10}
                3 & On-site documentation & 1 \\
                \rowcolor{gray!10}
                 & Carrying sandbags & 25 \\
                \rowcolor{gray!10}
                 & Meal distribution & 3 \\
                
                4 & On-site documentation & 1 \\
                 & Carrying sandbags & 25 \\
                 & Meal distribution & 3 \\
                
                \rowcolor{gray!10}
                5 & On-site documentation & 1 \\
                \rowcolor{gray!10}
                 & Carrying sandbags & 25 \\
                \rowcolor{gray!10}
                 & Meal distribution & 3 \\
                
                6 & On-site documentation & 1 \\
                 & Carrying sandbags & 25 \\
                 & Meal distribution & 3 \\
                
                \rowcolor{gray!10}
                7 & On-site documentation & 1 \\
                \rowcolor{gray!10}
                 & Carrying sandbags & 40 \\
                \rowcolor{gray!10}
                 & Meal distribution & 4 \\
                
                8 & On-site documentation & 1 \\
                 & Carrying sandbags & 25 \\
                 & Meal distribution & 3 \\
                
                \rowcolor{gray!10}
                9 & On-site documentation & 1 \\
                \rowcolor{gray!10}
                 & Carrying sandbags & 25 \\
                \rowcolor{gray!10}
                 & Meal distribution & 3 \\
                
                10 & On-site documentation & 1 \\
                 & Carrying sandbags & 30 \\
                 & Meal distribution & 3 \\
                
                \rowcolor{gray!10}
                11 & On-site documentation & 9 \\
                \rowcolor{gray!10}
                 & Carrying sandbags & 440 \\
                \rowcolor{gray!10}
                 & Meal distribution & 44 \\
                
                12 & On-site documentation & 1 \\
                 & Carrying sandbags & 25 \\
                 & Meal distribution & 3 \\
                
                \rowcolor{gray!10}
                13 & On-site documentation & 1 \\
                \rowcolor{gray!10}
                 & Carrying sandbags & 25 \\
                \rowcolor{gray!10}
                 & Meal distribution & 3 \\
                
                14 & On-site documentation & 1 \\
                 & Carrying sandbags & 25 \\
                 & Meal distribution & 3 \\
                
                \rowcolor{gray!10}
                15 & On-site documentation & 1 \\
                \rowcolor{gray!10}
                 & Carrying sandbags & 25 \\
                \rowcolor{gray!10}
                 & Meal distribution & 3 \\
                
                16 & On-site documentation & 1 \\
                 & Carrying sandbags & 25 \\
                 & Meal distribution & 3 \\
                
                \rowcolor{gray!10}
                17 & On-site documentation & 8 \\
                \rowcolor{gray!10}
                 & Filling sandbags & 90 \\
                \rowcolor{gray!10}
                 & Carrying sandbags & 270 \\
                \rowcolor{gray!10}
                 & Meal distribution & 36 \\
                
                18 & On-site documentation & 22 \\
                 & Filling sandbags & 270 \\
                 & Carrying sandbags & 810 \\
                 & Meal distribution & 108 \\
                
                \rowcolor{gray!10}
                19 & On-site documentation & 3 \\
                \rowcolor{gray!10}
                 & Filling sandbags & 30 \\
                \rowcolor{gray!10}
                 & Carrying sandbags & 90 \\
                \rowcolor{gray!10}
                 & Meal distribution & 12 \\
                
                20 & On-site documentation & 1 \\
                 & Information dissemination & 8 \\
                 & Transport missions & 8 \\
                 & Meal distribution & 8 \\
                 & Care support & 16 \\
                
                \rowcolor{gray!10}
                21 & On-site documentation & 1 \\
                \rowcolor{gray!10}
                 & Information dissemination & 8 \\
                \rowcolor{gray!10}
                 & Transport missions & 8 \\
                \rowcolor{gray!10}
                 & Meal distribution & 8 \\
                \rowcolor{gray!10}
                 & Care support & 16 \\
                
                22 & Filling water containers & 8 \\
                 & On-site documentation & 1 \\
                 & Information dissemination & 8 \\
                 & Meal distribution & 8 \\
                 & Care support & 16 \\
                
                \rowcolor{gray!10}
                23 & On-site documentation & 1 \\
                \rowcolor{gray!10}
                 & Meal distribution & 25 \\
                
                24 & On-site documentation & 1 \\
                 & Meal distribution & 25 \\
                
                \rowcolor{gray!10}
                25 & On-site documentation & 1 \\
                \rowcolor{gray!10}
                 & Meal distribution & 25 \\
                
                26 & On-site documentation & 1 \\
                 & Meal distribution & 25 \\
                
                \rowcolor{gray!10}
                27 & On-site documentation & 1 \\
                \rowcolor{gray!10}
                 & Meal distribution & 25 \\
                \bottomrule
            \end{longtable}

\end{document}